\documentclass{siamonline171218}

\usepackage{lipsum}
\usepackage{amsfonts}
\usepackage{graphicx}
\usepackage{epstopdf}
\usepackage{algorithmic}
\ifpdf
  \DeclareGraphicsExtensions{.eps,.pdf,.png,.jpg}
\else
  \DeclareGraphicsExtensions{.eps}
\fi

\numberwithin{theorem}{section}

\newcommand{\TheTitle}{Higher-Order Total Directional Variation: Imaging Applications} 
\newcommand{\TheAuthors}{S. Parisotto, J. Lellmann, S. Masnou and C.-B. Sch{\"{o}}nlieb}

\headers{\TheTitle}{\TheAuthors}

\title{
{\TheTitle}
\thanks{
Submitted to the editors DATE.
\funding{
SP acknowledges UK EPSRC grant EP/L016516/1 for the University of Cambridge, Cambridge Centre for Analysis DTC. 
CBS acknowledges support from
the EPSRC grants Nr. EP/M00483X/1, EP/K009745/1,
the EPSRC centre EP/N014588/1,
the Leverhulme Trust project 'Breaking the non-convexity barrier',
the Alan Turing Institute TU/B/000071, 
the CHiPS (Horizon 2020 RISE project grant), 
the Isaac Newton Institute 
and 
the Cantab Capital Institute for the Mathematics of Information.
SM acknowledges support from
the Labex MILYON/ANR-10-LABX-0070 and the ANR-14-CE27-0019 MIRIAM project grant
 }
}
}

\author{
Simone Parisotto\thanks{CCA, University of Cambridge, Wilberforce Road,
Cambridge CB3 0WA, UK (\email{sp751@cam.ac.uk})}
\and
Jan Lellmann\thanks{MIC, University of L\"{u}beck, Maria-Goeppert-Stra\ss e 3,
23562 L\"{u}beck, DE (\url{jan.lellmann@mic.uni-luebeck.de})}
\and
Simon Masnou\thanks{Universit\'{e} Claude Bernard Lyon 1, Institut Camille Jordan, Lyon, France (\email{masnou@math.univ-lyon1.fr})} 
\and 
Carola-Bibiane Sch\"onlieb\thanks{DAMTP, University of Cambridge, Wilberforce Road,
Cambridge CB3 0WA, UK (\url{cbs31@cam.ac.uk})}
}

\usepackage{amsopn}
\DeclareMathOperator{\diag}{diag}

\ifpdf
\hypersetup{
  pdftitle={\TheTitle},
  pdfauthor={\TheAuthors}
}
\fi

\usepackage[english]{babel}
\usepackage[utf8]{inputenc} 

\usepackage{tipa}
\usepackage{amsmath}
\usepackage{amsfonts}
\usepackage{amssymb}
\usepackage{mathtools}

\usepackage[optno]{optional}

\usepackage{epigraph}
\usepackage{pifont}

\usepackage[super]{nth}

\ifx\doctype\doctypepresentation
\else
\usepackage{footnote}
\makesavenoteenv{tabular}
\makesavenoteenv{table}
\fi

\usepackage{scalerel,stackengine}
\stackMath
\newcommand\widecheck[1]{\savestack{\tmpbox}{\stretchto{\scaleto{\scalerel*[\widthof{\ensuremath{#1}}]{\kern-.6pt\bigwedge\kern-.6pt}{\rule[-\textheight/2]{1ex}{\textheight}}}{\textheight}}{0.5ex}}\stackon[1pt]{#1}{\scalebox{-1}{\tmpbox}}}

\renewcommand\widehat[1]{\savestack{\tmpbox}{\stretchto{\scaleto{\scalerel*[\widthof{\ensuremath{#1}}]{\kern-.6pt\bigwedge\kern-.6pt}{\rule[-\textheight/2]{1ex}{\textheight}}}{\textheight}}{0.5ex}}\stackon[1pt]{#1}{\tmpbox}}

\AtBeginDocument{
  
  \let\argmin\relax
  \let\diag\relax

  \let\div\relax

  \let\grad\relax

  \let\prox\relax

  \let\tr\relax

  \DeclareMathOperator*{\argmin}{arg\,min}
  \DeclareMathOperator{\diag}{diag}

  \DeclareMathOperator{\div}{div}

  \DeclareMathOperator{\grad}{\bm{\nabla}}

  \DeclareMathOperator{\prox}{\mathbf{prox}}

  \DeclareMathOperator{\tr}{trace}
}
 \usepackage{calrsfs}
\DeclareMathAlphabet{\pazocal}{OMS}{zplm}{m}{n}

\newcommand{\Ccal}{\pazocal{C}}
\newcommand{\Dcal}{\pazocal{D}}

\newcommand{\Kcal}{\pazocal{K}}
\newcommand{\Ical}{\pazocal{I}}

\newcommand{\Mcal}{\pazocal{M}}

\newcommand{\Pcal}{\pazocal{P}}

\newcommand{\Rcal}{\pazocal{R}}
\newcommand{\Scal}{\pazocal{S}}
\newcommand{\Tcal}{\pazocal{T}}

\newcommand{\Ycal}{\pazocal{Y}}
\newcommand{\Wcal}{\pazocal{W}}

\newcommand{\Qrm}{\mathrm{Q}}

\newcommand{\Xrm}{\mathrm{X}}

\newcommand{\Abold}{\mathbf{A}}

\newcommand{\Jbold}{\mathbf{J}}

\newcommand{\Ibold}{\mathbf{I}}

\newcommand{\Lbold}{\mathbf{L}}
\newcommand{\Mbold}{\mathbf{M}}

\newcommand{\Qbold}{\mathbf{Q}}
\newcommand{\Rbold}{\mathbf{R}}
\newcommand{\Sbold}{\mathbf{S}}

\newcommand{\Ubold}{\mathbf{U}}

\newcommand{\Xbold}{\mathbf{X}}

\newcommand{\Wbold}{\mathbf{W}}
\newcommand{\Zbold}{\mathbf{Z}}

\newcommand{\abold}{{\bm{a}}}
\newcommand{\bbold}{{\bm{b}}}

\newcommand{\ebold}{{\bm{e}}}

\newcommand{\kbold}{{\bm{k}}}

\newcommand{\pbold}{{\bm{p}}}

\newcommand{\rbold}{{\bm{r}}}
\newcommand{\sbold}{{\bm{s}}}
\newcommand{\tbold}{{\bm{t}}}
\newcommand{\ubold}{{\bm{u}}}
\newcommand{\vbold}{{\bm{v}}}
\newcommand{\xbold}{{\bm{x}}}
\newcommand{\ybold}{{\bm{y}}}
\newcommand{\wbold}{{\bm{w}}}
\newcommand{\zbold}{{\bm{z}}}

\usepackage{bm} \newcommand{\alphabold}{{\bm{\alpha}}}
\newcommand{\etabold}{{\bm{\eta}}}

\newcommand{\thetabold}{{\bm{\theta}}}

\newcommand{\Sigmabold}{{\bm{\Sigma}}}
\newcommand{\xibold}{{\bm{\xi}}}
\newcommand{\iotabold}{{\bm{\iota}}}

\newcommand{\Psibold}{{\bm{\Psi}}}

\newcommand{\Lambdabold}{{\bm{\Lambda}}}

 \ifx\doctype\doctypepresentation
\usepackage{color}
\else
\usepackage[table]{xcolor}
\fi

\definecolor{darkred}{rgb}{0.55,0.0,0.0}
\definecolor{darkgreen}{rgb}{0,0.55,0.0}
\definecolor{darkblue}{rgb}{0,0.0,0.55}
\definecolor{RoyalRed}{rgb}{0.6179,0.0236,0.0894} 
\definecolor{deepblue}{rgb}{0,0,0.5}
\definecolor{deepred}{rgb}{0.6,0,0}
\definecolor{deepgreen}{rgb}{0,0.5,0}
\definecolor{grey}{gray}{0.5}
\definecolor{lightgrey}{gray}{0.8}
\definecolor{lightblue}{rgb}{0.4,0.6,0.8}
\definecolor{lightred}{rgb}{0.8,0.6,0.4}

\colorlet{commentcolour}{green!50!black}
\colorlet{stringcolour}{red!60!black}
\colorlet{keywordcolour}{magenta!90!black}
\colorlet{exceptioncolour}{yellow!50!red}
\colorlet{commandcolour}{blue!60!black}
\colorlet{numpycolour}{blue!60!green}
\colorlet{literatecolour}{magenta!90!black}
\colorlet{promptcolour}{green!50!black}
\colorlet{specmethodcolour}{violet}
\colorlet{indendifiercolour}{green!70!white} 
\usepackage{xparse}
\DeclareDocumentCommand\CCCspace{mo}
 {\IfValueTF{#2}
   {\mathrm{{C}}^{#1}(#2)}
   {\mathrm{{C}}^{#1}}}
 \DeclareDocumentCommand\CCospace{mo}
 {\IfValueTF{#2}
   {\mathrm{{C}}_0^{#1}(#2)}
   {\mathrm{{C}}_0^{#1}}}
 
 \DeclareDocumentCommand\CCCcomp{mo}
 {\IfValueTF{#2}
   {\mathrm{{C}}_c^{#1}(#2)}
   {\mathrm{{C}}_c^{#1}}}

\DeclareDocumentCommand\WW{mo}
 {\IfValueTF{#2}
   {\mathrm{{W}}^{#1}(#2)}
   {\mathrm{{W}}^{#1}}}
 
\DeclareDocumentCommand\WWloc{mo}
 {\IfValueTF{#2}
   {\mathrm{{W}}_{\text{loc}}^{#1}(#2)}
   {\mathrm{{W}}_{\text{loc}}^{#1}}}
 
\DeclareDocumentCommand\WWo{mo}
 {\IfValueTF{#2}
   {\mathrm{{W}}_0^{#1}(#2)}
   {\mathrm{{W}}_0^{#1}}}
 
\DeclareDocumentCommand\HH{mo}
 {\IfValueTF{#2}
   {\mathrm{{H}}^{#1}(#2)}
   {\mathrm{{H}}^{#1}}}
 
\DeclareDocumentCommand\HHo{mo}
 {\IfValueTF{#2}
   {\mathrm{{H}}_0^{#1}(#2)}
   {\mathrm{{H}}_0^{#1}}}

\DeclareDocumentCommand\LL{mo}
 {\IfValueTF{#2}
   {\mathrm{{L}}^{#1}(#2)}
   {\mathrm{{L}}^{#1}}}
 
\DeclareDocumentCommand\LLc{mo}
 {\IfValueTF{#2}
   {\mathrm{{L}}_\mathrm{c}^{#1}(#2)}
   {\mathrm{{L}}_\mathrm{c}^{#1}}}
 
\DeclareDocumentCommand\LLloc{mo}
 {\IfValueTF{#2}
   {\mathrm{{L}}_{\text{loc}}^{#1}(#2)}
   {\mathrm{{L}}_{\text{loc}}^{#1}}} 
 
\DeclareDocumentCommand\BH{mo}
 {\IfValueTF{#2}
   {\mathrm{{BH}}^{#1}(#2)}
   {\mathrm{{BH}}^{#1}}}
 
\DeclareDocumentCommand\BV{mo}
 {\IfValueTF{#2}
   {\mathrm{{BV}}^{#1}(#2)}
   {\mathrm{{BV}}^{#1}}}
 
\DeclareDocumentCommand\BDV{mmo}
 {\IfValueTF{#3}
   {\mathrm{{BDV}}^{(#1,#2)}(#3)}
   {\mathrm{{BDV}}^{(#1,#2)}}}
 
 \DeclareDocumentCommand\BDVM{mmo}
 {\IfValueTF{#3}
   {\mathrm{{BDV}}^{#1}_{#2}(#3)}
   {\mathrm{{BDV}}^{#1}_{#2}}}
 
\DeclareDocumentCommand\BDVb{mmmo}
 {\IfValueTF{#4}
   {\mathrm{{BDV}}_{#3}^{(#1,#2)}(#4)}
   {\mathrm{{BDV}}_{#3}^{(#1,#2)}}}
 
\DeclareDocumentCommand\BDVloc{mmo}
 {\IfValueTF{#3}
   {\mathrm{{BDV}}_{\text{loc}}^{(#1,#2)}(#3)}
   {\mathrm{{BDV}}_{\text{loc}}^{(#1,#2)}}} 
 
\DeclareDocumentCommand\BGDV{mmo}
 {\IfValueTF{#3}
   {\mathrm{{BGDV}}^{(#1,#2)}(#3)}
   {\mathrm{{BGDV}}^{(#1,#2)}}}
 
\DeclareDocumentCommand\BGV{mo}
 {\IfValueTF{#2}
   {\mathrm{{BGV}}^{#1}(#2)}
   {\mathrm{{BGV}}^{#1}}}
 
\DeclareDocumentCommand\BD{mo}
 {\IfValueTF{#2}
   {\mathrm{{BD}}^{#1}(#2)}
   {\mathrm{{BD}}^{#1}}} 
 
\DeclareDocumentCommand\BVloc{mo}
 {\IfValueTF{#2}
   {\mathrm{{BV}}_{\text{loc}}^{#1}(#2)}
   {\mathrm{{BV}}_{\text{loc}}^{#1}}}
\DeclareDocumentCommand\SBV{mo}
 {\IfValueTF{#2}
   {\mathrm{{SBV}}^{#1}(#2)}
   {\mathrm{{SBV}}^{#1}}}
 
\DeclareDocumentCommand\TV{mo}
 {\IfValueTF{#2}
   {\mathrm{{TV}}^{#1}(#2)}
   {\mathrm{{TV}}^{#1}}}
 
  \DeclareDocumentCommand\Sym{mo}
 {\IfValueTF{#2}
   {\mathrm{{Sym}}^{#1}(#2)}
   {\mathrm{{Sym}}^{#1}}}
\DeclareDocumentCommand\TGV{mmo}
 {\IfValueTF{#3}
   {\mathrm{{TGV}}_{#1}^{#2}(#3)}
   {\mathrm{{TGV}}_{#1}^{#2}}}
\DeclareDocumentCommand\TGDV{mmmo}
 {\IfValueTF{#4}
   {\mathrm{{TGDV}}^{#1}_{#2, #3}(#4)}
   {\mathrm{{TGDV}}^{#1}_{#2, #3}}} 
 
 \DeclareDocumentCommand\PTV{mo}
 {\IfValueTF{#2}
   {\mathrm{{PTV}}^{#1}(#2)}
   {\mathrm{{PTV}}^{#1}}}

 \DeclareDocumentCommand\TDV{mmo}
 {\IfValueTF{#3}
   {\mathrm{{TDV}}^{(#1,#2)}({#3})}
   {\mathrm{{TDV}}^{(#1,#2)}}}
 
  \DeclareDocumentCommand\TDVb{mmmo}
 {\IfValueTF{#4}
   {\mathrm{{TDV}}_{{#3}}^{(#1,#2)}({#4})}
   {\mathrm{{TDV}}_{{#3}}^{(#1,#2)}}}
 
   \DeclareDocumentCommand\TDVM{mmo}
 {\IfValueTF{#3}
   {\mathrm{{TDV}}_{{#2}}^{#1}({#3})}
   {\mathrm{{TDV}}_{{#2}}^{#1}}}
 
   \DeclareDocumentCommand\DVM{mmo}
 {\IfValueTF{#3}
   {\mathrm{{DV}}_{{#2}}^{#1}({#3})}
   {\mathrm{{DV}}_{{#2}}^{#1}}} 
 
  \DeclareDocumentCommand\TDVh{mmo}
 {\IfValueTF{#3}
   {\mathrm{{TDV}}_h^{(#1,#2)}(#3)}
   {\mathrm{{TDV}}_h^{(#1,#2)}}}
 
\DeclareDocumentCommand\DTV{mo}
 {\IfValueTF{#2}
   {\mathrm{{DTV}}^{#1}(#2)}
   {\mathrm{{DTV}}^{#1}}} 
 
 \DeclareDocumentCommand\DTGV{mmo}
 {\IfValueTF{#3}
   {\mathrm{{DTGV}}^{#1}_{#2}(#3)}
   {\mathrm{{DTGV}}^{#1}_{#2}}} 
 
 \DeclareDocumentCommand\EXP{mo}
 {\IfValueTF{#2}
   {\EE_{#2}\left[#1\right]}
   {\EE\left[#1\right]}}

 \usepackage[algo2e]{algorithm2e}

\usepackage[formats]{listings}

\newcommand{\codetitlestyle}[1]{\small\textit{#1}}
\newcommand{\belowtitleskip}{2pt}

\usepackage{listingsutf8}
\usepackage[framed,numbered,autolinebreaks,useliterate]{mcode}

\DeclareFixedFont{\ttb}{T1}{txtt}{bx}{n}{12} \DeclareFixedFont{\ttm}{T1}{txtt}{m}{n}{12}  

\lstdefinestyle{pythonstyle}{
language=Python,
basicstyle=\ttfamily\small,
frame=tb,                         }

\lstdefinestyle{matlabstyle}{
language=matlab,
basicstyle=\ttfamily\small,
frame=tb,
rulecolor=\color{black!40},
emphstyle=\color{blue},
commentstyle=\color{commentcolour}\slshape,
}

\lstdefinestyle{matlab}{frame=single,                           basicstyle=\scriptsize\ttfamily,             keywordstyle=[1]\color{darkblue}\bfseries,        keywordstyle=[2]\color{darkred},         keywordstyle=[3]\color{darkred}\underbar,  identifierstyle=,                       commentstyle=\usefont{T1}{pcr}{m}{sl}\color{darkred}\small,
        stringstyle=\color{darkred},             showstringspaces=false,                 tabsize=5,                              morekeywords={xlim,ylim,var,alpha,factorial,poissrnd,normpdf,normcdf},
morekeywords=[2]{on, off, interp},
morekeywords=[3]{FindESS, homework_example},
morecomment=[l][\color{darkred}]{...},     numbers=left,                           firstnumber=1,                          numberstyle=\tiny\color{darkred},          stepnumber=5                            }

\lstnewenvironment{python}[1][]
{
\pythonstyle
\lstset{#1}
}
{}

\newcommand\pythoninline[1]{{\pythonstyle\lstinline!#1!}}

\lstnewenvironment{matlab}[1][]{\lstset{style=matlabstyle,title={\codetitlestyle{MATLAB code }},belowcaptionskip=\belowtitleskip}}{}

\newcommand{\RR}{\mathbb{R}}

\newcommand{\ZZ}{\mathbb{Z}}
\newcommand{\NN}{\mathbb{N}}

\newcommand{\EE}{\mathbb{E}}

\newcommand{\diff}{\mathop{}\mathrm{d}}

\newcommand{\abs}[1]{\left\lvert #1 \right\rvert}
\newcommand{\norm}[1]{\left\lVert #1 \right\rVert}

\newcommand{\blank}{\,{\cdot}\,}

   \newcommand{\T}{\mathrm{T}} \newcommand{\SO}[1]{\mathrm{SO(#1)}}

\newlength{\imagewidth} 

\usepackage[optno]{optional}

\usepackage[binary-units=true]{siunitx}

\usepackage{mathtools}

\usepackage{tikz}
\usetikzlibrary{arrows,matrix,positioning,scopes,calc}

\usepackage{graphicx}
\usepackage[font=scriptsize]{subcaption}
\graphicspath{{./../images/}}
\usepackage{tabularx,longtable,booktabs,wrapfig,multirow}

\usepackage{epstopdf}
\AppendGraphicsExtensions{.tif}

\AppendGraphicsExtensions{.tif}

\newsiamremark{example}{Example}
\newsiamremark{remark}{Remark}

\allowdisplaybreaks

\crefformat{figure}{#2Figure~#1#3}
\crefformat{algorithm}{#2Algorithm~#1#3}
\crefformat{section}{#2Section~#1#3}
\crefformat{subsection}{#2Section~#1#3}
\crefformat{chapter}{#2Chapter~#1#3}
\crefformat{equation}{#2(#1)#3}
\crefformat{table}{#2Table~#1#3}
\crefformat{theorem}{#2Theorem~#1#3}
\crefformat{lemma}{#2Lemma~#1#3}
\crefformat{remark}{#2Remark~#1#3}
\crefformat{example}{#2Example~#1#3}
\crefformat{proposition}{#2Proposition~#1#3}
\crefformat{property}{#2Propery~#1#3}
\crefformat{definition}{#2Definition~#1#3}
\crefrangeformat{figure}{Figures~(#3#1#4--#5#2#6)}
\crefrangeformat{equation}{(#3#1#4)--(#5#2#6)} 
\usepackage{anyfontsize}
\allowdisplaybreaks

\begin{document}

\maketitle

\begin{abstract}
We introduce a class of higher-order anisotropic total variation regularisers, which are defined for possibly inhomogeneous, smooth elliptic anisotropies, that extends the Total Generalized Variation (TGV) regulariser and its variants. 
We propose a primal-dual hybrid gradient approach to approximate numerically the associated gradient flow. This choice of regularisers allows to preserve and enhance intrinsic anisotropic features in images. This is illustrated on various examples from different imaging applications: image denoising, wavelet-based image zooming, and reconstruction of surfaces from scattered height measurements.
\end{abstract}



\section{Introduction}
In the last decades, total variation ($\TV{}$) regularisation has been successfully applied to a variety of imaging problems. In particular since \cite{ROF},  $\TV{}$ plays a crucial role for variational image denoising, deblurring, inpainting, segmentation, magnetic resonance image (MRI) reconstruction and many others, see \cite{ChCaCrNoPo10}. While the $\TV{}$ regulariser successfully eliminates noise and at the same time preserves characteristic image features like edges, it still has some shortcomings. A major one is the staircasing effect, resulting into blocky-like images \cite{CasChaNov2007,Nikolova2000}. One approach to mitigate this effect is based on higher order total variation regularisers, see e.g.~\cite{ChaEsePar2010,ChanMarquinaMulet2000,Papafitsoros2014, SetzerSteidl,WuTai}, aiming to eliminate the staircasing effect by higher regularity in homogeneous regions of the image while still allowing for discontinuities in the presence of edges.  The total generalized variation ($\TGV{\alphabold}{q}$) regulariser has been proposed in \cite{BreKunPoc2010} to balance the first $q$ derivatives of $u$ with a regularisation parameter vector $\alphabold$.
Another modification of  the $\TV{}$ regulariser has been the introduction of directional information in the regularisation, allowing to smooth images in an anisotropic fashion favouring preferred directions, e.g.\ \cite{bayram2012directional, EADTV, BeBuDr06, Dong2009, Steidl2009, GL10, lenzen2013class, STV, EstSoaBre15, Ehrhardt16}. A recent combination of directional $\TV{}$ and higher-order derivatives is the directional total generalized variation \cite{directionaltv} that equips the $\TGV{\alphabold}{q}$ regulariser with one constant preferred smoothing direction.

\begin{figure}[tbp]
\centering
\begin{subfigure}[t]{0.325\textwidth}
\centering\captionsetup{justification=centering}
\includegraphics[width=\textwidth]{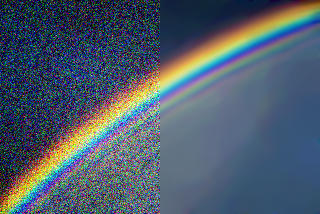}
\caption{Denoising (20\% Gaussian noise).}
\end{subfigure}
\hfill
\begin{subfigure}[t]{0.325\textwidth}
\centering\captionsetup{justification=centering}
\includegraphics[width=0.665\textwidth]{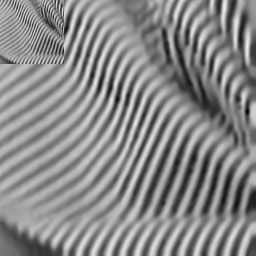}
\caption{4X Wavelet-based zooming.}
\end{subfigure}
\hfill
\begin{subfigure}[t]{0.325\textwidth}
\centering\captionsetup{justification=centering}
\includegraphics[width=0.82\textwidth]{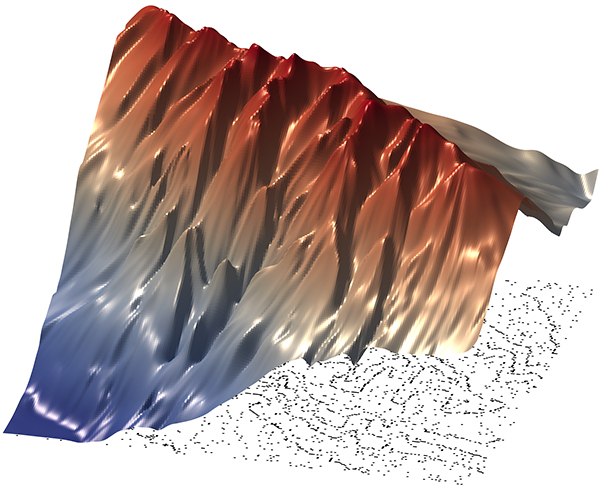}
\caption{DEM interpolation of 7\% data.}
\end{subfigure}
\caption{Imaging applications for the class of directional regularisers introduced in this paper.
}
\label{fig: application preview intro}
\end{figure}

In this paper we extend the directional total generalized variation introduced in \cite{directionaltv} and the third-order directional total variation regulariser introduced in \cite{LelMorSch2013} to a new class of directional total variation regularisers that can feature a combination of orders of derivatives as well as spatially-varying directional information by means of weighting the derivatives in the $\TGV{\alphabold}{q}$ regulariser with $2$-tensors. 
We introduce this class of total generalized variation regularisers, discuss its numerical solution by a tailored primal-dual hybrid gradient scheme (which replaces the commercial CVX+Mosek solver used in \cite{LelMorSch2013}), and showcase its performance and regularisation properties for a range of imaging problems. For the latter, we show the effect of this generalized class of regularisers in different imaging applications where the introduced anisotropy plays a crucial role: image denoising, wavelet-based zooming and digital elevation map (DEM) interpolation with applications to atomic force microscopy (AFM) data, see \cref{fig: application preview intro}. For the application of our regularisers to video denoising, we refer to \cite{ParSch2019} and in general to \cite{Parisotto2019}. The theoretical foundation for this general class of directional, higher-order regularisers is presented in our companion work \cite{ParMasSch18analysis}.

Let us go into more details. Let $\Omega\subset\RR^d$ be a bounded Lipschitz domain for $d\geq 1$ and $u:\Omega\to\RR$ a function, we define the the higher-order directional total variation of $u$ as
\begin{equation}
\TDVM{q}{\alphabold}[u,\Mcal]
:= 
\sup_{\Psibold}\left\{ \int_\Omega u \div_{\Mcal}^q \Psibold \diff \xbold \,\Big\lvert\,\text{ for all } \Psibold \in\Ycal_{\Mcal,\alphabold}^q \right\},
\label{eq: TDV intro def}
\end{equation}
where we call $q$ the order of the regularisation, $\Mcal$ is a collection of weighting fields, and
\begin{equation}
\Ycal_{\Mcal,\alphabold}^q 
= 
\left\{\Psibold\,:\, \Psibold \in \CCCcomp{q}[\Omega,\Tcal^q(\RR^d)],\, \norm{\div_{\Mcal}^j\Psibold}_\infty\leq \alpha_j,\forall\,j=0,\dots q-1\right\},
\end{equation}
where $\Tcal^q(\RR^d)$ is the vector space of $q$-tensors in $\RR^d$ and $\alphabold$ is a vector of regularisation parameters. We will provide the rigorous definition for \cref{eq: TDV intro def} in \cref{sec: higher-order tdv}. 
We comment for now that the regulariser in \cref{eq: TDV intro def} is designed for introducing weighted directional derivatives in the classical definition of $\TGV{\alphabold}{q}$. The anisotropy is introduced by a family of weights $\Mcal$ and a thereby suitably weighted divergence $\div_\Mcal^q$ of order $q$, defined in \cref{eq: weighted divergence}.

\subsection{Related work}
In what follows we review the state-of-the-art that is most relevant for the proposed higher-order directional total variation regulariser. 
We focus in particular on functional regularisers but it is worth mentioning that there is a rich literature on fairly general anisotropic PDE models, mainly of first and second order, see for instance \cite{SapRin1996,TshDer2005,WeickertBook,KimMalSoc2000,BloCha1998,CasTanSap1999,AlvGuiLioMor1993,DIZENZO1986116} and the references therein. Our model handles a more limited class of anisotropies but it can do it at any order of derivatives, particularly useful in various applications.

The idea of anisotropic smoothing for imaging has certainly been popularized by the book of Weickert on anisotropic diffusion equations \cite{WeickertBook} based on the key notion of structure tensor to encode directional information, see also~\cite{HarrisStephens88,Forstner86,KassWitkin85}. Weickert's structure tensor for a continuous imaging function $u:\Omega \rightarrow\mathbb R$ and non-negative parameters $\sigma,\rho$ is defined as
\begin{equation}
\Jbold_\rho(u) :=  K_\rho \ast \left(\grad u_\sigma \otimes \grad u_\sigma\right),
\label{eq: structure tensor weickert}
\end{equation}
where $u_\sigma = K_\sigma\ast u$ and $K_\sigma, K_\rho$ are Gaussian kernels with standard deviations $\sigma, \rho$, respectively. 
For  $\grad u_\sigma \neq 0$
the structure tensor $\Jbold_\rho(u)$ has two orthogonal eigenvectors $\vbold_1$ and $\vbold_2$ with corresponding non-zero real eigenvalues $\lambda_1(\xbold)$ and $\lambda_2(\xbold)$. 
Here $\vbold_1$ and $\vbold_2$ approximately point in the direction $\grad u_\sigma$ and $\grad^\perp u_\sigma$. From this, diffusion tensors can be constructed inheriting $\vbold_1$ and $\vbold_2$ as eigenvectors but whose eigenvalues are expressions of $\lambda_1$ and $\lambda_2$ so as to increase or reduce smoothing in these directions, compare for instance coherence-enhancing diffusion \cite{weickert1999coherence}.
The concept of  structure tensor is used for variational regularisation in \cite{Steidl2009} in the framework of a \emph{single orientation estimation} approach.  More precisely, the authors consider a regulariser of the type
\[
\int_\Omega w(\xbold) \Jbold_\rho(u)\diff\xbold,
\]
for a non-negative weight function $w:\Omega\to\RR$ and a continuous imaging function $u:\Omega\to\RR$, for smoothing an image into a dominant single direction. 
For smoothing a noisy data $u^\diamond$ in two directions, the authors propose to estimate directions $\vbold_1$ and  $\vbold_2$
as in \cite{aach2006analysis} and, in their \emph{double orientation estimation} approach, decompose $u(\xbold)=u_1(\xbold)+u_2(\xbold)$ via
\begin{equation}
\min_{\begin{subarray}{c}
u_1,u_2\\
u_1+u_1=u\end{subarray}} \alpha\left( \norm{\vbold_1^\T \grad u_1}_1 + \norm{\vbold_2^\T \grad u_2}_1\right)+ \frac{1}{2}\norm{u - u^\diamond}_2^2.
\end{equation}

An approach based on the analysis of eigenvalues and eigenvectors of the structure tensors can be found in \cite{GL10} while in \cite{lenzen2013class} the admissible set of test functions are locally adapted to the geometry of $\ubold$ via the support function regulariser.
Furthermore, in \cite{STV}, the \emph{structure tensor total variation} (STV) focuses on the nuclear norm of the structure tensor \cref{eq: structure tensor weickert} in order to measure the local image variation:
\begin{equation}
\mathrm{STV}_p(u) = \norm{
\left(
\sqrt{\lambda_1},\,
\sqrt{\lambda_2}\right)}_p.
\label{eq: STV}
\end{equation}
Also in \cite{EstSoaBre15}, a regulariser is proposed whose smoothing directions vary according to the image content, leading to the analysis of
\[
\Rcal_{\Jbold}(u) = g_\gamma(\Sbold \grad u),
\quad 
\Sbold = \max\left(
\sqrt{\gamma}, \sqrt[4]{\lambda_1 + \lambda_2} \Lambdabold^{-0.5} \Qbold^\T\right),
\]
where $g_\gamma$ is the Huber regularisation with parameter $\gamma>0$ and the structure tensor $\Jbold_\rho$ is eigen-decomposed as $\Jbold_\rho = \Qbold\Lambdabold\Qbold^\T$, with eigenvectors stored in the matrix $\Qbold$ and the eigenvalues in the diagonal matrix $\Lambdabold$. 

Let us also mention that early works where the gradient is weighted date back to \cite{BeBuDr06}, where the oriented local image structure is extracted from images by the regulariser
\begin{equation}
\TV{\alpha}[u] 
= 
\int_\Omega \abs{\Mbold_\alpha\grad u}\diff \xbold,
\label{eq: cartoonTV}
\end{equation}
with $\Mbold_\alpha$ being the orthogonal rotation matrix for an angle $\alpha >0$.

Further, in \cite{bayram2012directional}, a discrete \emph{directional total variation} ($\TV{}_{a,\theta}$) regulariser for denoising discrete images $\ubold$ with a single dominant direction (directional images) is introduced via affine transformations of test functions: the circular unit ball generated by the $\LL{2}$-norm is transformed into an ellipse $E^{a,\theta}$, with major semiaxis $a>1$ rotated by $\theta$, penalizing variations for large $a$ along $\theta$:
\begin{equation}
\TV{}_{a,\theta}(\ubold) = \sum_{i,j}  \sup_{\Psibold\in E^{a,\theta}}\langle (\grad \ubold)_{i,j},\,\Psibold\rangle.
\label{eq: TValphatheta}
\end{equation}

In a straightforward generalization of \cref{eq: TValphatheta}, inhomogeneous fields $\theta$ are allowed, namely $\thetabold:=\theta(\xbold)$.
In \cite{EADTV} the authors propose to adapt $\thetabold$ to the edge directions,
\begin{equation}
(\cos\thetabold,\,\sin\thetabold)= \frac{\left(\grad u_\sigma\right)^\perp}{\norm{\grad u_\sigma}_2},
\label{eq: cossinEADTV}
\end{equation}
so as to associate at each pixel position $(i,j)$ a specific ellipsoid ball $E^{a,\theta_{ij}}$ for the test functions, leading to the discrete \emph{edge adaptive directional total variation} (EADTV) regulariser:
\begin{equation}
\mathrm{EADTV}_{\alpha,\thetabold}(\ubold) 
= 
\sum_{i,j}  \sup_{\Psibold\in E^{\alpha,\thetabold_{ij}}}\langle (\grad \ubold)_{i,j},\,\Psibold\rangle.
\label{eq: EADTV}
\end{equation}

In \cite{Ehrhardt16}, a discrete weighted \emph{directional Total Variation} (dTV) regulariser is introduced as
\begin{equation}
\mathrm{dTV}(\ubold) 
= 
\sum_{n=1}^N\abs{\Pcal_{\xibold_n}\grad u_n},\quad\text{with}\quad \xibold_n:= \grad u_n / (K_\sigma\ast\abs{\grad u_n}),
\label{eq: dTV}
\end{equation}
by projecting onto the complementary part of a vector field $\Pcal_{\xibold_n} \xbold = \xbold - \langle \xibold_n,\,\xbold \rangle \xibold_n$. 

In \cite{directionaltv}, the continuous \emph{directional total variation} ($\DTV{}$) and \emph{directional total generalized variation} ($\DTGV{}{}$) are analysed for a single homogeneously fixed angle $\theta$ and for the minor semi-axis (now denoted with $a$) of the ellipse $E^{a,\theta}$. There, $\DTV{}$ and $\DTGV{}{}$ are built upon test functions $\Psibold$ in the isotropic ball $B_1$, so as to constrain $\widetilde{\Psibold}$ to the ellipsoid ball, i.e.\ $\widetilde{\Psibold}=\Rbold_\theta \Lambdabold_a \Psibold \in E^{a,\theta}$ and where $\Rbold_\theta$ and $\Lambdabold_a$ are rotation and contraction matrices, respectively, similarly to our setting explained later, see Equation \cref{eq: controrot matrices}:
\begin{align}
\DTV{}[u] 
&= 
\sup_{\widetilde{\Psibold}} \left\{\int_\Omega u\div\widetilde{\Psibold}\diff \xbold\,\Big\lvert\, \widetilde{\Psibold}\in\CCCcomp{1}[\Omega,\RR^2]\right\},
\label{eq: DTV}
\\
\DTGV{q}{\alphabold}[u] 
&=
\sup_{\widetilde{\Psibold}}\left\{ \int_\Omega u\div^q\widetilde{\Psibold}\diff \xbold\,\left|\, 
\begin{aligned}
&\widetilde{\Psibold}\in \CCCcomp{q}[\Omega,\Sym{q}(\RR^2)],\,\norm{\div^j{\widetilde{\Psibold}(\xbold)}}\leq \alpha_j,\\
&\forall\xbold\in\Omega\text{ and }\forall j= 0,\dots q-1
\end{aligned}
\right.\right\}.
\label{eq: DTGV}
\end{align}
By comparing \cref{eq: DTGV} with \cref{eq: TDV intro def}, we immediately note that our proposed setting deals with non-symmetric test functions and directional information inhomogenously varying in $\Omega$, encoded in a weighted divergence term.

\subsection{Our proposal}
Our work extends the regularisers in \cref{eq: DTV}-\cref{eq: DTGV} for handling spatially varying directions $\thetabold$ in $\Omega\subset\RR^2$ instead of a fixed scalar direction $\theta$. 
We investigate the directional total variation regulariser of \cref{eq: TDV intro def} and study its performance for a variety of image processing problems by solving
\begin{equation}
u^\ast 
\in 
\argmin_u \sum_{q=1}^\Qrm \TDVM{q}{\alphabold_q}[u,\Mcal_q] + \frac{\eta}{2}\norm{\Scal u-u^\diamond}_2^2,
\label{eq: minimization intro}
\end{equation}
where $u^\diamond\in\LL{2}[\Omega]$ is a given, imperfect and possibly incomplete imaging data, and $\Scal:\LL{2}[\Omega]\to\LL{2}[\Omega]$ a linear operator. We consider the cases for which $\Mbold_j$ in $\Mcal_q=(\Mbold_j)_{j=1}^q$ from \cref{eq: TDV intro def} is
\[
\text{either}\quad\Mbold_j 
=
\Ibold\quad\text{or}\quad\Mbold_j=\Lambdabold_\bbold (\Rbold_{\thetabold})^\T,
\]
with $\Ibold$ the identity, $\Lambdabold_\bbold$ the contraction and $\Rbold_{\thetabold}$ the rotation matrices, defined as:
\begin{equation}
\Lambdabold_\bbold = \begin{pmatrix}
b_1(\xbold) & 0\\
0 & b_2(\xbold)\\
\end{pmatrix}\quad\text{and}\quad
\Rbold_\thetabold=\begin{pmatrix}
\cos\theta(\xbold)  & -\sin\theta(\xbold) \\
\sin\theta(\xbold) & \cos\theta(\xbold)
\end{pmatrix},
\label{eq: controrot matrices}
\end{equation}
and with $\bbold:=(b_1(\xbold),b_2(\xbold))^\T\in[0,1]^2$, $\thetabold:=\theta(\xbold)\in[0,2\pi)$.
Occasionally, we will use the vector field $\vbold=(\cos\thetabold,\sin\thetabold)^\T$ and its orthogonal $\vbold_\perp=(-\sin\thetabold,\cos\thetabold)^\T$. 
Thus, we interpret the core operation of the dual version of the regulariser in \cref{eq: TDV intro def}, 
$\Mbold_1\grad \otimes u$, as weighted directional derivatives of $u$ along $\vbold$ and $\vbold_\perp$ since
\begin{equation}
\Mbold_1\grad \otimes u 
= 
\Lambdabold_\bbold (\Rbold_\thetabold)^\T  \grad \otimes u
= 
\begin{pmatrix}
b_1\grad_\vbold u\\
b_2\grad_{\vbold_\perp} u
\end{pmatrix}.
\label{eq: directional equivalence}
\end{equation}
We will in particular focus on the case $\bbold=(1,\beta(\xbold))$ for $\beta(\xbold)\in[0,1]$ being either inhomogeneous or constant in $\Omega$, see \cref{rem: ellipses} for the geometrical interpretation when different choices are made for the constant.

\begin{remark}\label{rem: ellipses}
In \cref{fig: ellipse2} we simulate the two dimensional behaviour of \cref{eq: directional equivalence} for different choices of $b_2$. More precisely, for a continuous imaging function $u:\Omega\to\RR$ we represent a possible situation at the position $\xbold\in\Omega$ of the vectors $\pbold=(p_1,p_2)=\grad u$ and $\vbold=(v_1,v_2)$, depicted with red and blue arrows, respectively. 
We also represent the components $\rbold=(r_1,r_2)=(p_1v_1,p_2v_2)$ of $\grad_\vbold u = r_1+r_2 = p_1v_1+p_2v_2$  by a green arrow.  
The vectors and the corresponding arrows are the same in all \cref{fig: ellipse_b11} to \cref{fig: ellipse_b10}.
Moreover, the test functions $\Psibold=(\Psi^1,\Psi^2)$ lie on the black circle due to the constraint $\norm{\Psibold}_2\leq 1$. Note that in the 2D domain we have
\[
\Mbold_1\grad  u \cdot \Psibold= \grad u\cdot \Mbold_1^\T \Psibold,
\]
which allows to change the metric space of the test functions into an elliptic ball in magenta. 
Being fixed $b_1=1$,  each figure corresponds to a particular choice of $b_2=\beta$ between 0 and 1. 
Finally,  the magenta arrow corresponds to the direction of $\Psibold$ which realizes the supremum of the regulariser $\TDVM{1}{\alphabold}$ in Equation \cref{eq: TDV intro def}. 
We observe in \cref{fig: ellipse_b10} the limit case $\bbold=(1,0)$ where $\TDVM{1}{\alphabold}[u,\Mcal]$ penalizes the rate of change of $u$ only along $\vbold$ without orthogonal $\vbold_\perp$ contribution. In all the other circumstances, $\vbold_\perp$ acts as quality estimation of $\vbold$, leading to a full isotropic approach in the case  $\bbold=(1,1)$ of \cref{fig: ellipse_b11}, since the magenta arrow is bended in the direction of the gradient $\grad u$ rather than the direction of $\grad_\vbold u$.
\end{remark}

\begin{figure}[tbh]
\centering
\begin{subfigure}[t]{0.16\textwidth}
\includegraphics[width=\textwidth,trim={5.5cm 2cm 4.5cm 1cm},clip=true]{{{./images/intro/TDVexample/1}.png}}
\caption{$\bbold=(1,1)$}
\label{fig: ellipse_b11}
\end{subfigure}
\hfill
\begin{subfigure}[t]{0.16\textwidth}
\includegraphics[width=\textwidth,trim={5.5cm 2cm 4.5cm 1cm},clip=true]{{{./images/intro/TDVexample/0.8}.png}}
\caption{$\bbold=(1,0.8)$}
\end{subfigure}
\hfill
\begin{subfigure}[t]{0.16\textwidth}
\includegraphics[width=\textwidth,trim={5.5cm 2cm 4.5cm 1cm},clip=true]{{{./images/intro/TDVexample/0.6}.png}}
\caption{$\bbold=(1,0.6)$}
\end{subfigure}
\hfill
\begin{subfigure}[t]{0.16\textwidth}
\includegraphics[width=\textwidth,trim={5.5cm 2cm 4.5cm 1cm},clip=true]{{{./images/intro/TDVexample/0.4}.png}}
\caption{$\bbold=(1,0.4)$}
\end{subfigure}
\hfill
\begin{subfigure}[t]{0.16\textwidth}
\includegraphics[width=\textwidth,trim={5.5cm 2cm 4.5cm 1cm},clip=true]{{{./images/intro/TDVexample/0.2}.png}}
\caption{$\bbold=(1,0.2)$}
\end{subfigure}
\hfill
\begin{subfigure}[t]{0.16\textwidth}
\includegraphics[width=\textwidth,trim={5.5cm 2cm 4.5cm 1cm},clip=true]{{{./images/intro/TDVexample/0}.png}}
\caption{$\bbold=(1,0)$}
\label{fig: ellipse_b10}
\end{subfigure}
\caption{Different choices for $\bbold=(1,\beta)$ in $\TDVM{1}{\alphabold}[u,\Mcal]$, with $\Mcal=\Lambdabold_\bbold (\Rbold_{\theta})^\T$ and $\theta$ fixed.}
\label{fig: ellipse2}
\end{figure}

\subsection{Contribution of the paper}
In what follows we will derive:
\begin{itemize}
\item a rigorous definition of the total directional regulariser \cref{eq: TDV intro def};
\item a characterisation of \cref{eq: TDV intro def} that turns \cref{eq: minimization intro} into a form that is amenable for numerical solution. For this we propose a primal-dual algorithm and present certain instances for different combinations of orders $q=1,\dots,\Qrm$, up to $\Qrm=3$ in \cref{eq: TDV intro def};
\item a number of numerical experiments with this new regulariser for image denoising, image zooming and interpolation of two-dimensional surfaces from a sparse number of given height values.
\end{itemize}

\subsection{Organization of the paper}
In \cref{sec: TDV} we discuss the higher-order total directional variation regularisers with anisotropy.
The numerical details of the discretisation are introduced in \cref{sec: numerical details}, with the primal-dual algorithm and the numerical optimisation described in \cref{sec: numerical optimization}. Imaging applications to denoising, wavelet-based zooming and surface interpolation, e.g.\ in atomic force microscopy imaging, are discussed in \cref{sec: model-description-denoising} and \cref{sec: imaging applications}. 
\section{Higher-order total directional variation}\label{sec: TDV}
In this section we introduce the rigorous definition of \cref{eq: TDV intro def}. 
To do so, we first introduce the terminology of tensors and their mathematical manipulation.

\subsection{Tensors}
Following \cite{BreKunPoc2010}, let $\Tcal^\ell(\RR^d)$ 
be the vector space of $\ell$-tensors
defined as
\begin{align*}
\Tcal^\ell(\RR^d) 
&:=
\left\{\xibold:\underbrace{\RR^d\times \dots \times \RR^d}_{\ell\text{-times}}\to\RR,\text{ such that } \xibold \text{ is } \ell\text{-linear}\right\}.
\end{align*}
On $\Tcal^\ell(\RR^d)$, we have the following operations:
\begin{itemize}
\item let $\otimes$ be the tensor product for $\xibold_1\in\Tcal^{\ell_1}(\RR^d)$, $\xibold_2\in\Tcal^{\ell_2}(\RR^d)$, with $\xibold_1\otimes\xibold_2\in\Tcal^{\ell_1+\ell_2}(\RR^d)$:
\[
(\xibold_1\otimes\xibold_2)(\abold_1,\dots,\abold_{\ell_1 + \ell_2}) = \xibold_1(\abold_1,\dots,\abold_{\ell_1})\xibold_2(\abold_{\ell_1+1},\dots,\abold_{\ell_1+\ell_2});
\]
\item let $\tr(\xibold)\in\Tcal^{\ell-2}(\RR^d)$ be the trace of $\xibold\in\Tcal^\ell(\RR^d)$, with $\ell\geq 2$, defined by
\[
\tr(\xibold)(\abold_1,\dots,\abold_{\ell-2}) = \sum_{i=1}^d \xibold(\ebold_i,\abold_1,\dots,\abold_{\ell-2},\ebold_i),
\]
where $\ebold_i$ is the $i$-th standard basis vector;
\item let $(\blank)^\sim$ be such that if $\xibold\in\Tcal^\ell(\RR^d)$, then
$\xibold^\sim(\abold_1,\dots \abold_\ell) = \xibold(\abold_\ell, \abold_1,\dots,\abold_{\ell-1})$;
\item let $\overline{(\blank)}$ be such that if $\xibold\in\Tcal^\ell(\RR^d)$, then
$
\overline{\xibold}(\abold_1,\dots \abold_\ell) = \xibold(\abold_\ell,\dots,\abold_1);
$
\item let $\xibold,\etabold\in\Tcal^\ell(\RR^d)$. The space $\Tcal^\ell(\RR^d)$ is equipped with the scalar product defined as
\[
\xibold\cdot\etabold = \sum_{p\in\{1,\dots,d\}^\ell} \xibold_{p_1,\dots,p_\ell} \etabold_{p_1,\dots,p_\ell}.
\]
\end{itemize}

We now introduce the derivative operator for tensors and its weighted version.
\begin{definition}\label{def: weighted derivative tensors}
Let $\grad=(\partial_1,\dots, \partial_d)^\T$ be the derivative operator and $\xibold\in\Tcal^\ell(\RR^d)$.
The derivative of $\xibold$ is defined as $(\grad\otimes \xibold)\in\Tcal^{\ell+1}(\RR^d)$ via the following:
\[
\grad\otimes \xibold
:= 
\left(\partial_{j} \xibold_{i_1,\dots, i_\ell}\right)_{j,i_1,\dots, i_\ell}.
\]
Let $\etabold\in\Tcal^2(\RR^d)$. The derivative operator weighted by $\etabold$ is defined as $\etabold\grad\in\Tcal^{1}(\RR^d)$ and the derivative of $\xibold\in\Tcal^\ell(\RR^d)$ weighted by $\etabold$ is defined as $(\etabold\grad\otimes \xibold)\in\Tcal^{\ell+1}(\RR^d)$ via the following:
\begin{equation}
\etabold\grad\otimes \xibold 
:= 
\left( 
\sum_{k=1}^d \etabold_{j,k} \partial_k \xibold_{i_1,\dots, i_\ell}
\right)_{j,i_1,\dots,i_\ell}.
\label{eq: sum}
\end{equation}
\end{definition}

\begin{remark}
For notational purposes, the sum in \cref{eq: sum} will be shortened using Einstein notation over the repeated subscript, meaning that each element of the tensor $\etabold\grad\otimes\xibold$ is written as $\etabold_{j,k} \partial_k \xibold_{i_1,\dots, i_\ell}$.
\end{remark}

In what follows, we will also denote the space
of  $q$-times uniformly continuously differentiable $\Tcal^\ell(\RR^d)$-valued tensors  as $\CCCspace{q}[\overline{\Omega},\Tcal^\ell(\RR^d)]$ which is a Banach space with the norm
\[
\norm{u}_{\infty,q} = \max_{\ell=0,\dots,q} \sup_{x\in\Omega} \abs{\grad^\ell\otimes u(\xbold)},
\]
where $(\grad^q\otimes u): \Omega\to\Tcal^{q+\ell}(\RR^d)$, and we will consider also the space
$
\CCCcomp{q}[\Omega,\Tcal^\ell(\RR^d)] 
$ of $\Tcal^\ell(\RR^d)$-valued tensors which are $q$-times continuously differentiable with compact support in $\Omega$.

\subsection{Definition of total directional variation}\label{sec: higher-order tdv}
For making sense of the distributional formulation of higher-order directional variation in \cref{eq: TDV intro def} we need an integration by parts formula for the weighted derivative of tensors in \cref{def: weighted derivative tensors}. Namely we consider
\[
\int_\Omega (\Mbold\grad\otimes \Abold)\cdot\Psibold \diff \xbold,
\] 
with $\Omega\subset\RR^d$ being a bounded Lipschitz domain, $\Mbold\in\CCCspace{1}[\Omega,\Tcal^2(\RR^d)]$, $\Abold\in\CCCspace{1}[\Omega,\Tcal^\ell(\RR^d)]$ and $\Psibold\in\CCCcomp{1}[\Omega,\Tcal^{\ell+1}(\RR^d)]$.
We report in this section the main results from the second part of our companion work \cite{ParMasSch18analysis}, where detailed proofs can be found.
First, we give an integration by parts formula where only $\Mbold$ switches:
\begin{lemma}\label{lem: move M}
Let $\Omega$, $\Mbold$, $\Abold$ and $\Psibold$ as above. Then:
\begin{equation}
\int_\Omega (\Mbold\grad\otimes\Abold)\cdot\Psibold \diff\xbold =\int_\Omega (\grad\otimes \Abold)\cdot \tr\left(\Mbold\otimes\Psibold^\sim\right)\diff\xbold,\,\text{ for all }\Mbold,\,\Abold,\,\Psibold.
\end{equation}
\end{lemma}
Then a general adjoint property follows:
\begin{lemma}\label{lem: divergence modified}
Let $\Omega$, $\Mbold$, $\Abold$ and $\Psibold$ as above. Then:
\begin{equation}
\int_\Omega (\Mbold\grad\otimes \Abold)\cdot\Psibold \diff \xbold = 
-\int_\Omega \Abold \cdot \div_{\Mbold}\Psibold\diff\xbold,\,\text{ for all }\Mbold,\,\Abold,\,\Psibold,
\end{equation}
where
$
\div_{\Mbold}\Psibold := \tr\left(\grad\otimes \left[\tr\left(\Mbold\otimes\Psibold^\sim\right)\right]^\sim\right).
$
\end{lemma}

We can now define the \emph{total directional variation} of order $q$ with weights $\alphabold\in\RR_+^q$.
\begin{definition}\label{def: TDV}
Let $\Omega\subset\RR^d$, $u\in\LL{1}[\Omega,\RR]$, $q\in\NN$, $\Mcal := (\Mbold_j)_{j=1}^q$ be a collection of fields in $\CCCspace{\infty}[\Omega,\Tcal^2(\RR^d)]$ and $\alphabold:=(\alpha_0,\dots,\alpha_{q-1})$ be a positive weight vector. Then, the total directional variation of order $q$, associated to $\Mcal$ and $\alphabold$, is defined as:
\begin{equation}
\TDVM{q}{\alphabold}[u,\Mcal]
:= 
\sup_\Psibold 
\left\{ 
\int_\Omega u \div_{\Mcal}^q \Psibold \diff \xbold \,\Big\lvert\,\text{for all } \Psibold \in\Ycal_{\Mcal,\alphabold}^q \right\},
\label{eq: TDV definition M}
\end{equation}
where
\begin{equation}
\Ycal_{\Mcal,\alphabold}^q 
= 
\left\{
\Psibold\,:\, \Psibold \in \CCCcomp{q}[\Omega,\Tcal^q(\RR^d)],\, \norm{\div_{\Mcal}^j\Psibold}_\infty\leq \alpha_j,\forall\,j=0,\dots q-1
\right\},
\label{eq: Ycal}
\end{equation}
and the weighted divergence of order $q$ is defined recursively, from \cref{lem: divergence modified}, as:
\begin{equation}
\begin{aligned}
\div^0_{\Mcal}\Psibold &:=\Psibold, &\text{if } j&=0,
\\
\div^1_{\Mcal}\Psibold &:= \div_{\Mbold_q}\Psibold, & \text{if } j&=1,
\\
&\vdots & &\vdots\\
\div_{\Mcal}^q(\Psibold) &:=  
\div_{\Mbold_{q-j+1}}
\left(
\div_{\Mcal}^{j-1} \Psibold
\right) 
&\text{if } j&=2,\dots,q.
\end{aligned}
\label{eq: weighted divergence}
\end{equation}
\end{definition}

\begin{remark}
For $\Mcal=(\Ibold)_{j=1}^{q}$, then $\TDVM{q}{\alphabold}[u,\Mcal]\equiv\neg\mathrm{sym}\TGV{\alphabold}{q}[u]$, see \cite{ParMasSch18analysis} and \cite[Remark 3.10]{BreKunPoc2010}.
\end{remark}

\subsection{Directional matrices for applications}
In what follows, we introduce a particular parametrisation of directional matrices for fields $\Mcal$ in \cref{eq: TDV definition M}.
For standard imaging applications,  we will usually deal with grey-scale images $u:\Omega\to\RR$,  $\Omega\subset\RR^2$, i.e.\ $d=2$.

\begin{definition}[Directional matrices]\label{def: directional matrices}
Let $\left(\bbold^j\right)_{j=1}^q$, $\bbold^j:\Omega\to [0,1]^2$, be a collection of so-called contraction weights (being each element of modulus $\leq 1$), $(\thetabold^j)_{j=1}^q$, $\thetabold^j : \Omega \to [0,2\pi)$, be a collection of angles, and $\Lambdabold_{\bbold}^j$ and $\Rbold_{\thetabold}^j$ the associated contraction and rotation matrices defined, respectively, as
\[
\Lambdabold_{\bbold}^j
:= 
\begin{pmatrix}
\bbold_1^j & 0 \\
0 & \bbold_2^j \\
\end{pmatrix}, \quad 
\Rbold_{\thetabold}^j
:=
\begin{pmatrix}
\cos\thetabold^j & -\sin\thetabold^j\\
\sin\thetabold^j & \cos\thetabold^j
\end{pmatrix}\in\SO{2}.
\]
Then we define $\Mcal:=(\Mbold_{j})_{j=1}^q$ to be a collection of contraction-rotation matrices (in Einstein notation) as
\[
\Mbold_{j} 
:= 
\Lambdabold_{\bbold}^j (\Rbold_{\thetabold}^j)^\T 
=  
\left(
\lambda_{pk}^j r_{\ell k}^j \right)_{p,\ell},
\]
where $\lambda_{pk}^j, r_{k\ell}^j$ are the element-wise entries of the matrices $\Lambdabold_\bbold^j$, $\Rbold_\thetabold^j$, respectively.
\end{definition}

\begin{definition}[Weighted derivatives of order 1]\label{def: weighted 1 derivatives}
Let $\grad=(\partial_1,\partial_2)^\T$ be the derivative operator. The gradient of a differentiable imaging function $u$ is given by $
\grad u := \grad\otimes u = 
\left( \partial_1 u, \partial_2 u\right)^\T
$
and the weighted derivative operator of order 1 associated to the directional matrix $\Mbold_{1}$ from \cref{def: directional matrices} is
\[
\Mbold_1 \grad u 
:= 
\left(\Mbold_1\grad \otimes u\right)_p
=
\left(
\lambda_{pk}^1 r_{\ell k}^1 \partial_\ell u
\right)_p.
\]
\end{definition}

\begin{remark}
If $\Mbold_1=\Ibold$ (i.e.\ $b_1^1\equiv b_2^1\equiv1$ and $\thetabold^q\equiv 0$ for all $\xbold\in\Omega$), then $\Mbold_1\grad u \equiv \grad u$.
\end{remark}

\begin{remark}\label{rem: equiv overline D}
Given $\thetabold_1$, let $\vbold^1=(\cos\thetabold^1,\sin\thetabold^1)$ and $\vbold_\perp^1=(-\sin\thetabold^1,\cos\thetabold^1)$. Then
\[
\Mbold_1\grad \otimes u :=
\begin{pmatrix}
\bbold_1^1 \grad_{\vbold^1} u \\
\bbold_2^1 \grad_{\vbold_\perp^1} u
\end{pmatrix},
\]
where $\grad_\zbold u$ represents the directional derivative along a vector field $\zbold$, defined as
\[
\grad_{\zbold}{u}(\xbold) = \grad u(\xbold) \cdot \zbold = \sum_{i=1}^2 \partial_i u\,z_i.
\]
\end{remark}

\begin{definition}[Weighted derivatives of order \texorpdfstring{$q$}{q}]
\label{eq: weighted Q gradient}
We define the derivative of order $q$ of $u$ using \cref{def: weighted 1 derivatives} recursively as
\[
\grad^q u 
:= 
\left( \grad\otimes \grad^{q-1}u \right)_{p_q,\dots,p_1} 
= 
\left(
\partial_{p_q} \dots \partial_{p_1} u
\right)_{p_q,\dots,p_1} .
\]
We define the weighted derivative of order $q$ of $u$ with respect to $\Mcal$ recursively as
\begingroup\makeatletter\def\f@size{10.5}\check@mathfonts
\def\maketag@@@#1{\hbox{\m@th\normalsize\normalfont#1}}
\[
\begin{aligned}
\grad_\Mcal^0 u &:= u,\\
\grad_\Mcal^1 u &:= (\Mbold_1\grad \otimes u)_{p_1} = 
\left(\lambda_{p_1 k}^1 r_{\ell k}^1 \partial_\ell u\right)_{p_1},\\
&\vdots\\
\grad_\Mcal^q u  
&:= 
\left( (\Mbold_q \grad)\otimes\left(
\Mbold_{q-1}\grad^{q-1} u\right) \right
)_{p_{q},\dots,p_{1}}
=
\left(
\lambda_{p_q k_q}^q r_{\ell_q k_q}^q \partial_{\ell_q}\left(\Mbold_{q-1}\grad^{q-1} \otimes u\right)\right)_{p_q,p_{q-1},\dots,p_{1}}.
\end{aligned}
\]
\endgroup
\end{definition}

\subsection{Examples}
We present some examples of the total directional variation of order $q$ for $q=1,2,3$, $\alphabold=(\alpha_j)_{j=0}^{q-1}$ and a collection of directional matrices $\Mcal=(\Mbold_j)_{j=1}^q$: 
\begin{itemize}
\item order $q=1$ and $\Mcal=(\Mbold_1)$:
\[
\begin{aligned}
\TDVM{1}{\alphabold}[u,\Mcal]&:= \sup_\Psibold \left\{ \int_\Omega u \div _{\Mbold_1}\Psibold \diff \xbold \,\Big\lvert\,\text{for all } \Psibold \in \Ycal_{\Mcal,\alphabold}^1 \right\},\\
\Ycal_{\Mcal,\alphabold}^1 
&= 
\left\{ \Psibold \,:\, \Psibold \in\CCCcomp{\infty}[\Omega,\RR^2]\text{ s.t.\ } \norm{\Psibold}_\infty\leq \alpha_0 \right\};
\end{aligned}
\]
\item order $q=2$ and $\Mcal=(\Mbold_1,\Mbold_2)$:
\[
\begin{aligned}
\TDVM{2}{\alphabold}[u,\Mcal] 
&:= 
\sup_\Psibold \left\{ \int_\Omega u \div_{\Mbold_1}(\div_{\Mbold_2} \Psibold) \diff \xbold \,\Big\lvert\,\text{for all } \Psibold \in \Ycal_{\Mcal,\alphabold}^2 \right\},
\\
\Ycal_{\Mcal,\alphabold}^2 
&= 
\left\{ \Psibold \,:\, \Psibold \in\CCCcomp{\infty}[\Omega,\RR^{2\times 2}]\text{ s.t.\ } \norm{\Psibold}_\infty\leq \alpha_0,\, \norm{\div_{\Mbold_2}\Psibold}_\infty\leq \alpha_1 \right\};
\end{aligned}
\]
\item order $q=3$ and $\Mcal=(\Mbold_1,\Mbold_2,\Mbold_3)$:
\[
\begin{aligned}
\TDVM{3}{\alphabold}[u,\Mcal] 
&:= 
\sup_\Psibold \left\{ \int_\Omega u \div_{\Mbold_1}\left(\div_{\Mbold_2}(\div_{\Mbold_3} \Psibold)\right) \diff \xbold \,\Big\lvert\,\text{for all } \Psibold \in \Ycal_{\Mcal,\alphabold}^3 \right\},
\\
\Ycal_{\Mcal,\alphabold}^3
&= 
\left\{ \Psibold \,:\, \Psibold \in\CCCcomp{\infty}[\Omega,\RR^{2\times 2\times 2}]\text{ s.t.\ }
\begin{pmatrix}
\norm{\Psibold}_\infty\leq \alpha_0
\\
\norm{\div_{\Mbold_3}\Psibold}_\infty\leq\alpha_1
\\
\norm{\div_{\Mbold_2}(\div_{\Mbold_3}\Psibold)}_\infty\leq \alpha_2
\end{pmatrix}
\right\}.
\end{aligned}
\]
\end{itemize}

\section{Numerical discretisation}\label{sec: numerical details}
The rest of the paper focuses on the discretised formulation of \cref{eq: minimization intro}, and its numerical solutions and performances on a number of image processing variational examples. 
We start by discretising the $\TDVM{q}{\alphabold}-\LL{2}$ problem in \cref{eq: minimization intro}.

\subsection{Staggered grids}
The discretisation of the $\TDVM{q}{\alphabold}-\LL{2}$ in
\cref{eq: minimization intro} is based on  finite-difference schemes for derivatives on staggered regular Cartesian grids of width $h>0$:
\begin{itemize}
\item the \emph{grid of pixels} $\Omega^h$, of axes $x_1$ and $x_2$ for a 2-dimensional domain and size $M\times N$, is defined as
\[
\Omega^h = \left\{ (k h, lh)\,\lvert\, (1,1)\leq(k,l)\leq (M,N)\right\};
\]
\item the \emph{grid of cell centres} $\Gamma^h$, of size $(M-1)\times(N-1)$ and used to perform the weighted derivative operation (i.e.\ for introducing the anisotropy, see the grid associated to the blue squares in \cref{fig: staggered grid}), is defined as:
\[
\Gamma^h = 
\left\{ 
(
\widetilde{k} h, \widetilde{l} h
)
\,\lvert\, 
(\widetilde{k}, \widetilde{l})=\left(k+\frac{1}{2},l+\frac{1}{2}\right),\,
\left(
1,1\right)\leq (k, l) < \left(M,N\right)
\right\};
\]
\item the collection of \emph{grids} $\left(\Xrm^{j,h}_{\iotabold}\right)_{j=1}^q$ associated to the differential operators involved, where $\iotabold=(\iota_1,\dots,\iota_j)$ is a multi-index variable and $\iota_s\in\{1,2\}$ for each $s=1,\dots,j$ indicates the partial derivative involved ($1$ for $\partial_{x_1}$ and $2$ for $\partial_{x_2}$). Every $\Xrm^{j,h}_{\iotabold}$ is a sub-collection of $2^j$ grids, each one of size $M\times N$ and denoted by $\Xrm^{j,h}_{(\iota_1,\dots,\iota_j)}$, each one associated to a fixed choice for the derivative operator $\left(\partial_{x_{\iota_j}}\cdots\partial_{x_{\iota_1}}\right)$ considered:
\[
\Xrm^{j,h}_{(i_1,\dots,i_j)} = 
\left\{ 
(
m h, n h
)
\,\lvert\, 
(m, n)=
\left(k+\frac{\abs{I_1}_\#}{2},l+\frac{\abs{I_2}_\#}{2}\right),\,
\left(
1,1\right)\leq (k,l) \leq \left(M,N \right)
\right\},
\]
where $\abs{I_1}_\#$ and $\abs{I_2}_\#$ are the cardinality of the sets $I_1$ and $I_2$ containing as many elements as the number of derivatives along the axes $x_1$ and $x_2$, respectively. A visual representation of such grids is given in \cref{fig: staggered grid}. 
%
For example:
\begin{itemize}
\item with a bit of abuse of notation, if $j=0$ then $\Xrm^{0,h}_{(-)}$ coincides with $\Omega^h$;
\item $\Xrm^{1,h}_{(1)},\Xrm^{1,h}_{(2)}$ result in $\Omega^h$ shifted by $h/2$ along $x_1$ and $x_2$ axes, respectively;
\item if $q=3$ and $\iotabold=(2,1,1)$, then we are referring to the grid associated to one out of the eight possible combinations for the third order derivative $\grad^3$, namely $\partial_{x_2}\partial_{x_1}\partial_{x_1}$, which is located on the grid identified by our notation $\Xrm_{(1,1,2)}^{3,h}$.
\end{itemize}
\end{itemize}

\subsection{Discretised objects}\label{sec: discretised objects}
Let the order of derivatives $q>0$ be fixed.
By means of the superscript $h$, we define the finite-dimensional approximation of the following quantities, where $|\Omega^h|$, $|\Gamma^h|$ and $|\Xrm^{j,h}_{(\iota_1,\dots,\iota_j)}|$ are the number of grid points in $\Omega^h$, $\Gamma^h$ and $\Xrm^{j,h}_{(\iota_1,\dots,\iota_j)}$, respectively:
\begin{itemize}
\item $\ubold^h\in\RR^{|\Omega^h|}$ is the discretisation of the function $u$;\item $\ubold^{\diamond,h}$ is the discretisation of the observed imaging data $u^\diamond$; 
\item $\vbold^h=(\vbold_1^h,\vbold_2^h)\in \RR^{|\Gamma^h|\times|\Gamma^h|}$ is a discrete vector field;
\item $\bbold^h=(\bbold_1^h,\bbold_2^h)\in \RR^{|\Gamma^h|\times|\Gamma^h|}$ are discrete contraction weights for $\Lambda_\bbold$;
\item $\Mbold_{j}^h\in \RR^{|\Gamma^h|\times |\Gamma^h|\times|\Gamma^h| \times |\Gamma^h|}$ discretises the weights $\Mbold_{j}\in\Tcal^2(\RR^2)$, for each $j=1,\dots,q$;
\item $\Mcal^h=(\Mbold_{j}^{h})_{j=1}^q$ 
discretises the collection of weights $\Mcal=(\Mbold_{j})_{j=1}^q$;
\item $\Psibold^{h}=(\Psibold_1^{h},\dots,\Psibold_{2^q}^{h})\in \RR^{|\Gamma^h| \times \dots \times |\Gamma^h|}$ discretises the test functions $\Psibold\in \Tcal^q(\RR^2)$;
\item $\zbold^h=(\zbold_j^h)_{j=0}^{q-1}$ discretises the primal variables $\zbold$, 
with $\zbold_0^h=\ubold^h\in\RR^{|\Omega^h|}$, $\zbold_0^{\diamond,h}=\ubold^{\diamond,h}$ 
and each $\zbold_{j}^h \in \RR^{|\Xrm_{(1,\dots,1)}^{j,h}| \times\dots\times |\Xrm_{(2,\dots,2)}^{j,h}|}$,
for $j=1,\dots, q-1$;
\item  $\wbold^h = (\wbold_{j}^h)_{j=1}^q$ discretises the dual variables $\wbold$, with $\wbold_{j}^h \in \RR^{|\Xrm_{(1,\dots,1)}^{j,h}| \times\dots\times |\Xrm_{(2,\dots,2)}^{j,h}|}$ for $j=1,\dots,q$.
\end{itemize}

\subsection{Isotropic differential operators}
Here we discuss the discretization of the adjoint unweighted operators $\grad$ and $\div$. 
For $\ubold^h\in\RR^{|\Omega^h|}$, the \emph{discrete gradient} operator is defined as
\[
\begin{aligned}
\grad^h\,:\,\RR^{|\Omega^h|} &\to \RR^{|\Xrm_{(1)}^{1,h}|\times|\Xrm_{(2)}^{1,h}|}
\\
\ubold^h &\xmapsto{} (\partial_1^h \ubold^h, \partial_2^h \ubold^h),
\end{aligned}
\]
where we use the \emph{central second-order finite difference scheme} on the grids $\Xrm_{(1)}^{1,h},\,\Xrm_{(2)}^{1,h}$:
\[
{(\partial_1^h \ubold^h)}_{k+\frac{1}{2},l}
=
\begin{dcases}
\frac{u^h_{k+1,l}-u^h_{k,l}}{2\left(\frac{h}{2}\right)} &\text{if } k<M,\\
0                      &\text{if } k=M,
\end{dcases}
\quad 
\text{and}
\quad
{(\partial_2^h \ubold^h)}_{k,l+\frac{1}{2}}
=
\begin{dcases}
\frac{u^h_{k,l+1}-u^h_{k,l}}{2\left(\frac{h}{2}\right)} &\text{if } l<N,\\
0                      &\text{if } l=N.
\end{dcases}
\]
Let $\pbold^h=(\pbold_1^h, \pbold_2^h)\in \RR^{|\Xrm_{(1)}^{1,h}|\times|\Xrm_{(2)}^{1,h}|}$ and let
the \emph{discrete divergence} operator
\[
\begin{aligned}
\div^{h} \,:\,\RR^{|\Xrm_{(1)}^{1,h}|\times|\Xrm_{(2)}^{1,h}|} &\to\quad \RR^{|\Omega^h|}
\\
\pbold^h \qquad &\xmapsto{} \grad^h \cdot\, \pbold^h,
\end{aligned}
\]
be defined for each pixel $(k,l)$ via the \emph{central second-order difference scheme} on $\Omega^h$:
\[
(\div^{h}\pbold^h)_{k,l} = 
\begin{dcases}
\frac{(\pbold_1^h)_{k+\frac{1}{2},l}}{2\left(\frac{h}{2}\right)} &\text{if } k=1,\\
\frac{ (\pbold_1^h)_{k+\frac{1}{2},l}  -  (\pbold_1^h)_{k-\frac{1}{2},l}}{2\left(\frac{h}{2}\right)}                     &\text{if } k\in(1,M),
\\
-\frac{(\pbold_1^h)_{k-\frac{1}{2},l}}{2\left(\frac{h}{2}\right)} &\text{if } k=M,
\end{dcases} 
+ 
\begin{dcases}
\frac{(\pbold_2^h)_{k,l+\frac{1}{2}}}{2\left(\frac{h}{2}\right)}     &\text{if } l=1,\\
\frac{(\pbold_2^h)_{k,l+\frac{1}{2}}  -  (\pbold_2^h)_{k,l-\frac{1}{2}}}{2\left(\frac{h}{2}\right)}                     &\text{if } l\in(1,N),
\\
- \frac{(\pbold_2^h)_{k,l-\frac{1}{2}}}{2\left(\frac{h}{2}\right)} &\text{if } l=N.
\end{dcases} 
\] 

Thus, the isotropic \emph{discrete gradient} and \emph{discrete divergence} are designed to fulfil the discrete adjointness property, for every 
$\ubold^h\in\RR^{|\Omega^h|}$ and $\pbold^h\in \RR^{|\Xrm_{(1)}^{1,h}|\times |\Xrm_{(2)}^{1,h}|}$:
\begin{equation}
\langle \grad^h \ubold^h,\, \pbold^h \rangle_{\Gamma^h}
= 
\langle  \ubold^h, \div^{h} \pbold^h \rangle_{\Omega^h},
\end{equation}
where $\langle\blank,\,\blank\rangle_{\Gamma^h}:\RR^{|\Xrm_{(1)}^{1,h}|\times |\Xrm_{(2)}^{1,h}|} \times \RR^{|\Xrm_{(1)}^{1,h}|\times |\Xrm_{(2)}^{1,h}|} \to \RR$ and $\langle\blank,\,\blank\rangle_{\Omega^h}:\RR^{|\Omega^h|} \times \RR^{|\Omega^h|} \to \RR$. 

For higher-order derivatives of order $q$ we denote the isotropic \emph{discrete gradient} and \emph{discrete divergence} operator by $\grad^{q,h}$ and $\div^{q,h}$ and write
\[
\begin{aligned}
\grad^{q,h}\,:\,\RR^{|\Omega^h|} 
&\to 
\RR^{|\Xrm_{(1,\dots,1)}^{q,h}|\times\dots\times |\Xrm_{(\iota_1,\dots,\iota_q)}^{q,h}|\times\dots\times|\Xrm_{(2,\dots,2)}^{q,h}|}
\\
\ubold^h 
&\xmapsto{} 
\left(
\partial_1^h 
(\partial_1^h(\dots (\partial_1^h \ubold^h))),\,
\dots,
\, 
\partial_{\iota_q}^h(\partial_{\iota_{q-1}}\dots (\partial_{\iota_1}^h \ubold^h))),\,
\dots,
\,
\partial_2^h (\partial_2^h(\dots \partial_2^h \ubold^h)))
\right),
\end{aligned}
\]
and
\[
\begin{aligned}
\div^{q,h} \,:\,\RR^{|\Xrm_{(1,\dots,1)}^{q,h}|\times\dots\times|\Xrm_{(2,\dots,2)}^{q,h}|} &\to\qquad\RR^{|\Omega^h|}
\\
\pbold^h \qquad &\xmapsto{}  \grad^{q,h} \cdot \pbold^h.
\end{aligned}
\]
The adjointness property is fulfilled for every $\ubold^h\in\RR^{|\Omega^h|}$ and $\pbold^h\in \RR^{|\Xrm_{(1,\dots,1)}^{q,h}|\times\dots\times|\Xrm_{(2,\dots,2)}^{q,h}|}$:
\begin{equation}
\langle \grad^{q,h} \ubold^h,\, \pbold^h \rangle_{\Gamma^h}
=
\langle \ubold^h,\,\div^{q,h} \pbold^h\rangle_{\Omega^h},
\end{equation}
with $\langle\blank,\blank\rangle_{\Gamma^h}: \RR^{|\Xrm_{(1,\dots,1)}^{q,h}|\times\dots\times |\Xrm_{(2,\dots,2)}^{q,h}|} \times \RR^{|\Xrm_{1,\dots,1}^{q,h}|\times\dots\times |\Xrm_{(2,\dots,2)}^{q,h}|} \to \RR$ and $\langle\blank,\blank\rangle_{\Omega^h}:\RR^{|\Omega^h|} \times \RR^{|\Omega^h|} \to \RR$. 

\subsection{Transfer operators}
The offset in the location between $\zbold_0^h,\, \grad^{1,h} \zbold_0^h,\,  \dots,\,\grad^{q,h} \zbold_0^h$ and the fields $\Mcal^h$ associated to $\vbold^h$ requires the introduction of \emph{transfer operators}, a concept from multigrid methods \cite{trottenberg2000multigrid}, so as to make the quantities computable in the same location. 
In what follows, we will provide some insights for the general case.

Let $\Wcal:=(\Wcal^j)_{j=1}^q$ be a family of transfer operators $\Wcal^j=(\Wbold_\iotabold^j)$, with $\Wbold_\iotabold^j: \RR^{|\Xrm_\iotabold^{j,h}|}\to \RR^{|\Gamma^h|}$ and $\iotabold$ a multi-index variable, with entries in $\{1,2\}$ similarly as for the staggered grids $\Xrm_{\iotabold}^{j,h}$. 
The idea is that $\Wcal^j$ interpolates the data from the grids of $j$-th order derivatives $\Xrm_\iotabold^{j,h}$ to the grid of cell centres $\Gamma^h$, e.g.\ $\Wbold^j_\iotabold$ is the operator made by partition of unit weights.
Since it is an averaging matrix, its adjoint operation is denoted by $(\Wcal^j)^{\T}$, where the extension from $\Gamma^h$ to the boundary of $\Xrm_\iotabold^{j,h}$ is made possible by mirroring the data as appropriate.\begin{example}
For $\zbold_0^h\in\RR^{|\Omega^h|}$ and $q$ fixed,
the derivatives of $\zbold_0^h$ (up to order $q$) are
\[\begin{aligned}
\wbold_{1}^{h} := \grad^{1,h} \zbold_0^h 
&= 
(\partial_1^h \zbold_0^h,\,\partial_2^h \zbold_0^h) 
\in \RR^{|\Xrm^{1,h}_{(1)}|\times |\Xrm^{1,h}_{(2)}|},
\\
\wbold_{2}^{h} := \grad^{2,h} \zbold_0^h 
&= 
(\partial_{1}^h\partial_{1}^h \zbold_0^h,\, \partial_{1}^h\partial_{2}^h \zbold_0^h,\, \partial_{2}^h\partial_{1}^h \zbold_0^h,\,\partial_{2}^h\partial_{2}^h \zbold_0^h) 
\in \RR^{|\Xrm^{2,h}_{(1,1)}|\times |\Xrm^{2,h}_{(2,1)}| \times |\Xrm^{2,h}_{(1,2)}|\times |\Xrm^{2,h}_{(2,2)}|},\\
&\vdots\\
\wbold_{q}^{h} : = \grad^{q,h} \zbold_0^h
&=
(\partial_{1}^h\dots \partial_{1}^h \zbold_0^h,\, \dots,\,\partial_{2}^h\dots \partial_{2}^h \zbold_0^h)
\in 
\RR^{|\Xrm^{q,h}_{(1,\dots, 1)}|\times\dots\times |\Xrm^{q,h}_{(2,\dots, 2)}|},
\end{aligned}
\]
the transfer operators $\Wcal=(\Wcal^j)_{j=1}^q$ are
\[
\Wcal^j: =
\begin{dcases}
\Wbold_{(1,\dots, 1)}^j : \RR^{|\Xrm_{(1,\dots, 1)}^{j,h}|} \to \RR^{|\Gamma^h|},\\
\qquad\qquad\qquad\vdots\\
\Wbold_{(\iota_1,\dots, \iota_j)}^j : \RR^{|\Xrm_{(\iota_1,\dots, \iota_j)}^{j,h}|} \to \RR^{|\Gamma^h|},\\
\qquad\qquad\qquad\vdots\\
\Wbold_{(2,\dots, 2)}^j : \RR^{|\Xrm_{(2,\dots,2)}^{j,h}|} \to \RR^{|\Gamma^h|},
\end{dcases}
\]
and each $\Wbold^j_{(\iota_1,\dots,\iota_j)}$ is the interpolation matrix that interpolates the values of $(\partial_{\iota_j}\dots\partial_{\iota_1} \zbold_0^h)\in\RR^{|\Xrm_{(\iota_1,\dots,\iota_j)}^{j,h}|}$ to $\RR^{|\Gamma^h|}$ by an arithmetic mean.
For example, for the first order derivatives we have
\[
\Wcal^1\grad^{1,h}\otimes \ubold^h
=
(\Wbold_{(1)}^1 \partial_1^h\ubold^h,\Wbold_{(2)}^1 \partial_2^h\ubold^h)\in\RR^{|\Gamma^h|\times|\Gamma^h|},
\]
where, for $k=1,\dots,M-1$ and $l=1,\dots N-1$,
\[
\begin{aligned}
(\Wbold_{(1)}^1 \partial_1^h\ubold^h)_{k+\frac{1}{2},l+\frac{1}{2}} 
&=
\frac{(\partial_1^h\ubold^h)_{k+\frac{1}{2},l} + (\partial_1^h\ubold^h)_{k+\frac{1}{2},l+1}}{2},
\\
(\Wbold_{(2)}^1 \partial_2^h\ubold^h)_{k+\frac{1}{2},l+\frac{1}{2}} 
&= 
\frac{(\partial_2^h\ubold^h)_{k,l+\frac{1}{2}} + (\partial_2^h\ubold^h)_{k+1,l+\frac{1}{2}}}{2}.
\end{aligned}
\]
As a further example for the second order derivative case and $\Mbold_1^h=\Ibold$ (which implies no averaging on the first derivatives for our construction), we have
\[
\Wcal^2\grad^{1,h}\otimes \grad^{1,h}\otimes\ubold^h
=
(
\Wbold_{(1,1)}^2 \partial_1^h\partial_1^h\ubold^h,
\Wbold_{(2,1)}^2 \partial_1^h\partial_2^h\ubold^h,
\Wbold_{(1,2)}^2 \partial_2^h\partial_1^h\ubold^h,
\Wbold_{(2,2)}^2 \partial_2^h\partial_2^h\ubold^h),
\]
where $\Wcal^2\grad^{1,h}\otimes \grad^{1,h}\otimes \ubold^h\in\RR^{|\Gamma^h|\times|\Gamma^h|\times|\Gamma^h|\times|\Gamma^h|}$ and, for $k=1,\dots,M-1$ and $l=1,\dots N-1$,
\[
\begin{aligned}
(\Wbold_{(1,1)}^2 \partial_1^h\partial_1^h\ubold^h)_{k+\frac{1}{2},l+\frac{1}{2}} 
&=
\frac{(\partial_1^h\partial_1^h\ubold^h)_{k,l} + (\partial_1^h\partial_1^h\ubold^h)_{k-1,l+1}}{2},
\\
(\Wbold_{(2,1)}^2 \partial_1^h\partial_2^h\ubold^h)_{k+\frac{1}{2},l+\frac{1}{2}} 
&= 
(\partial_1^h\partial_2^h\ubold^h)_{k,l},
\\
(\Wbold_{(1,2)}^2 \partial_2^h\partial_1^h\ubold^h)_{k+\frac{1}{2},l+\frac{1}{2}} 
&=
(\partial_2^h\partial_1^h\ubold^h)_{k,l},
\\
(\Wbold_{(2,2)}^2 \partial_2^h\partial_2^h\ubold^h)_{k+\frac{1}{2},l+\frac{1}{2}} 
&=
\frac{(\partial_2^h\partial_2^h\ubold^h)_{k,l} + (\partial_2^h\partial_2^h\ubold^h)_{k+1,l-1}}{2},
\end{aligned}
\]
with $(\partial_1^h\partial_1^h\ubold^h)_{0,l} = (\partial_1^h\partial_1^h\ubold^h)_{1,l}$  and
$(\partial_2^h\partial_2^h\ubold^h)_{k,0} = (\partial_2^h\partial_2^h\ubold^h)_{k,1}$.
\end{example}
\begin{remark}
The choice of the staggered grid increases the accuracy of the solution and allows to compute the inner products between gradients and the vector fields onto a unique regular Cartesian grid of reference, avoiding offsets. Moreover, the transfer operators $\Wcal$ reduce the bandwidth of higher order finite difference matrices, improving the quality of the result and reducing the smoothing due to large stencils.
\end{remark}

\begin{remark}
Note that when $\Mbold_{j}^h=\Ibold$ and $\alpha_{j}=+\infty$ for every $j<q$, as in the applications described in \cref{sec: imaging applications} of this paper, the use of transfer operators is needed only for the outer derivative, i.e.\ the one associated to the weighting field $\Mbold_{q}^h$. 
\end{remark}

\begin{figure}[tbhp]
\centering
\begin{subfigure}[t]{0.325\textwidth}\centering
\includegraphics[width=0.9\textwidth,trim=4.5cm 2.5cm 2.5cm 4.5cm,clip=true]{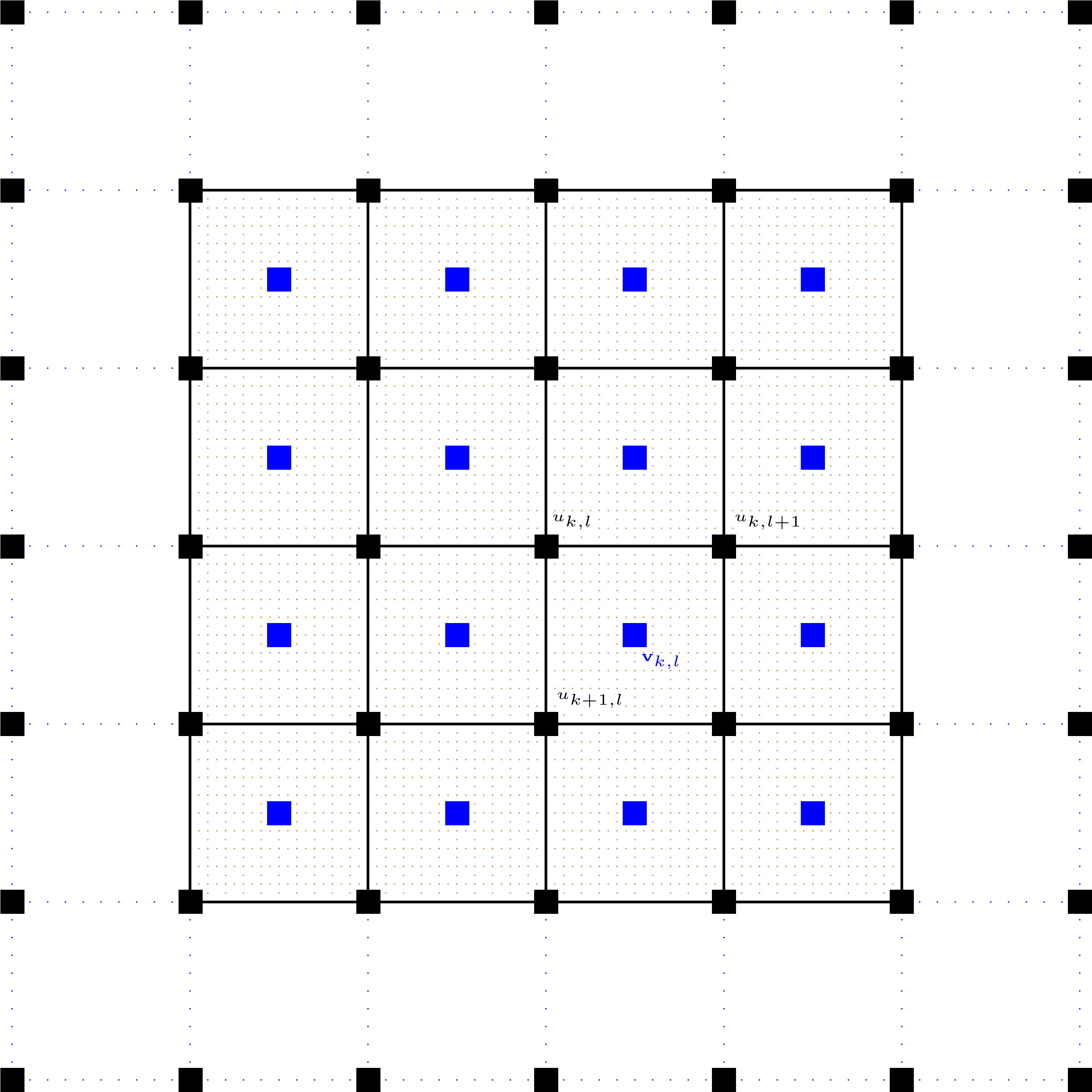}
\captionsetup{justification=centering}
\caption{$\ubold^h\in\RR^{|\Omega^h|}$ in $\textcolor{black}{\blacksquare}$, $\vbold^h\in\RR^{|\Gamma^h|}$ in $\textcolor{blue}{\blacksquare}$,\newline$(\bbold^h, \Mcal^h, \Psibold^h$ also in $\textcolor{blue}{\blacksquare}$).}
\end{subfigure}
\hfill
\begin{subfigure}[t]{0.325\textwidth}\centering
\includegraphics[width=0.9\textwidth,trim=4.5cm 2.5cm 2.5cm 4.5cm,clip=true]{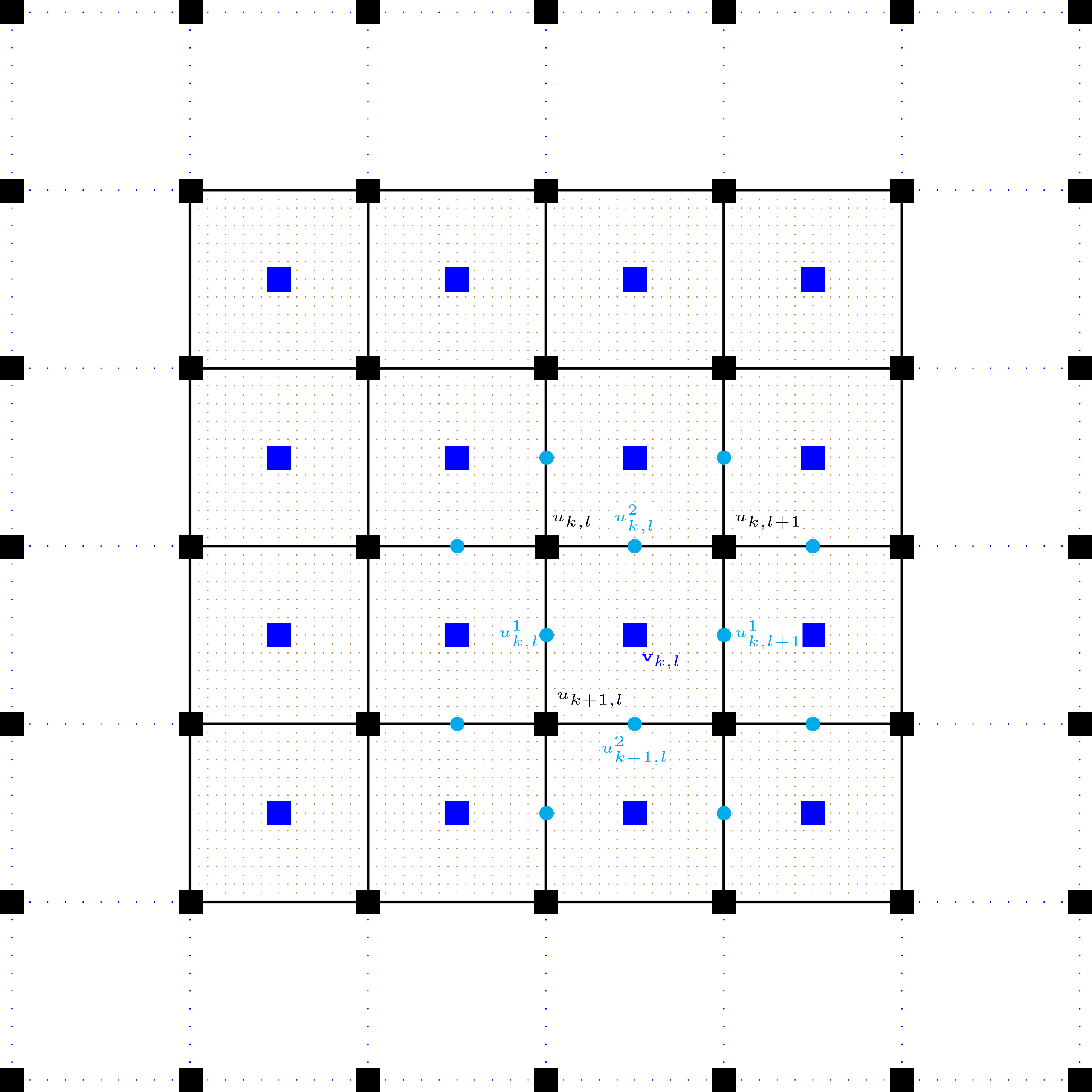}
\captionsetup{justification=centering}
\caption{First order derivatives,\newline$\grad \ubold^h =(\ubold^1,\ubold^2)$  in $\textcolor{cyan}{\bullet}$.} 
\end{subfigure}
\hfill
\begin{subfigure}[t]{0.325\textwidth}\centering
\includegraphics[width=0.9\textwidth,trim=4.5cm 2.5cm 2.5cm 4.5cm,clip=true]{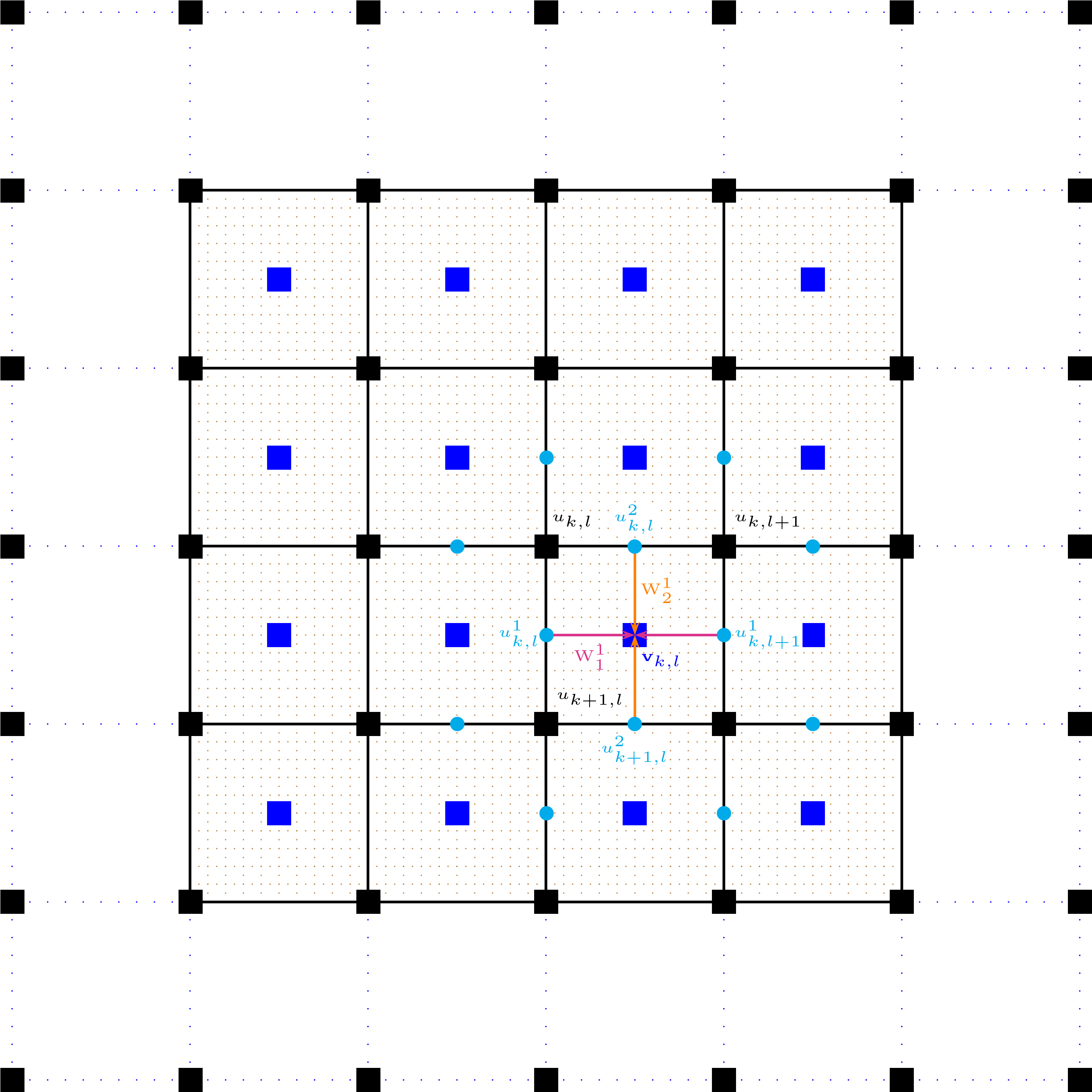}
\captionsetup{justification=centering}
\caption{Interpolating $\grad \ubold^h$ on $\Gamma^h$ via $\Wcal^1$:\newline $\textcolor{magenta}{\Wbold_1^1}$ acts on $\textcolor{cyan}{u^1}$, $\textcolor{orange}{\Wbold_2^1}$ acts on $\textcolor{cyan}{u^2}$.} 
\end{subfigure}
\\
\begin{subfigure}[t]{0.325\textwidth}\centering
\includegraphics[width=0.9\textwidth,trim=4.5cm 2.5cm 2.5cm 4.5cm,clip=true]{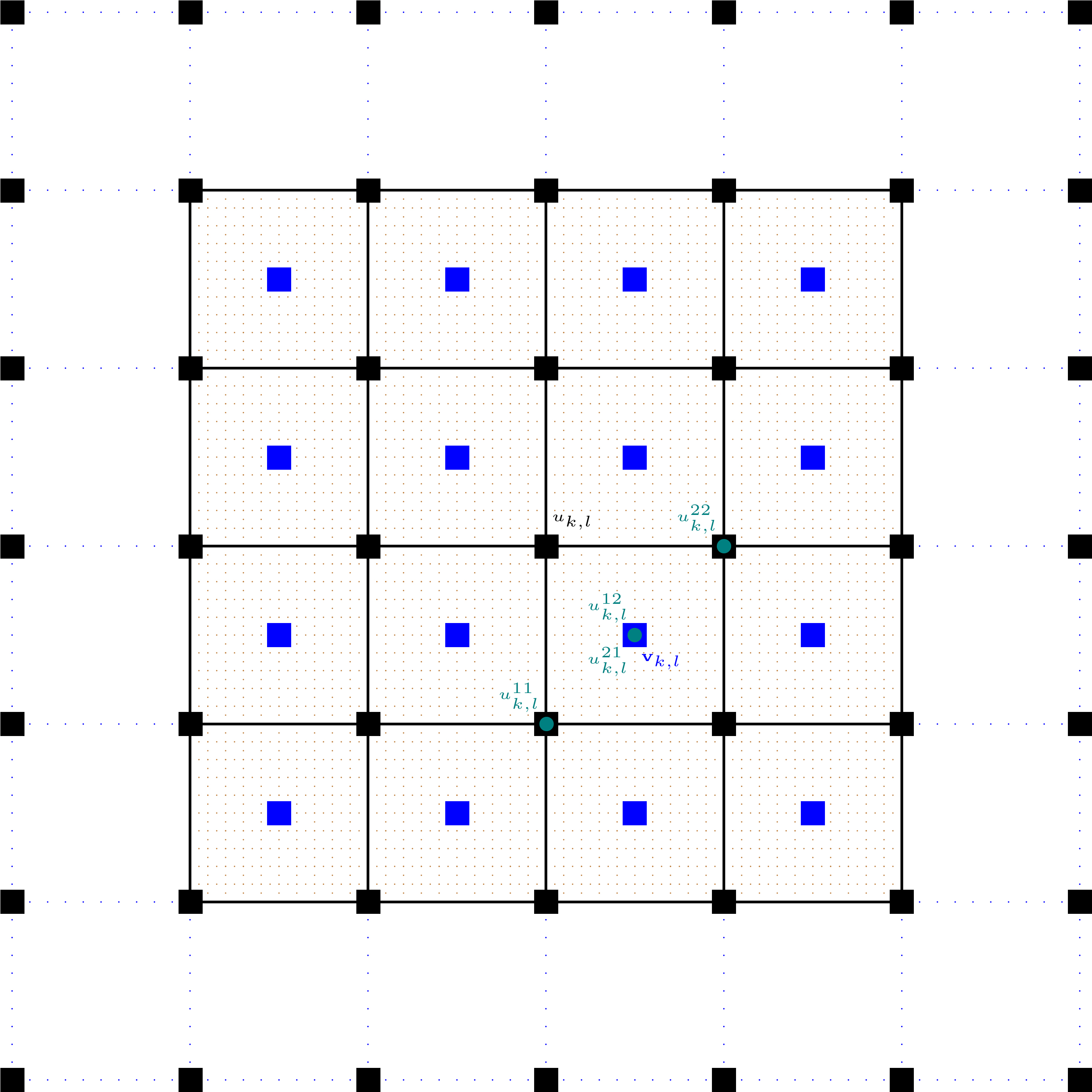}
\captionsetup{justification=centering}
\caption{Second order derivatives for $u_{k,l}$,\newline$\grad^2 \ubold^h =(\ubold^{11},\ubold^{12},\ubold^{21},\ubold^{22})$  in $\textcolor{teal}{\bullet}$.} 
\end{subfigure}
\hfill
\begin{subfigure}[t]{0.325\textwidth}\centering
\includegraphics[width=0.9\textwidth,trim=4.5cm 2.5cm 2.5cm 4.5cm,clip=true]{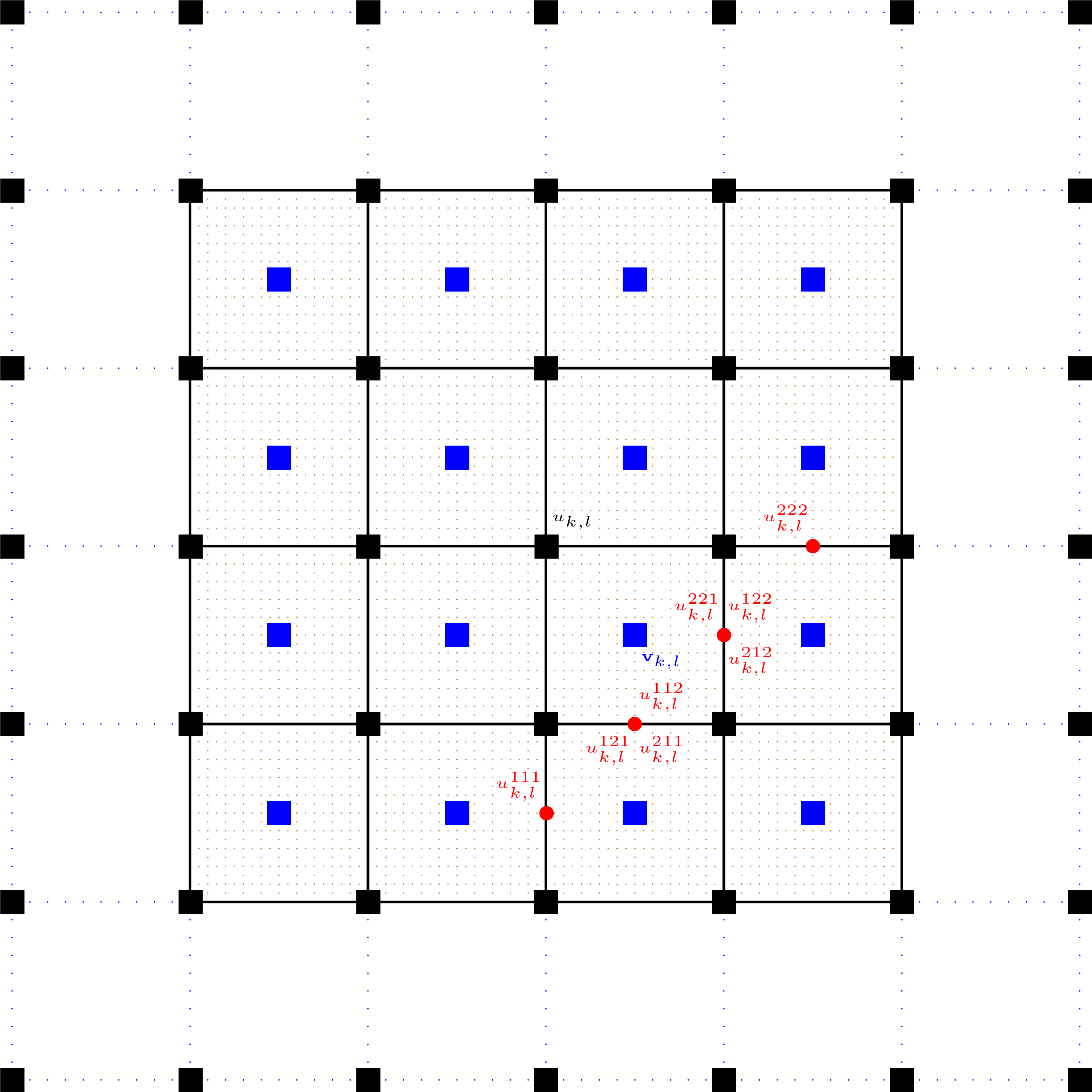}
\captionsetup{justification=centering}
\caption{Third order derivatives for $u_{k,l}$,\newline$\grad^3 \ubold^h =(\ubold^{111},\dots,\ubold^{222})$  in $\textcolor{red}{\bullet}$.} 
\end{subfigure}
\begin{subfigure}[t]{0.325\textwidth}\centering
\includegraphics[width=0.9\textwidth,trim=4.5cm 2.5cm 2.5cm 4.5cm,clip=true]{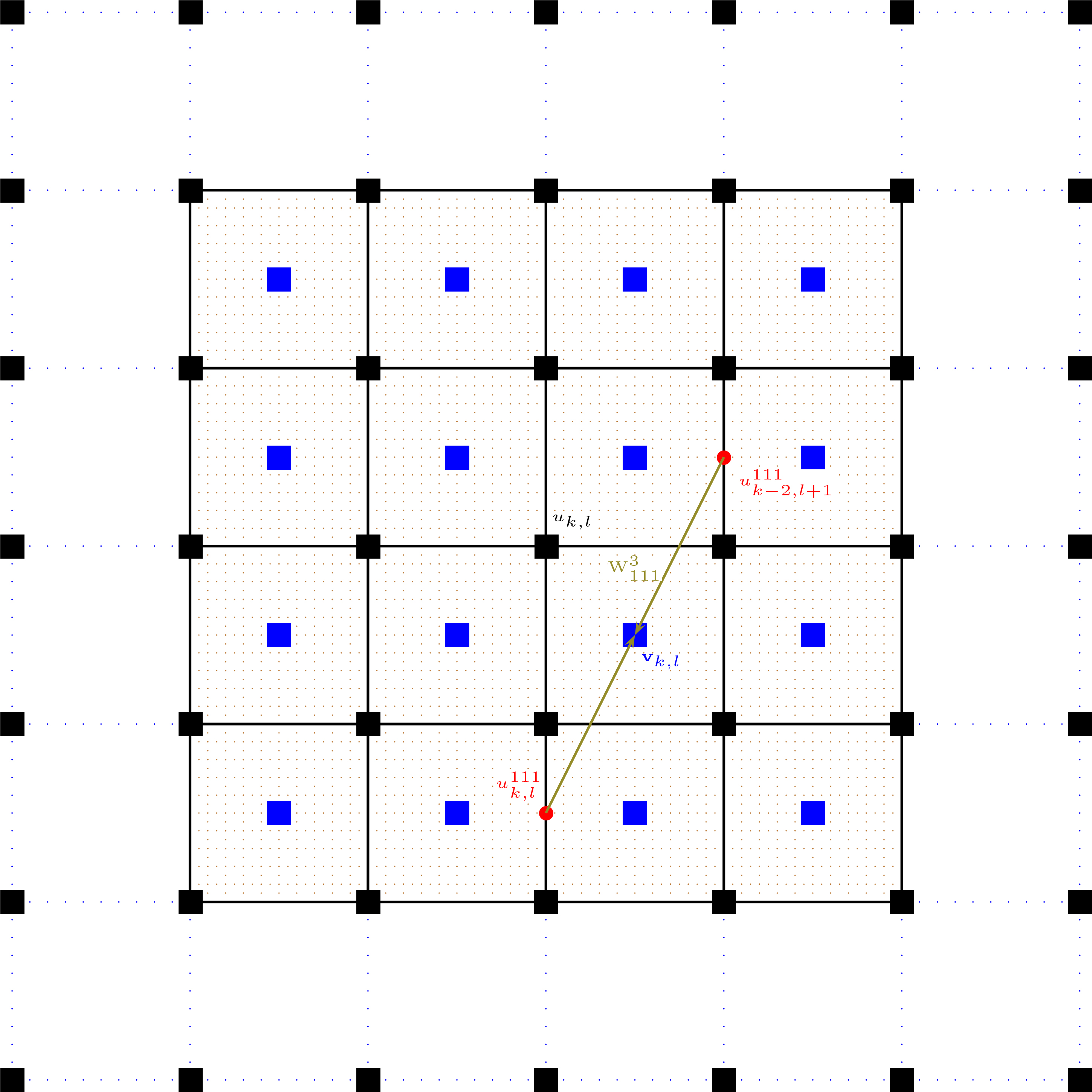}
\captionsetup{justification=centering}
\caption{Interpolating $\partial_1^3\ubold=\ubold^{111}$ on $\Gamma^h$\newline via $\textcolor{olive}{\Wbold_{111}^3}$ (part of $\Wcal^3$).} 
\end{subfigure}
\caption{Staggered grids with higher-order derivatives defined via iterative central finite differences of half stepsize, visualisation of variables and examples of transfer operators. 
Derivatives of $\ubold^h$, denoted by superscripts, are located on coloured bullets.
The transfer operators $\Wcal$ interpolate values from coloured bullets to blue squares and vice-versa.
}
\label{fig: staggered grid}
\end{figure}

We report in \cref{fig: staggered grid} the positions of $\grad^{q,h}$, up to order $q=3$, 
in order to illustrate how transfer operators $\Wcal$ work in interpolating the data on $\Gamma^h$.

\subsection{Anisotropic differential operators}
By construction, $\Mbold_1^h\in \RR^{|\Gamma^h|\times|\Gamma^h|\times|\Gamma^h|\times|\Gamma^h|}$ and $\grad^{1,h} \ubold^h\in \RR^{|\Xrm_{(1)}^{1,h}|\times|\Xrm_{(2)}^{1,h}|}$ so the grids $\Xrm_{(1)}^{1,h}, \Xrm_{(2)}^{1,h}$ have an ($h/2$)-offset with respect to $\Gamma^h$. 
In this case, locations of 
$\grad^{1,h} \ubold^h$ 
and 
$\Mbold_1^h$
are matched via the transfer operators $\Wcal^1=(\Wbold_{(1)}^1,\Wbold_{(2)}^1)$.

From \cref{rem: equiv overline D}, $\Mbold_1\grad \otimes u$ can be discretised in the correct grid position by the operator
\[
\Mbold_1^h\Wcal^1\grad^{1,h}: 
\RR^{|\Omega^h|}\to \RR^{|\Gamma^h|\times |\Gamma^h|},
\]
and the discretisation reads as
\[
\Mbold_1^h\Wcal^1\grad^{1,h} \otimes \ubold^h = 
\begin{pmatrix}
\bbold_1^h \Wcal^1 \grad^{1,h} \ubold^{h} \cdot \vbold^h \\
\bbold_2^h \Wcal^1\grad^{1,h} \ubold^{h} \cdot \vbold_\perp^h
\end{pmatrix}
=
\begin{pmatrix}
\bbold_1^h \left( \Wbold_{(1)}^1 \partial_1^h \ubold^h \vbold_1^h + \Wbold_{(2)}^1 \partial_2^h \ubold^h \vbold_2^h \right) \\
\bbold_2^h \left( \Wbold_{(2)}^1 \partial_2^h \ubold^h \vbold_1^h - \Wbold_{(1)}^1 \partial_1^h \ubold^h \vbold_2^h\right)
\end{pmatrix}.
\]
Therefore, the discrete weighted divergence
$\div_{\Mbold_1^h,\Wcal^{1}}^{h} \pbold^h:
\RR^{|\Gamma^h| \times |\Gamma^h|}\to \RR^{|\Omega^h|}$ 
is
\begin{equation}
\div_{\Mbold_1^h,\Wcal^1}^{h} \pbold^h 
=
\grad^{1,h}
\cdot 
\begin{pmatrix}
(\Wbold^{1}_{(1)})^\T \left( \bbold_1^h \pbold_1^h \vbold_1^h - \bbold_2^h \pbold_2^h \vbold_2^h \right) 
\\
(\Wbold^{1}_{(2)})^\T \left( \bbold_1^h \pbold_1^h \vbold_2^h + \bbold_2^h \pbold_2^h \vbold_1^h \right)
\end{pmatrix}.
\label{eq: div1 discrete}
\end{equation}
This leads to the discrete adjointness property, for every  
$\ubold^h\in\RR^{|\Omega^h|}, \pbold^h\in \RR^{|\Gamma^h|\times|\Gamma^h|}$:
\begin{equation}
\langle \Mbold_1^h \Wcal^1\grad^{1,h} \ubold^h,\, \pbold^h \rangle_{\Gamma^h}
= 
\langle \ubold^h,\,\div_{\Mbold_1^h,\Wcal^1}^{h}  \pbold^h \rangle_{\Omega^h},
\label{eq: discrete adjoint weighted derivative}
\end{equation}
where $\langle\blank,\,\blank\rangle_{\Gamma^h}: \RR^{|\Gamma^h|\times |\Gamma^h|} \times \RR^{|\Gamma^h|\times|\Gamma^h|}\to \RR$ and $\langle\blank,\,\blank\rangle_{\Omega^h}: \RR^{|\Omega^h|} \times \RR^{|\Omega^h|}\to \RR$. 

When considering higher order derivatives for a generic order $q$, the adjoint formula \cref{eq: discrete adjoint weighted derivative} is slightly more complicated due to the recursive definition of the weighted gradient and the location of the nested multiplication.  Indeed by \cref{eq: weighted Q gradient} we have for a fixed $q$
$
\grad^q_{\Mcal} u 
:= 
\Mbold_{q}\grad\otimes\dots\otimes\Mbold_{1}\grad\otimes u
$, whose finite-dimensional approximation is formally written as
$\grad^{q,h}_{\Mcal^h,\Wcal} \ubold^h\, : \, \RR^{|\Omega^h|} \to \RR^{|\Gamma^h|\times \dots \times|\Gamma^h|}$ via the recursion rule
\begin{equation}
\grad^{j,h}_{\Mcal^h,\Wcal} \ubold^h
:=
\begin{dcases}
\ubold^h
&\text{if } j=0,\\
(\Wcal^j)^\T \Mbold_{j}^h \Wcal^j \grad^{1,h} \otimes 
\grad^{j-1,h}_{\Mcal^h,\Wcal} \ubold^h
&
\text{if } j=1,\dots,q-1,\\
\Mbold_{q}^h \Wcal^q \grad^{1,h} \otimes 
\grad^{j-1,h}_{\Mcal^h,\Wcal} \ubold^h
&
\text{if } j=q.
\end{dcases}
\label{eq: higher order operator primal}
\end{equation}

The finite-dimensional approximation of the adjoint $\div_{\Mcal}^{q}$ is denoted with $\div_{\Mcal^h,\Wcal}^{q,h}
\,:\, 
\RR^{|\Gamma^h| \times \dots\times |\Gamma^h|} \to \RR^{|\Omega^h|}$ and defined via the recursion rule
\begingroup\makeatletter\def\f@size{10.5}\check@mathfonts
\def\maketag@@@#1{\hbox{\m@th\normalsize\normalfont#1}}
\begin{equation}
\div_{\Mcal^h,\Wcal}^{j,h} 
\Psibold^{h} 
:=
\begin{dcases}
\Psibold^{h}
&
\text{if } j=0,\\
\Wcal^{q-j}\grad^{1,h}\,\cdot\,(\Wcal^{q-j+1})^\T(\Mbold_{q-j+1}^h)^\T
\div_{\Mcal^h,\Wcal}^{j-1,h} 
\Psibold^{h}
&\text{if } j=1,\dots,q-1,
\\
\grad^{1,h}\,\cdot\,(\Wcal^{1})^\T(\Mbold_{1}^h)^\T
\div_{\Mcal^h,\Wcal}^{q-1,h} 
\Psibold^{h}
&
\text{if } j=q.
\end{dcases}
\label{eq: higher order operator dual}
\end{equation}
\endgroup

\begin{remark}\label{rem: trasop rem}
Note that in \cref{eq: higher order operator primal} for $j=q$ we omitted the inverse transfer operator $(\Wcal^q)^\T$ so as to force the highest derivative
$
\grad^{q,h}_{\Mcal^h,\Wcal} \ubold^h$ to be located in $\RR^{|\Gamma^h|\times \dots\times |\Gamma^h|}
$ and match with the position of $\Psibold^{h}$.
For the same reason, in \cref{eq: higher order operator dual} we omitted the transfer operator for $j=q$ so as to match the quantity with $\ubold^h$ in the grid $\Omega$. This operation is performed in view of the adjointness property stated below in \cref{eq: full adjointness}.
\end{remark}

For every $\ubold^h\in\RR^{|\Omega^h|}, \pbold^h\in \RR^{|\Gamma^h|\times \dots\times |\Gamma^h|}$, the discrete adjointness property holds:
\begin{equation}
\langle \grad_{\Mcal^h,\Wcal}^{q,h} \ubold^h,\, \pbold^h \rangle_{\Gamma^h}
= 
\langle \ubold^h,\, \div_{\Mcal^h,\Wcal}^{q,h}  \pbold^h \rangle_{\Omega^h},
\label{eq: full adjointness}
\end{equation}
where $\langle\blank,\,\blank\rangle_{\Gamma^h}: \RR^{|\Gamma^h| \times \dots\times |\Gamma^h|} \times \RR^{|\Gamma^h|\times \dots\times |\Gamma^h|}\to \RR$ and $\langle\blank,\,\blank\rangle_{\Omega^h}: \RR^{|\Omega^h|} \times \RR^{|\Omega^h|}\to \RR$. 
 
\section{Numerical optimisation}\label{sec: numerical optimization}
In what follows, we solve in the first instance the single line version of \cref{eq: minimization intro}, namely for a fixed $q>0$, a fixed $\alphabold=(\alpha_{j})_{j=0}^{q-1}$, a fixed collection of weighting matrices $\Mcal=(\Mbold_{j})_{j=1}^q$ and $\Scal$ the operator associated to the problem to solve, we aim to tackle the problem
\begin{equation}
\min_{\ubold\in\RR^{|\Omega^h|}}
\TDVM{q}{\alphabold}[\ubold,\Mcal]
+ 
\frac{\eta}{2}\norm{\Scal\ubold_0 - \ubold_0^{\diamond}}_2^2,
\label{eq: discretised problem fixed q}
\end{equation}
by means of a primal-dual hybrid gradient method \cite{ChaPoc2011,ChaPoc2016} and  following \cite[Equation 4.4]{BreKunPoc2010}.
With all discrete objects in place, we have
\[
\min_{\zbold_0^h \in \RR^{|\Omega^h|}} 
\TDVM{q,h}{\alphabold}[\zbold_0^h,\Mcal^h]
+ 
\frac{\eta}{2}\norm{\Scal\zbold_0^h - \zbold_0^{\diamond,h}}_2^2,
\] 
where $\zbold_0^{\diamond,h}=\ubold^{\diamond,h}\in\RR^{|\Omega^h|}$, $\div_{\Mcal^h,\Wcal}^{q,h}$ is the discretized weighted divergence w.r.t.\ the weights $\Mcal^h$ and the transfer operators $\Wcal$, $\TDVM{q,h}{\alphabold}$ is the discrete version of $\TDVM{q}{\alphabold}$ defined as
\[
\TDVM{q,h}{\alphabold}[\zbold_0^h,\Mcal^h] 
= 
\sup_{\Psibold^{h}\in\Ycal_{\Mcal^h,\alphabold}^{q,h}} 
\langle
\zbold_0^h,\,\div_{\Mcal^h,\Wcal}^{q,h}\Psibold^{h}
\rangle,
\]
and $\Ycal_{\Mcal^h,\alphabold}^{q,h}$ is the discretization of $\Ycal_{\Mcal,\alphabold}^q$ in \cref{eq: Ycal}, defined as
\begin{equation}
\Ycal_{\Mcal^h,\alphabold}^{q,h} 
= 
\left\{ \Psibold^{h} \in \RR^{2q|\Gamma^h|}\,\Big\lvert\, \norm{\div^{j,h}_{\Mcal^h,\Wcal} \Psibold^{h}}_\infty\leq \alpha_{j},\text{ for all }j=0,\dots,q-1\right\}.
\label{eq: test space discrete}
\end{equation}

\subsection{Discrete characterisation of TDV}
For a fixed $q>0$, the regulariser $\TDVM{q,h}{\alphabold}$ can be characterised as follows.
From the discrete version of $\TGV{\alphabold}{q}$ in \cite[Section 4.1]{Setzer11infimalconvolution} and following the characterization of $\TGV{\alphabold}{q}$ in \cite[Remark 3.8 and Remark 3.10]{BreKunPoc2010}, we can write the equivalent discrete definition of $\TDVM{q}{\alphabold}[u,\Mcal]$ for $\ubold^h\in \RR^{|\Omega^h|}$ and $\Mcal^h=(\Mcal_{j}^h)_{j=1}^q$ as
\begin{equation}
\TDVM{q,h}{\alphabold}[\ubold^h,\Mcal^h] = 
\inf_{
\substack{
\zbold_j^h \in \RR^{|\Xrm_{\iotabold}^{j,h}|}
\\ 
j=0,\dots,q, 
\\
\zbold_0=\ubold^h,\, \zbold_q^h=\bm{0}}
}
\sum_{j=1}^q \alpha_{q-j}\norm{
\Kcal^h_{j,j}
\zbold_{j-1}^h - \zbold_j^h}_{1},
\label{eq: equivalence TDV}
\end{equation}
where 
\begin{equation}
\Kcal^h_{j,j}
=
\begin{dcases}
(\Wcal^j)^\T\Mbold_{j}^h\Wcal^j\grad^{h}
&\text{if } j=1,\dots q-1,\\
\Mbold_{q}^h\Wcal^q\grad^{h} &\text{if } j=q.
\end{dcases}
\label{eq: Kcal discrete cases}
\end{equation}
Indeed, in the following let $j=1,...,q$, $\wbold_{q-j+1}^h\in\RR^{|\Xrm_\iotabold^{q-j+1}|}$ and $\zbold_{q-j}^h\in\RR^{|\Xrm_\iotabold^{q-j,h}|}$.
We call
\[
\DVM{j,h}{\alphabold}[\wbold_{q-j+1}^h,\Mcal^h] 
= 
\sup_{
\Psibold^h\in \Ycal_{\Mcal^h,\alphabold}^{q,h}
}
-\langle
 (\Wcal^{q-j+1})^\T 
 \Mbold_{q-j+1}^h \Wcal^{q-j+1} 
 \wbold_{q-j+1}^h,
 \,
 \div_{\Mcal^h,\Wcal}^{j-1}
 \Psibold^h
\rangle,
\]
where $\div_{\Mcal^{h},\Wcal}^{j-1}$ as in \cref{eq: higher order operator dual} 
and 
$\Ycal_{\Mcal^h,\alphabold}^{q,h}$ 
as in \cref{eq: test space discrete}.
Note that the sup is finite by definition of $\Ycal_{\Mcal^h,\alphabold}^{q,h}$ . Thus
$
\TDVM{q,h}{\alphabold}[\ubold^h,\Mcal^h]
=
\DVM{q,h}{\alphabold}(\grad^{h} \ubold^h,\Mcal^h)$ and we define
\[
K_\ell^h = 
\left\{
\Psibold^h\,:\,\Psibold^h\in \Ycal_{\Mcal^h,\alphabold}^{j,h},\, \norm{\div_{\Mcal^h,\Wcal}^{\ell,h}\Psibold^h}_\infty\leq \alpha_\ell \right\}
\quad\text{for every }\ell=0,\dots,q-1.
\]
With $
\sbold_{q,j}^h
= 
\Mbold_{q}^h \Wcal^q \grad^h 
\otimes 
\dots 
\otimes 
(\Wcal^{q-j+2})^\T
\Mbold_{q-j+2}^h
\Wcal^{q-j+2}
\grad^h$ and $(j-1)$-times integration by parts, the functional becomes
\[
\begin{aligned}
\DVM{j,h}{\alphabold}[\wbold_{q-j+1}^h,\Mcal^h] 
&= 
\left(\sum_{\ell=0}^{j-1}
\delta_{K_\ell^h}\right)^\ast 
\left(
(-1)^j 
\sbold_{q,j}^h \otimes 
(\Wcal^{q-j+1})^\T
\Mbold_{q-j+1}^h 
\Wcal^{q-j+1}
\wbold_{q-j+1}^h
\right),
\end{aligned}
\]
where $\delta_{K_\ell^h}$ is the characteristic function with values either $0$ or $+\infty$, and where the adjoint $\ast$ is the same as in \cite[Remark 3.9]{BreKunPoc2010}.
By Fenchel duality for the operator $\div_{\Mcal^h,\Wcal}^{j-1,h}$ we have:
\begingroup\makeatletter\def\f@size{9.5}\check@mathfonts
\def\maketag@@@#1{\hbox{\m@th\normalsize\normalfont#1}}
\[
\begin{aligned}
\DVM{j,h}{\alphabold}[\grad^h \zbold_{q-j}^h,\Mcal^h]
&=
\sup_{\Psibold^h \in\RR^{2^q\abs{\Gamma^h}}
} 
\left(
\begin{aligned}
-&
\langle
(\Wcal^{q-j+1})^\T\Mbold_{q-j+1}^h \Wcal^{q-j+1} \grad^h  \zbold_{q-j}^h
,\,
\div_{\Mcal^h,\Wcal}^{j-1,h}
\Psibold^h
\rangle
\\
&- \delta_{K_{j-1}^h}(\Psibold^h) - \sum_{\ell=0}^{j-2}\delta_{K_\ell}(\Psibold^h)
\end{aligned}
\right)
\\
&=
\sup_{\Psibold^h \in\RR^{2^q\abs{\Gamma^h}}
} 
\left(
\begin{aligned}
-&
\sum_{\ell=0}^{j-2}\delta_{K_\ell^h}(\Psibold^h)
-
\delta_{
\{
\norm{\blank}\leq\alpha_{j-1}
\}
}
(\div^{j-1,h}_{\Mcal^h,\Wcal} \Psibold^h) 
\\
-&
\langle 
(\Wcal^{q-j+1})^\T \Mbold_{q-j+1}^h \Wcal^{q-j+1} \zbold_{q-j}^h,\,\div^{j-1,h}_{\Mcal^h,\Wcal} \Psibold^h\rangle 
\end{aligned}
\right)
\\
&=
\inf_{ 
\zbold_{q-j+1}^h \in \RR^{\abs{\Xbold_{\iotabold}^{q-j+1,h}}}
} 
\left(
\begin{aligned}
\alpha_{j-1} &
\norm{(\Wcal^{q-j+1})^\T\Mbold_{q-j+1}^h\Wcal^{q-j+1}\grad^h \zbold_{q-j}^h- \zbold_{q-j+1}^h}_1 
\\
&+
\left(
\sum_{\ell=0}^{j-2}\delta_{K_\ell^h}\right)^\ast
\left(
(-1)^{j-1} 
\sbold_{q,j+1}^h \zbold_{q-j+1}^h
\right)
\end{aligned}
\right)
\\
&=
\inf_{ \zbold_{q-j+1}^h \in \RR^{\abs{\Xbold_{\iotabold}^{q-j+1,h}}}} 
\left(
\begin{aligned}
\alpha_{j-1}&
\norm{
(\Wcal^{q-j+1})^\T
\Mbold_{q-j+1}^h
\Wcal^{q-j+1}
\grad^h \zbold_{q-j}^h- \zbold_{q-j+1}^h}_1  \\
&+ \DVM{j-1,h}{\alphabold}[\grad^h \zbold_{q-j+1},\Mcal^h]
\end{aligned}
\right).
\end{aligned}
\]
\endgroup
Iterating the procedure for $j=q,\dots,2$ and by the identity
\[
\DVM{1,h}{\alphabold}(\grad^h \zbold_{q-1}^h) 
= 
\alpha_0 \norm{\Mbold_{q}^h\Wcal^q\grad^h\zbold_{q-1}^h}_1,
\]
we get
\begingroup\makeatletter\def\f@size{9.5}\check@mathfonts
\def\maketag@@@#1{\hbox{\m@th\normalsize\normalfont#1}}
\[
\DVM{q,h}{\alphabold}[\grad^h\zbold_0^h] 
= 
\inf_{
\substack{
\zbold_j^h\in \RR^{|\Xrm_\iotabold^{j,h}|}
\\ j=1,\dots,q-1\\ \zbold_0^h =\ubold^h
}
}
\left(
\sum_{j=1}^{q-1} \alpha_{q-j}\norm{(\Wcal^j)^\T\Mbold_{j}^h\Wcal^j\grad^h\zbold_{j-1}^h - \zbold_j^h}_1
\right)
+
\alpha_0\norm{\Mbold_q^h\Wcal^q \grad^h \zbold_{q-1}^h}_1,
\]
\endgroup
and thus, with $\Kcal^h_{j,j}$ as in \cref{eq: Kcal discrete cases}, we conclude
\[
\TDVM{q,h}{\alphabold}[\ubold^h,\Mcal^h] = 
\inf_{
\substack{
\zbold_j^h \in \RR^{|\Xrm_\iotabold^{j,h}|}
\\ j=1,\dots,q,\\\zbold_0^h=\ubold^h,\, \zbold_q^h=\bm{0}
}
}
\left(
\sum_{j=1}^q \alpha_{q-j}
\norm{\Kcal^h_{j,j}\zbold_{j-1}^h - \zbold_j^h}_1
\right).
\]
A continuous version of \cref{eq: equivalence TDV} also holds. This is proved in the second part of this work \cite{ParMasSch18analysis}.

The characterisation of $\TDVM{q,h}{\alphabold_q}$ in \cref{eq: equivalence TDV} is fundamental for writing a suitable primal-dual algorithm for the minimization problem in \cref{eq: minimization intro}.

\subsection{Discretised single minimization problems}
Let $\ubold^{\diamond,h}$ be a given discrete imaging data. 
For a fixed order $q>0$, let $\TDVM{q,h}{\alphabold}(\ubold^h,\Mcal^h)$ be decomposed as in \cref{eq: equivalence TDV},
with $\Mcal^h=(\Mbold_{j}^h)_{j=1}^q$, $\alphabold=(\alpha_{j})_{j=0}^{q-1}$ and the $\norm{\blank}_1$ denoted in the discrete setting by $\norm{\Zbold}_{2,1} = \sum_{i,j}\sqrt{\sum_{k=1}^s{(\Zbold_k)_{i,j}^2}}$ for a generic tensor-valued object $\Zbold=(\Zbold_k)_{k=1}^s$, with each $\Zbold_k\in\Tcal^2(\RR^{M\times N})$. 

The discrete single minimization problems, for $\zbold^h$ defined as in \cref{sec: discretised objects} read as:
\begin{itemize}
\item for order $q=1$, $\Mcal^h=(\Mbold_{1}^h)$, $\alphabold=\alpha_{0}$, $\zbold^h=(\zbold_0^h)$:
\begin{equation}
\min_{\zbold^h} 
\begin{aligned}
&
\left(\alpha_{0}
\norm{
\Mbold_{1}^h\Wcal^1\grad^h \zbold_0^h}_{2,1} + 
\frac{\eta}{2}\norm{\Scal \zbold_0^h -\zbold_{0}^{\diamond,h}}_2^2\right);
\end{aligned}
\label{eq: min order 1 discretized}
\end{equation}
\item for order $q=2$, $\Mcal_2^h=(\Mbold_{1}^h,\Mbold_{2}^h)$, $\alphabold=(\alpha_{0},\alpha_{1})$, $\zbold^h=(\zbold_0^h,\zbold_1^h)$:
\begin{equation}
\min_{\zbold^h} \left(
\begin{aligned}
&\alpha_{1} \norm{(\Wcal^1)^\T\Mbold_{1}^h\Wcal^1\grad^h \zbold_{0}^h - \zbold_1^h}_{2,1} 
\\
&+
\alpha_{0} \norm{
\Mbold_{2}^h\Wcal^2\grad^h \zbold_1^h}_{2,1}
+ 
\frac{\eta}{2}\norm{\Scal \zbold_0^h -\zbold_{0}^{\diamond,h}}_2^2
\end{aligned}
\right);
\label{eq: min order 2 discretized}
\end{equation}
\item for order $q=3$, $\Mcal_3^h=(\Mbold_{1}^h,\Mbold_{2}^h,\Mbold_{3}^h)$, $\alphabold_3=(\alpha_{0},\alpha_{1},\alpha_{2})$, $\zbold^h=(\zbold_0^h,\zbold_1^h,\zbold_2^h)$:
\begin{equation}
\min_{\zbold^h} \left(
\begin{aligned}
&
\alpha_{2} \norm{(\Wcal^1)^\T\Mbold_{1}^h\Wcal^1\grad^h \zbold_{0}^h - \zbold_1^h}_{2,1} 
\\
&+
\alpha_{1} \norm{(\Wcal^2)^\T\Mbold_{2}^h\Wcal^2\grad^h \zbold_{1}^h - \zbold_2^h}_{2,1} 
\\
&+
\alpha_{0} \norm{
\Mbold_{3}^h\Wcal^3\grad^h \zbold_2^h}_{2,1} 
+
\frac{\eta}{2}\norm{\Scal \zbold_0^h -\zbold_{0}^{\diamond,h}}_2^2
\end{aligned}
\right).
\label{eq: min order 3 discretized}
\end{equation}
\end{itemize}

For a fixed $q>0$, we aim to provide a more concise formulation of \cref{eq: min order 1 discretized,eq: min order 2 discretized,eq: min order 3 discretized}. 
Let $\zbold^h$ be as above
and $\Kcal^h=(\Kcal_{j,\ell}^h)_{j,\ell=1}^q$ be a matrix of operators associated to $\TDVM{q,h}{\alphabold}$ and defined as
\begin{equation}
\Kcal^h = 
\begin{pmatrix}
\Kcal^h_{1,1} & -\Ibold & \bm{0} & \dots & \dots & \dots & \bm{0} \\ 
\bm{0} & \Kcal^h_{2,2} & -\Ibold & \bm{0} &  &  & \vdots \\ 
\vdots & \ddots & \ddots & \ddots & \ddots &  & \vdots \\ 
\vdots &  & \bm{0} & \Kcal^h_{j,j} & -\Ibold & \bm{0} & \vdots \\ 
\vdots &  & & \ddots & \ddots & \ddots & \bm{0} \\ 
\vdots & &  &  & \bm{0} & \Kcal^h_{q-1,q-1} & -\Ibold \\ 
\bm{0} & \dots & \dots & \dots & \dots & \bm{0} & \Kcal^h_{q,q}
\end{pmatrix},
\label{eq: stack operators discrete}
\end{equation}
i.e.\ with $\Kcal^h_{j,\ell} = \bm{0}$ if $\ell \neq \{j, j+1\}$, $\Kcal^h_{j,j+1} = -\Ibold$ and as $\Kcal^h_{j,j}$ in \cref{eq: Kcal discrete cases} for each $j=1,\dots,q$.
Then, solving \cref{eq: discretised problem fixed q} is equivalent to solving for $\zbold^h=(\zbold_0^h,\dots,\zbold_{q-1}^h)$, with $\zbold_0^h = \ubold^h$, the problem:
\begin{equation}
\ubold^\star 
\in 
\argmin_{\ubold^h\in\RR^{|\Omega^h|}}
\inf_{
\zbold^h
} \,
\sum_{j=1}^{q} \alpha_{q-j}\norm{ \sum_{\ell=1}^q \Kcal^h_{j,\ell } \zbold_{\ell-1}^h }_{2,1} 
+ 
\frac{\eta}{2} \norm{\Scal \zbold_0^h-\ubold^{\diamond,h}}_2^2.
\label{eq: extended joint TDVM-L2 problem}
\end{equation}

By duality of the $\norm{\blank}_{2,1}$ norm and recalling that $\wbold^h$ is the dual vector defined in \cref{sec: discretised objects}, we rewrite
\cref{eq: extended joint TDVM-L2 problem} 
into a saddle-point minimization problem:
\begin{equation}
\min_{\zbold^h}
\max_{\wbold^h} \,
\langle \sum_{\ell=1}^q \Kcal^h_{j,\ell} \zbold_{\ell-1}^h, \wbold_{j}^h \rangle
-\sum_{j=1}^q\delta_{\left\{\norm{\blank}_{2,\infty}\leq \alpha_{q-j}\right\}}(\wbold_{j}^h)
+\frac{\eta}{2}\norm{\Scal \zbold_0^h - \ubold^{\diamond,h}}_2^2,
\label{eq: joint general saddle index discretised}
\end{equation}
or, in short notation:
\begin{equation}
\min_{\zbold^h}
\max_{\wbold^h} 
\,
\langle \Kcal^h\zbold^h,\,\wbold^h\rangle - F^\ast(\wbold^h) + G(\zbold^h),
\label{eq: joint general saddle discretised}
\end{equation}
where $G(\zbold^h)$ is a partially strongly convex term, since it can be seen as $ G(\zbold^h)= G_0(\Pcal \zbold^h)$ for the projection $\Pcal$ of $\zbold^h$ onto the subspace of $\zbold_0^h$ and with $G_0$ being strongly convex, and where
\[
F^\ast(\wbold^h) 
= \sum_{j=1}^q
 \delta_{
 \left\{
 \norm{\blank}_{2,\infty}
 \leq {\alpha_{q-j
 }}
 \right\}
 }(\wbold_{j}^h),
\]
and
\[
G(\zbold^h) = 
G_0(\Pcal\zbold^h) =
\frac{\eta}{2}\norm{\Pcal\begin{pmatrix} \zbold_0^h\\\zbold_1^h\\\vdots\\\zbold_q^h\end{pmatrix} - \begin{pmatrix}\ubold^{\diamond,h}\\\bm{0}\\\vdots\\\bm{0}\end{pmatrix}}_2^2
\text{ with }
\Pcal=
\begin{pmatrix}
\Scal  & \bm{0} & \dots  & \bm{0}\\
\bm{0} & \bm{0} & \dots  & \bm{0}\\
\vdots  & \vdots  &\ddots  &\bm{0}\\
\bm{0} & \bm{0} &\dots  & \bm{0}
\end{pmatrix},
\]
and $\Pcal\zbold^h = \left(\Scal\zbold_0^h, \bm{0},  \cdots, \bm{0}\right)^\T$.

\subsection{Proximal operators}
We aim to solve the saddle point problem
\cref{eq: joint general saddle discretised}
with a Primal-Dual Hybrid Gradient (PDHG) algorithm. We need the proximal operators of $F^\ast$ and $G$.

The proximal map of $F^\ast$ evaluated at a point $\widehat{\wbold}^h$ is the sum of the projections onto the respective polar balls since $F^\ast$ is fully separable\footnote{meaning that $F^\ast$ is a function that can be written as a sum of functions in disjoint sets of variables.},
\[
\prox_{\sigma F^\ast}(\widehat{\wbold}^h) = \sum_{j=1}^q \prox_{\sigma F^\ast}(\widehat{\wbold}_{j}^h),
\quad\text{with}\quad
\prox_{\sigma F^\ast}(\widehat{\wbold}_{j}^h) =
\frac{\widehat{\wbold}_{j}^h}{\max\left(1, \alpha_{q-j}^{-1}\norm{\widehat{\wbold}_{j}^h}_2\right)}.
\]

The proximal map of $G$ should be evaluated at a point $\widehat{\zbold}^h$. 
Recalling that $\zbold_0^{h}=\ubold^{h}$ and that $\zbold$ is as in \cref{sec: discretised objects}, we have: 
\begin{equation}
\begin{aligned}
\prox_{\tau G}(\widehat{\zbold}^h) 
&=
\argmin_{\zbold^h}
\left(
\frac{\eta}{2}
\norm{
\begin{pmatrix}
\Scal \zbold_0^{h} - \ubold^{\diamond,h}\\
\bm{0}\\\vdots\\\bm{0}\end{pmatrix}
}_2^2 
+ 
\frac{1}{2\tau} \norm{\begin{pmatrix} \zbold_0^h\\\zbold_1^h\\\vdots\\\zbold_q^h\end{pmatrix}-\begin{pmatrix}\widehat{\zbold}_0^h\\\widehat{\zbold}_1^h\\\vdots\\\widehat{\zbold}_q^h\end{pmatrix}}_2^2
\right)
\\
&=
\begin{pmatrix}\displaystyle
\argmin_{\zbold_0^h\in\RR^{|\Omega^h|}}
\left(
\frac{\eta}{2}\norm{\Scal \zbold_0^h-\ubold^{\diamond,h}}_2^2 
+ 
\frac{1}{2\tau} \norm{\zbold_0^h-\widehat{\zbold}_0^h}_2^2
\right)
\\
\widehat{\zbold}_1^h\\
\vdots
\\
\widehat{\zbold}_q^h\\
\end{pmatrix}.
\label{eq: prox G}
\end{aligned}
\end{equation}
Let us focus on the first component of $\prox_{\tau G}$, in \cref{eq: prox G}: we have to solve
\[
\argmin_{\zbold_0^h\in\RR^{|\Omega^h|}} 
\frac{\eta}{2}\norm{\Scal \zbold_0^h-\ubold^{\diamond,h}}_2^2 
+ 
\frac{1}{2\tau} \norm{\zbold_0^h-\widehat{\zbold}_0^h}_2^2,
\]
whose minimum is achieved by a $\zbold_0^h$ that solves, for $\Scal^\ast$ adjoint of $\Scal$,
\[
\begin{aligned}
\tau\eta \Scal^\ast ( \Scal \zbold_0^h-\ubold^{\diamond,h}) + (\zbold_0^h-\widehat{\zbold}_0^h) &= 0
\\
(\Ibold + \tau\eta\Scal^\ast\Scal)\zbold_0^h-\tau\eta\Scal^\ast \ubold^{\diamond,h} -\widehat{\zbold}_0^h &= 0
\\
(\Ibold + \tau\eta\Scal^\ast\Scal)\zbold_0^h-\tau\eta\Scal^\ast \ubold^{\diamond,h} -\widehat{\zbold}_0^h  +\tau\eta\Scal^\ast\Scal \widehat{\zbold}_0^h-\tau\eta\Scal^\ast\Scal \widehat{\zbold}_0^h &= 0
\\
(\Ibold + \tau\eta\Scal^\ast\Scal) (\zbold_0^h - 
 \widehat{\zbold}_0^h )  &= \tau\eta\Scal^\ast (\ubold^{\diamond,h} - \Scal \widehat{\zbold}_0^h).
\end{aligned}
\]
Thus, the first component of $\prox_{\tau G}(\widehat{\zbold}^h)$, that is  $\prox_{\tau G_0}(\Pcal\widehat{\zbold}^h)$, is
\[
\zbold_0^h
=
\widehat{\zbold}_0^h  + (\Ibold + \tau\eta\Scal^\ast\Scal)^{-1}\tau\eta\Scal^\ast (\ubold^{\diamond,h} -\Scal \widehat{\zbold}_0^h).
\]
Note that for the Rudin-Osher-Fatemi denoising problem we have $\Scal=\Ibold$, thus $\Scal^\ast=\Ibold$, and the proximal map agrees with the one computed in \cite[pag.\ 133]{ChaPoc2011}.

\subsection{Operator norm}\label{sec: operator norm}
Following the approach in \cite[Section 4]{BreKunPoc2010} and \cite[Section 6.1]{ChaPoc2011}, we estimate a bound on the norm of the linear operator $\grad_{\Mcal^h,\Wcal}^{q,h}$ in \cref{eq: higher order operator primal} in view of the implementation of a suitable primal-dual algorithm. 
Let $\norm{\blank}$ be the operator norm, then for each $q=1,\dots,\Qrm$ we have
\begin{equation}
L_q^2
=
\norm{\grad_{\Mcal^h,\Wcal}^{q,h}}^2 
= 
\norm{\div_{\Mcal^h,\Wcal}^{q,h}}^2 
= 
\sup_{\Psibold^h\in
\Ycal_{\Mcal^h,\alphabold}^{q,h}}
\frac{\norm{\div_{\Mcal^h,\Wcal}^{q,h}\Psibold^h}_2^2}{\norm{\Psibold^h}_2^2}.
\label{eq: operator norm}
\end{equation}
In the two-dimensional setting, when $q=1$ then $\div_{\Mcal^h,\Wcal}^{1,h}\Psibold^h$ in \cref{eq: operator norm} reduces to
\[
\div_{\Mcal^h,\Wcal}^{1,h}\Psibold^h 
= \div^{1,h}( (\Mbold_{1}^h \Wcal^1)^\T\Psibold^h),
\]
and by applying the finite difference scheme in \cref{eq: div1 discrete}, from $\norm{\grad}\leq \sqrt{8}h^{-1}$ we estimate:
\[
\begin{aligned}
L_1^2
=
\norm{\div_{\Mcal^h,\Wcal}^{1,h}\Psibold^h}^2
&=
\div_{\Mcal^h,\Wcal}^{1,h}\Psibold^h\cdot\div_{\Mcal^h,\Wcal}^{1,h}\Psibold^h
\\ 
&=\div( (\Mbold_{1}^h\Wcal^1)^\T\Psibold^h)\cdot\div( (\Mbold_{1}^h\Wcal^1)^\T\Psibold^h)\\
&=
\Mbold_{1}^h\Wcal^{1}\grad\otimes\div \left((\Mbold_{1}^h\Wcal^{1})^\T\Psibold^h \right) \cdot \Psibold^h\leq 
\frac{8}{h^2} \norm{\Mbold_{1}^h\Wcal^{1}}_F^2\norm{\Psibold^h}_2^2.
\end{aligned}
\]
For a fixed $q$, since it holds
\[
\norm{(\Wcal^j)^\T\Mbold_{j}^h\Wcal^{j}}_F = \norm{( (\Wcal^j)^\T \Mbold_{j}^h \Wcal^{j})^\T}_F,\quad\text{ for each }j=1,\dots,q,
\]
then the operator norm $L_q$ is estimated via
\begin{equation}
L_q^2
\leq 
(8 h^{-2})^q \prod_{j=1}^q \norm{ (\Wcal^j)^\T\Mbold_{j}^h\Wcal^j}^2_2.
\label{eq: operator norm q}
\end{equation}
\begin{remark}
Since $\Wcal$ is made by partition of unit transfer operators and $\norm{\Mbold_{j}^h}_2\leq 1$ by construction, we can estimate the right-hand side of \cref{eq: operator norm q} as:
\[
L_q^2 \leq (8 h^{-2})^q,
\]
which agrees with the classic isotropic setting given by the choices $\Mbold_{j}^h=\Ibold$ without the use of $\Wcal^j$ for every $j=1,\dots,q$. Indeed, we have $L^2\leq 8h^{-2}$ for $\TV{}{}$ and $L^2\leq 64h^{-4}$ for $\TGV{\alphabold}{2}$.
\end{remark}

\subsection{Primal-Dual Hybrid Gradient algorithm}\label{sec: primal-dual section}
Now we are ready for solving  \cref{eq: joint general saddle discretised} with a Primal-Dual Hybrid Gradient (PDHG) algorithm following \cite{ChaPoc2011}.
Let $q$ be fixed and $L_q=\norm{\Kcal}$ be the operator norm as in \cref{sec: operator norm}, i.e.\
$
L_q:=\sup_\zbold\left\{ \norm{\Kcal\zbold}_2\text{ s.t.\ }\norm{\zbold}_2\leq 1\right\},
$
and let $\tau,\sigma>0$, $\omega\in[0,1]$ be fixed, such that $\tau\sigma L_q^2<1$. 
Then the PDHG algorithm \cite[Algorithm 1]{ChaPoc2011} reads as the iteration of
\begin{equation}
\begin{aligned}
\wbold^{n+1,h} &= \prox_{\sigma F^\ast}\left( \wbold^{n,h} + \sigma \Kcal^h\overline{\zbold}^{n,h}\right);\\
\zbold^{n+1,h} &= \prox_{\tau G} \left(\zbold^{n,h} - \tau(\Kcal^h)^\ast\wbold^{n+1,h} \right);
\\
\overline{\zbold}^{n+1,h} &= \zbold^{n+1,h} + \omega (\zbold^{n+1,h} - \zbold^{n,h}),
\end{aligned}
\label{eq: primal-dual accelerated}
\end{equation}
where we denoted with an index $n$ the iterations, starting from admissible $\zbold^{0,h}$ and $\wbold^{0,h}$.
The final solution is achieved by $\ubold^\star=\zbold_0^{n+1,h}$. Compare \cref{alg: main primal-dual} for details.

\begin{algorithm}[tbh]
\caption{PDHG algorithm for the minimization model in Equation \cref{eq: joint general saddle discretised}}
\SetAlgoLined\small
\label{alg: main primal-dual}
\SetKwProg{Fn}{Function}{:}{}
\SetKwFunction{PrimalDual}{PrimalDual}
\SetKwFunction{computev}{compute\_vector\_field}
\SetKwFunction{primaldual}{primal\_dual}
\SetKwFunction{computeopnorm}{compute\_operator\_norm}
\SetKwData{maxitert}{$T$}
\SetKwData{maxiterk}{\texttt{maxiter}}
\SetKwData{F}{F}
\SetKwData{G}{G}
\SetKwData{K}{K}
\SetKwData{KS}{KS}
\SetKwData{ProxF}{ProxF}
\SetKwData{ProxFS}{ProxFS}
\SetKwData{ProxG}{ProxG}
\SetKwInOut{Input}{Input}
\SetKwInOut{ParametersModel}{Parameters for the model}
\SetKwInOut{ParametersV}{Parameters for $\vbold$}
\SetKwInOut{ParametersPD}{Parameters for primal-dual}
\SetKwInOut{Operators}{Operators needed}
\SetKwInOut{Initialization}{Variables to initialize}
\SetKwInOut{Update}{Update}
\BlankLine
\ParametersModel{$q>0$, $\eta >0$ and $\alphabold=(\alpha_{j})_{j=0}^{q-1}$.}
\Operators{$\Kcal^h$, $\Scal$.}
\ParametersPD{$\sigma_0$, $\tau_0$, $\omega_0$ such that $\sigma_0\tau_0 L_q^2<1$ for $L_q=\norm{\Kcal^h}$, $\gamma>0$.}
\Initialization{$\zbold^{0,h} = (\zbold_{j-1}^{0,h})_{j=1}^\Qrm$, $\wbold^{0,h} = (\wbold_{j}^{0,h})_{j=1}^q$.}
\BlankLine
\Fn{\PrimalDual{$\ubold^{\diamond,h},\Kcal^h,\alphabold,\eta$}}{
\BlankLine
\While{\texttt{stopping criterion is not satisfied}}
{
\BlankLine
\tcp{Dual problem}
\BlankLine
\For{$j=1,\dots,{q-1}$}{
\BlankLine
$\wbold_{j}^{n,h} 
 \leftarrow
\wbold_{j}^{n,h} + \sigma_n \left( \Kcal^h_{j,j} \overline{\zbold}_{j-1}^{n,h} - \overline{\zbold}_{j}^{n,h} \right)$\tcp*{update the dual (inner)}
$\displaystyle\wbold_{j}^{n+1,h} = \wbold_{j}^{n,h}/\max\left(1,\alpha_{q-j}^{-1}\norm{\wbold_{j}^{n,h}}_2\right)$\tcp*{proximal operator of $F^\ast$ (inner)}
\BlankLine
}
$\wbold_{q}^{n,h}
 \leftarrow 
\wbold_{q}^{n,h} + \sigma_n \Kcal^h_{q,q} \overline{\zbold}_{q-1}^{n,h}$\tcp*{update the dual (outer)}
\BlankLine
$\displaystyle\wbold_{q}^{n+1,h} = \wbold_{q}^{n,h}/\max\left(1,\alpha_{0}^{-1}\norm{\wbold_{q}^{n,h}}_2\right)$\tcp*{proximal operator of $F^\ast$ (outer)}
\BlankLine
\tcp{Primal problem}
\For{$j=2,\dots,{q}$}{
$\zbold_{j-1}^{n,h}  \leftarrow \displaystyle\zbold_{j-1}^{n,h}  - \tau_n (\Kcal^h_{j,j})^\ast \wbold_{j}^{n+1,h} -\wbold_{j-1}^{n+1,h} $\tcp*{update the primal (inner)}
$\zbold_{j-1}^{n+1,h} = \zbold_{j-1}^{n,h}$\tcp*{proximal operator of $G$ (inner)}
}
$\displaystyle 
\zbold_0^{n,h}  \leftarrow \zbold_0^{n,h} - \tau_n  \left( (\Kcal^h_{1,1})^\ast \wbold_{1}^{n+1,h}\right) $\tcp*{update the primal (outer)}
$\zbold_0^{n+1,h} = \zbold_0^{n,h}  + (\Ibold + \tau\eta\Scal^\ast\Scal)^{-1}\tau\eta\Scal^\ast (\ubold^{\diamond,h} -\Scal \zbold_0^{n,h})$\tcp*{proximal operator of $G$ (outer)}
\BlankLine
\tcp{Parameters update}
$\omega_{n}=\omega_0$;
$\tau_{n+1}=\tau_0$; 
$\sigma_{n+1}=\sigma_0$
\tcp*{or update $\omega_n,\tau_n,\sigma_n$ if acceleration is possible}
\BlankLine
\tcp{Extrapolation step}
$\overline{\zbold}^{n+1,h} = \zbold^{n+1,h} + \omega_n(\zbold^{n+1,h}-\zbold^{n,h})$;
}
}
\Return{$\ubold^\star=\zbold_0^{n+1,h}$.}
\BlankLine
\end{algorithm}

\paragraph{Acceleration}
If $\Scal=\Ibold$ and only the first order regulariser is involved ($\Qrm=1$) then the fidelity term $G$ is strongly convex with convexity parameter $\eta$ (since it does not involve the terms of $\zbold^h$ related to the derivatives of order greater than 1) and the dual problem is smooth. 
Therefore, it is possible to accelerate the PDHG algorithm with \cite[Algorithm 2]{ChaPoc2011}: we can take $\tau_0\sigma_0 L_q^2\leq1$ and update $\tau_n,\sigma_n,\omega_n$ by taking $L_G$ as the Lipschitz constant of $G$, $\tau_0=L_q^{-1}$ and $\gamma=0.5 L_G^{-1}$, with the update rule in \cref{eq: primal-dual accelerated} before $\overline{\zbold}^{n+1,h}$ reading as
\[
\omega_n = \left(1+ 2\gamma\tau_n \right)^{-0.5};\, \tau_{n+1} = \omega_n\tau_n;\,\sigma_{n+1} = \sigma_n\omega_n^{-1}.
\]
When $\Scal=\Ibold$ and $\Qrm>1$ then $G$ is only partially strongly convex and one can use either \cite[Algorithm 1]{ChaPoc2011} or the acceleration proposed in \cite{Valkonen2017}.
In any case, when $\Scal$ makes $G$ not strongly convex then the use of \cite[Algorithm 1]{ChaPoc2011} is recommended: in such case, $\sigma_n$ and $\tau_n$ are fixed a-priori, e.g.\ the authors in \cite{KnoBrePocSto2011} adopted the parameters $\sigma_n=\tau_n=1/\sqrt{12}$ for the second-order regulariser $\TGV{}{2}$ and for a grid of grid-size $h=1$.

\subsection{Primal-Dual Gap}
As exit condition for the primal-dual algorithm of the $\TDVM{q,h}{\alphabold}-\LL{2}$ minimisation problem, it is possible to either define a maximum number of iterations to reach or to impose a threshold for the primal-dual gap, defined for the current solutions $\overline{\zbold}$ and $\overline{\wbold}$ as
\begin{equation}
\begin{aligned}
\texttt{pdhg\_gap}(\overline{\zbold},\overline{\wbold}) = 
&\,\TDVM{q,h}{\alphabold}[\overline{\zbold}_0,\Mcal^h] + \frac{\eta}{2}
\norm{
\overline{\zbold}_0 - \ubold^{\diamond,h}}_2^2 \\
&+ \sum_{j=1}^q 
\delta_{\left\{\norm{\blank}_{2,\infty}\leq \alpha_{q-j}\right\}} (\overline{\wbold}_{j}) 
+ \frac{\eta}{2}\norm{ (\Kcal^h)^\ast \overline{\wbold}}_2^2
- \eta \langle (\Kcal^h)^\ast \overline{\wbold},\, \begin{pmatrix} \zbold_0^\diamond\\\bm{0}\\\vdots\\\bm{0} \end{pmatrix}\rangle.
\end{aligned}
\end{equation}

\subsection{Joint minimisation problem}\label{sec: joint description}
In the continuous setting and for any fixed $q>0$, the particular choice of $\Mcal = (\Ibold,\dots,\Ibold,\Mbold)$ and $\alphabold=(\alpha_{q,0},+\infty,\dots,+\infty)$ deserves a separate discussion. 
Indeed, the decomposition of $\TDVM{q}{\alphabold}$ in the continuous setting results as 
\[
\TDVM{q}{\alphabold}[u,\Mcal]=\alpha_0 \norm{\Mbold\grad \zbold_{q-1}}_{2,1}+\sum_{j=1}^{q-1} +\infty \norm{\grad \zbold_{j-1} - \zbold_j } _{2,1},
\]
in which the sparsity of the inner order of derivatives would not be fully exploited due to the weight $+\infty$. The above decomposition is equivalent to impose $\zbold_j = \grad\zbold_{j-1}$ for any $j=1,\dots,q-1$, resulting in $\TDVM{q}{\alphabold}[u,\Mcal]\equiv \norm{\Mbold\grad^q u}_{2,1}$. As discussed in \cref{rem: trasop rem}, in this case the usage of the transfer operators in the discrete setting and for the inner order of derivates is useless and we can therefore discretise the above regulariser as
\[
\TDVM{q,h}{\alphabold}[\ubold^h,\Mcal^h] \equiv \alpha_0\norm{\Mbold^h\Wcal^q \grad^q \ubold^h}_{2,1}.
\]

In our applications, we will mainly focus on the effect of weighting the highest order derivative by means of taking the exact inner derivatives, but combined jointly with regularisers $\TDVM{q}{\alphabold_q}$ for different $q>0$, i.e.\ we aim to solve in the discrete setting
\begin{equation}
\ubold^\star
\in\argmin_{\ubold^h\in\RR^{\abs{\Omega^h}}}
\sum_{q=1}^\Qrm\TDVM{q,h}{\alphabold_q}[\ubold^h,\Mcal_q^h] + \frac{\eta}{2} \norm{\Scal \ubold^h -\ubold^{\diamond,h}}_2^2,
\label{eq: reduced energy simplified}
\end{equation}
with $\Mcal_q^h=(\Ibold,\dots,\Ibold,\Mbold_q^h)$ and $\alphabold_q = (\alpha_{q,0},+\infty,\dots,\infty)$ for each $q=1,\dots,\Qrm$, leading to
\[
\ubold^\star\in
\argmin_{\ubold^h\in\RR^{\abs{\Omega^h}}}
 \left(
\sum_{q=1}^\Qrm
\alpha_{q,0} 
\norm{\Mbold_q^h\Wcal^q\grad\otimes\grad^{q-1}\ubold^h}_{2,1} + \frac{\eta}{2} \norm{\Scal \ubold^h -\ubold^{\diamond,h}}_2^2
\right).
\]
Then it is possible to reduce \cref{alg: main primal-dual} to \cref{alg: reduced main primal-dual}, where   the accelerated PDHG \cite[Algorithm 2]{ChaPoc2011} can be used for any choice of $\Scal$ that makes $G$ strongly convex.

\begin{algorithm}[tbh]
\caption{PDHG algorithm for the joint minimization reduced model in Equation \cref{eq: reduced energy simplified}}
\SetAlgoLined\small
\label{alg: reduced main primal-dual}
\SetKwProg{Fn}{Function}{:}{}
\SetKwFunction{PrimalDual}{PrimalDual\_reduced}
\SetKwFunction{computev}{compute\_vector\_field}
\SetKwFunction{primaldual}{primal\_dual}
\SetKwFunction{computeopnorm}{compute\_operator\_norm}
\SetKwData{maxitert}{$T$}
\SetKwData{maxiterk}{\texttt{maxiter}}
\SetKwData{F}{F}
\SetKwData{G}{G}
\SetKwData{K}{K}
\SetKwData{KS}{KS}
\SetKwData{ProxF}{ProxF}
\SetKwData{ProxFS}{ProxFS}
\SetKwData{ProxG}{ProxG}
\SetKwInOut{Input}{Input}
\SetKwInOut{ParametersModel}{Parameters for the model}
\SetKwInOut{ParametersV}{Parameters for $\vbold$}
\SetKwInOut{ParametersPD}{Parameters for primal-dual}
\SetKwInOut{Operators}{Operators needed}
\SetKwInOut{Initialization}{Variables to initialize}
\SetKwInOut{Update}{Update}
\BlankLine
\ParametersModel{$\Qrm>0$, $\eta >0$  and $\abold=(\alpha_{q,0})_{q=1}^\Qrm$.}
\Operators{$\Kcal^h = (\Kcal_q^h)_{q=1}^\Qrm$ where $\Kcal_q^h=\Mbold \Wcal^q\grad^{q,h}$, $\Scal$.}
\ParametersPD{$\sigma_0$, $\tau_0$, $\omega_0$ such that $\sigma_0\tau_0 L_q^2<1$ for $L_q=\max_q\norm{\Kcal_q^h}$.} \Initialization{$\ubold^{0,h}$, $\wbold^{0,h} = (\wbold_{q,j}^{0,h})_{q,j=1}^\Qrm$.}
\BlankLine
\Fn{\PrimalDual{$\ubold^{\diamond,h},\Kcal^h,\abold,\eta$}}{
\BlankLine
\While{\texttt{stopping criterion is not satisfied}}
{
\BlankLine
\tcp{Dual problem}
\For{$q=1,\dots,\Qrm$}{
$\wbold_{q}^{n,h}\displaystyle
 \leftarrow 
\wbold_{q}^{n,h} + \sigma_n   \Kcal_q^h \overline{\ubold}^{n,h}$\tcp*{update the dual}
$\displaystyle\wbold_{q}^{n+1,h} = \wbold_{q}^{n,h}/\max\left(1,\alpha_{q,0}^{-1}\norm{\wbold_{q}^{n,h}}_2\right)$\tcp*{proximal operator of $F^\ast$}
}
\BlankLine
\tcp{Primal problem}
$\displaystyle 
\ubold^{n,h}  \leftarrow \ubold^{n,h} - \tau_n \left(\sum_{q=1}^\Qrm\left( \Kcal_q^h \right)^\ast \wbold_{q}^{n+1,h} \right)$\tcp*{update the primal}
\BlankLine
$\ubold^{n+1,h} = \ubold^{n,h}  + (\Ibold + \tau\eta\Scal^\ast\Scal)^{-1}\tau\eta\Scal^\ast (\ubold^{\diamond,h} -\Scal \ubold^{n,h})$\tcp*{proximal operator of $G$}
\BlankLine
\tcp{Parameters update}
$\omega_{n}=\omega_0$;
$\tau_{n+1}=\tau_0$; 
$\sigma_{n+1}=\sigma_0$
\tcp*{or update $\omega_n,\tau_n,\sigma_n$ if acceleration is possible}
\BlankLine
\tcp{Extrapolation step}
$\overline{\ubold}^{n+1,h} = \ubold^{n+1,h} + \omega_n(\ubold^{n+1,h}-\ubold^{n,h})$;
}
}
\Return{$\ubold^\star=\ubold^{n+1,h}$.}
\BlankLine
\end{algorithm}

\section{Imaging Denoising}\label{sec: model-description-denoising}
In what follows we demonstrate the performance of
the introduced regulariser $\TDVM{q}{\alphabold_q}$ for the applications of image denoising.
We focus in particular on the cases $q=1,2,3$ for the single \cref{eq: discretised problem fixed q} and joint minimisation model \cref{eq: reduced energy simplified}.
Results are computed on a standard laptop (MATLAB R2019a, MacBook Pro 13'', 2.9 GHz Intel Core i5, 8 GB 1867 MHz DDR3). The code is freely available at the authors' webpage\footnote{Code freely available at \url{https://github.com/simoneparisotto}}.

Let $\Omega\subset\RR^2$, $u:\Omega\to\RR$ be a grey-scale image (colour images $\ubold:\Omega\to\RR^3$ are considered one colour channel at time), $\vbold:\Omega\to\RR^2$ a field with $\norm{\vbold}_2=1$ and $u^\diamond$ a given noisy image.  

\subsection{Estimation of vector field \texorpdfstring{$\vbold$}{v}}\label{sec: estimate v}
For estimating $\vbold$ we use the following strategy.
Let $\sigma,\rho>0$. Let $\lambda_1(\xbold),\lambda_2(\xbold)$ be such that $\lambda_1(\xbold)\geq\lambda_2(\xbold)$, the ordered eigenvalues of
\begin{equation}
\Jbold_\rho(u) := K_\rho \ast (\grad u_\sigma \otimes \grad u_\sigma)\quad\text{for } u_\sigma :=K_\sigma \ast u,
\label{eq: structure tensor}
\end{equation}
and $\ebold_1,\ebold_2\in\RR^2$ the associated eigenvectors. 
Let $\widetilde{\vbold}(\xbold)=\ebold_2(\xbold)$ be the local direction of the anisotropy, corresponding to an approximation of $(\grad^\perp u)/\norm{\grad^\perp u}$ . In order to compute a vector field smoother than $\widetilde{\vbold}$, we adopt a further regularisation step, similarly as in \cite{LelMorSch2013}. 

Let $w(\xbold)\in[0,1]$. We aim to smooth the vector field where the anisotropy weight $w(\xbold)$ is close to 0 while keeping the already computed vector field in regions with strong anisotropy. This is equivalent to solving the following problem:
\begin{equation}
\vbold 
= 
\argmin_{\tbold \in\RR^{|\Gamma^h|}}
\frac{1}{2} \int_\Omega w(\xbold) \norm{\tbold(\xbold)-\widetilde{\vbold}(\xbold)}_2^2\diff\xbold + \frac{\gamma}{2}\int_\Omega \norm{\grad \tbold(\xbold)}_2^2\diff \xbold.
\label{eq: find v denoising}
\end{equation}
We use the local estimation of the anisotropy as weights $w(\xbold)$, for $\varepsilon>0$:
\begin{equation}
w(\xbold) = \frac{\lambda_1(\xbold)-\lambda_2(\xbold)}{\lambda_1(\xbold) + \lambda_2(\xbold)+\varepsilon}.
\label{eq: denoising weight 2}
\end{equation} 
We can use $w(\xbold)$ also to vary locally $\bbold=(1,\beta(\xbold))$: we have already seen that the process is more isotropic as $\beta(\xbold)$ is closer to 1. For this reason, a possible strategy to vary $\beta(\xbold)$ is:
\begin{itemize}
\item first, estimate the anisotropy (values close to 1 correspond to isotropic regions) by:
\begin{equation}
\operatorname{anisotropy}(\xbold) =  1-w(\xbold),
\label{eq: compute b vary step 1}
\end{equation}
\item second, rescale in $[0,1]$ to define $\beta$:
\begin{equation}
\beta(\xbold) = \frac{\operatorname{anisotropy}(\xbold)-\min \operatorname{anisotropy}(\xbold)}{\max \operatorname{anisotropy}(\xbold) - \min \operatorname{anisotropy}(\xbold)}.
\label{eq: compute b vary step 2}
\end{equation}
\end{itemize}
With this strategy, the higher the image anisotropy the closer $\bbold$ is to $(1,0)$: in such cases strong directional structures are emphasised by our directional regulariser. Conversely, when $\bbold=(1,1)$, isotropic smoothing is performed in flat regions. 
In the following, we will require that $\bbold=(1,\beta(\xbold))^\T$ for every $\xbold\in\Omega$, thus
\[
\Mbold\grad\otimes u(\xbold)=\Lambdabold_\bbold(\Rbold_\thetabold)^\T\grad \otimes u(\xbold)= \begin{pmatrix}
\grad_\vbold u(\xbold)\\ \beta(\xbold)\grad_{\vbold_\perp} u(\xbold)
\end{pmatrix},
\]
where $\vbold$ is a given estimated field.
Also, we may refine $\vbold$ by updating the
parameters $(\sigma,\rho)$ so as to restart the denoising problem with a better estimation of the vector field.

\subsection{Single minimisation model}\label{sec: smm}
Here we describe results with the single model $\TDVM{q}{\alphabold}-\LL{2}$ and fixed $q$ to be chosen between $1,2$ and $3$. We will explore all the possible combinations for $(\Mbold_1,\dots,\Mbold_q)$, where each $\Mbold_j$ will be either the identity $\Ibold$ or the anisotropic $\Mbold$ weights.

\subsubsection*{Numerical results}
\begin{figure}
\centering
\begin{subfigure}[t]{0.235\textwidth}\centering
\includegraphics[width=1\textwidth]{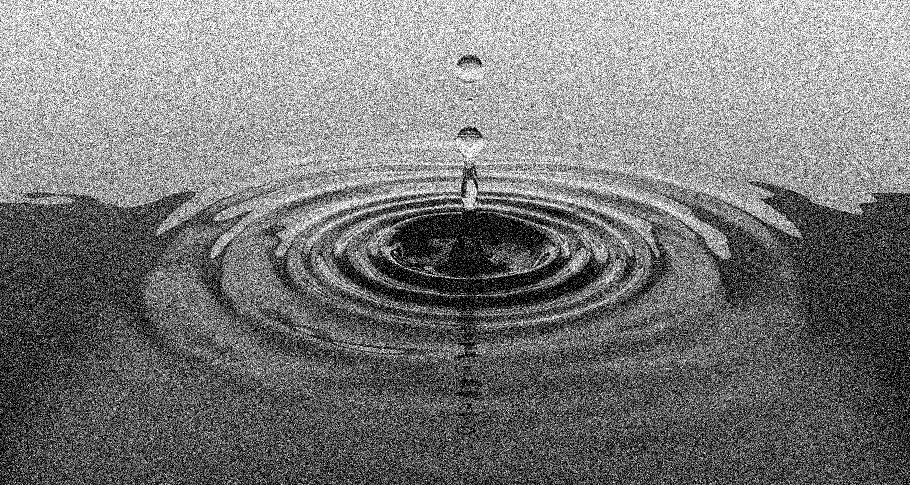}
\captionsetup{justification=centering}
\caption{Noisy $\ubold^\diamond$\\($20\%$ Gaussian)}
\end{subfigure}
\begin{subfigure}[t]{0.235\textwidth}\centering
\includegraphics[width=1\textwidth]{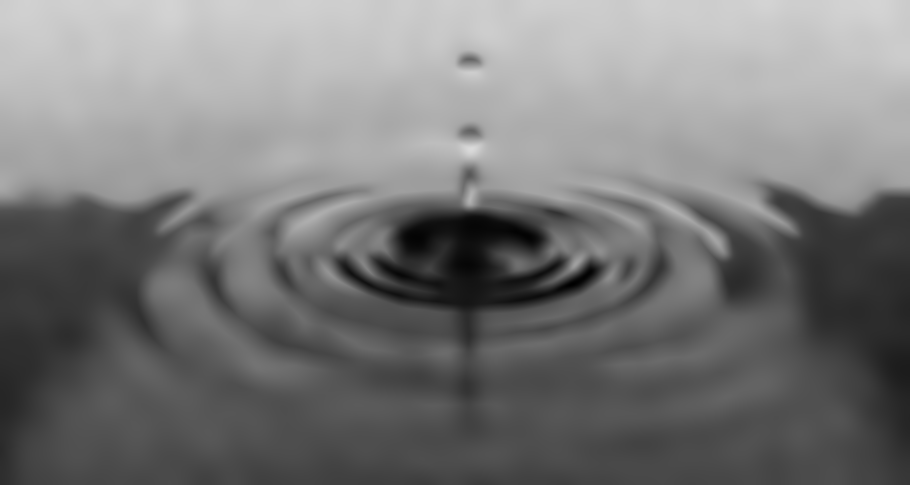}
\captionsetup{justification=centering}
\caption{$\TGV{\alphabold}{2}$\\PSNR: 25.09}
\end{subfigure}
\begin{subfigure}[t]{0.235\textwidth}\centering
\includegraphics[width=1\textwidth]{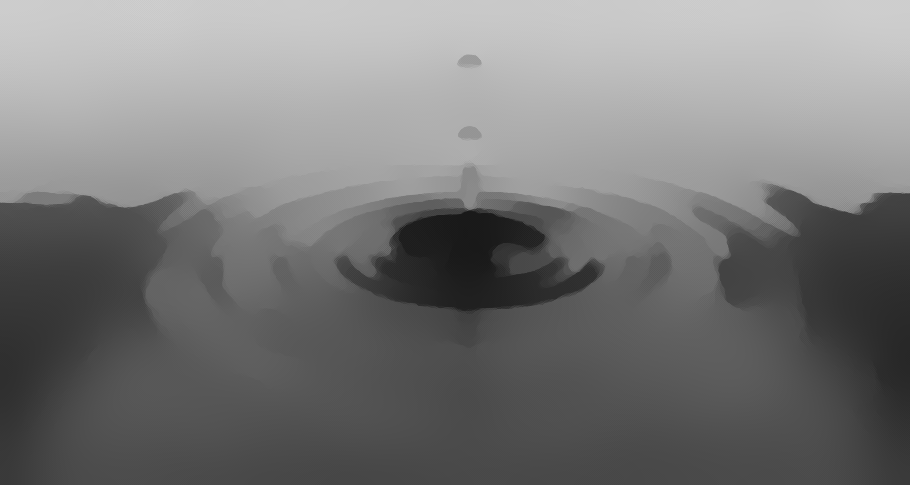}
\captionsetup{justification=centering}
\caption{$\TDVM{1}{\alphabold}[\blank,\Mcal_1]$\\$\Mcal_1=(\Ibold)$, PSNR: 23.78}
\end{subfigure}
\begin{subfigure}[t]{0.235\textwidth}\centering
\includegraphics[width=1\textwidth]{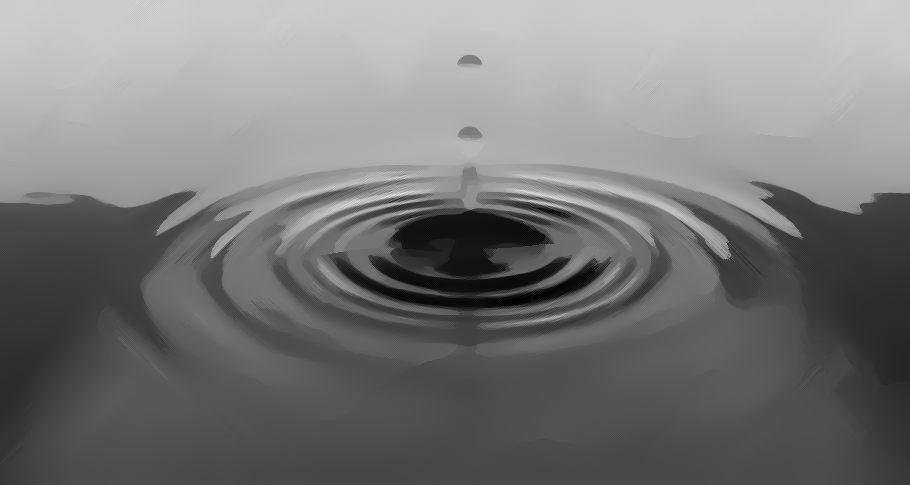}
\captionsetup{justification=centering}
\caption{$\TDVM{1}{\alphabold}[\blank,\Mcal_1]$\\$\Mcal_1=(\Mbold)$, PSNR: 27.19}
\end{subfigure}
\\
\begin{subfigure}[t]{0.235\textwidth}\centering
\captionsetup{justification=centering}
\includegraphics[width=1\textwidth]{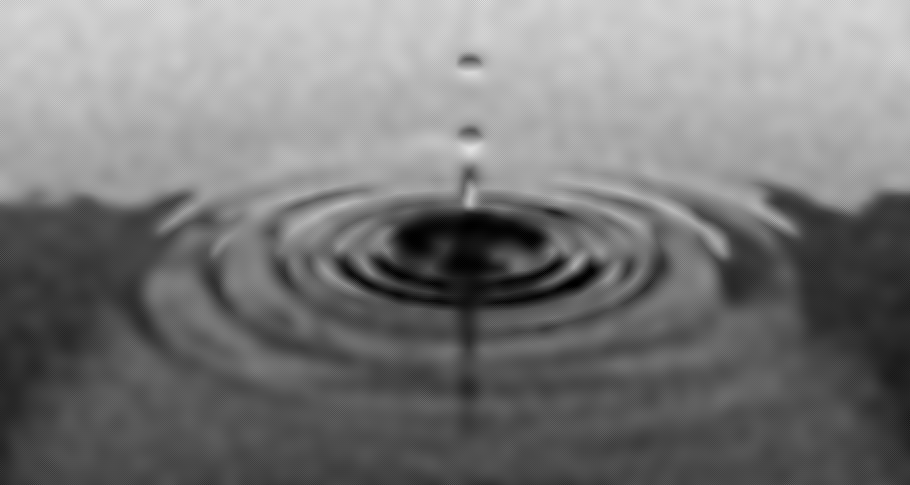}
\caption{$\TDVM{2}{\alphabold}[\blank,\Mcal_2]$\\$\Mcal_2=(\Ibold,\Ibold)$, PSNR: 25.37}
\label{fig: single tdv figure e}
\end{subfigure}
\begin{subfigure}[t]{0.235\textwidth}\centering
\includegraphics[width=1\textwidth]{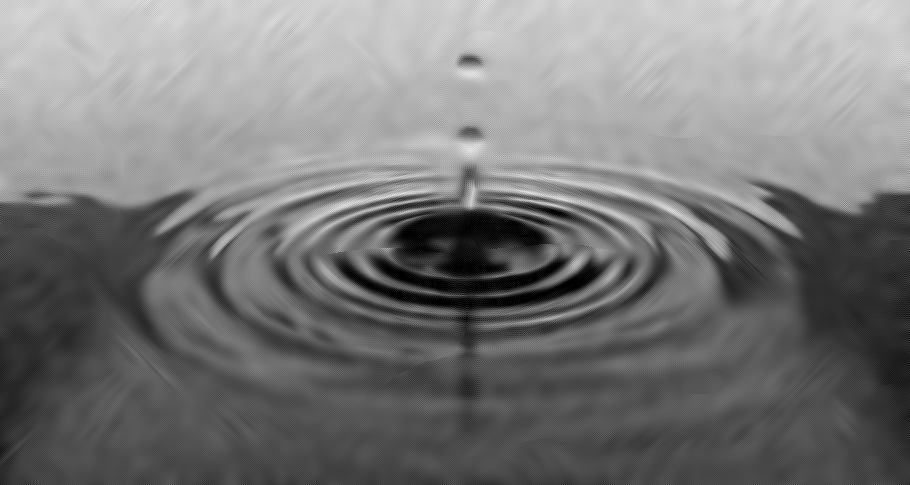}
\captionsetup{justification=centering}
\caption{$\TDVM{2}{\alphabold}[\blank,\Mcal_2]$\\$\Mcal_2=(\Mbold,\Ibold)$, PSNR: 27.37}
\label{fig: single tdv figure f}
\end{subfigure}
\begin{subfigure}[t]{0.235\textwidth}\centering
\captionsetup{justification=centering}
\includegraphics[width=1\textwidth]{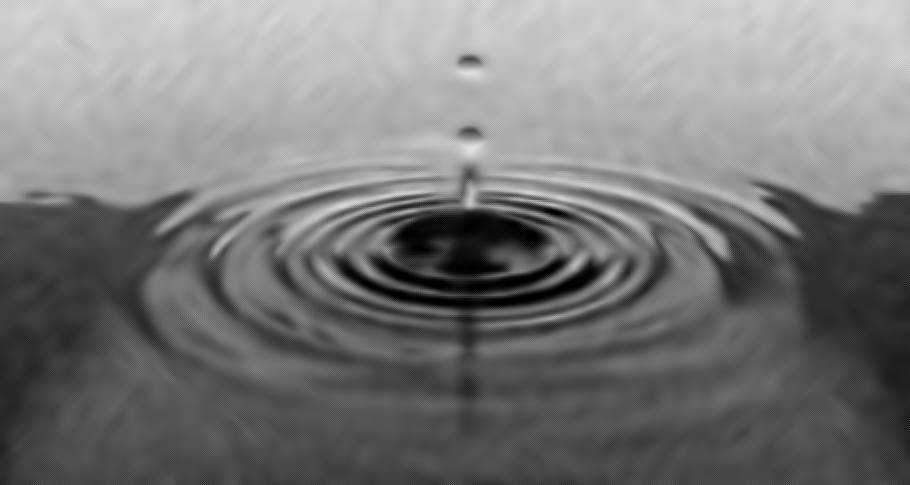}
\caption{$\TDVM{2}{\alphabold}[\blank,\Mcal_2]$\\$\Mcal_2=(\Ibold,\Mbold)$, PSNR: 27.38}
\label{fig: single tdv figure g}
\end{subfigure}
\begin{subfigure}[t]{0.235\textwidth}\centering
\captionsetup{justification=centering}
\includegraphics[width=1\textwidth]{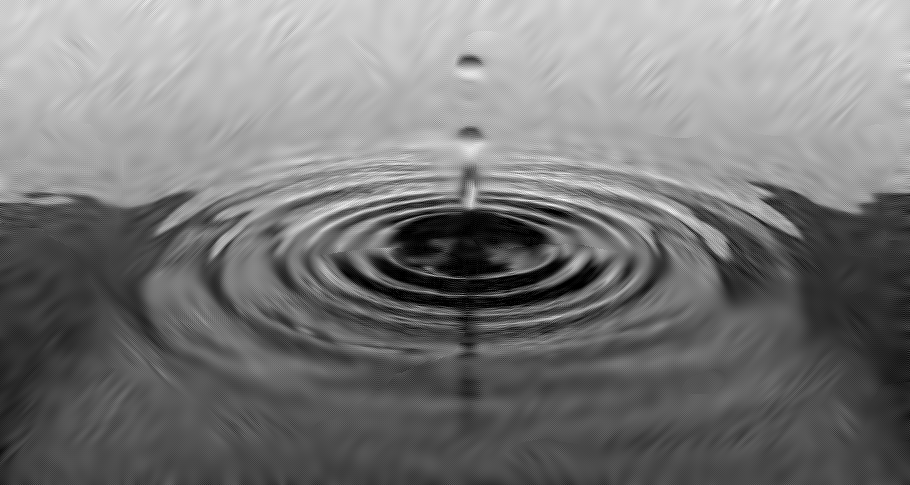}
\caption{$\TDVM{2}{\alphabold}[\blank,\Mcal_2]$\\$\Mcal_2=(\Mbold,\Mbold)$,PSNR: 27.30}
\label{fig: single tdv figure h}
\end{subfigure}
\\
\begin{subfigure}[t]{0.235\textwidth}\centering
\captionsetup{justification=centering}
\includegraphics[width=1\textwidth]{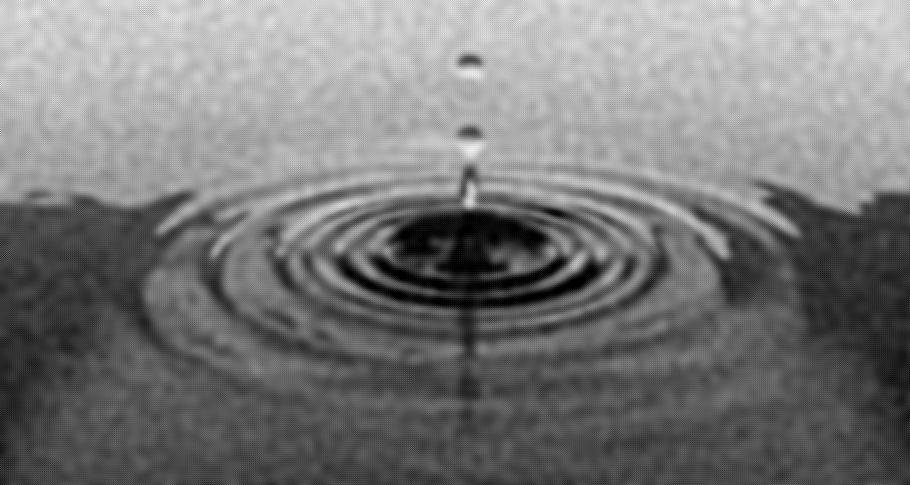}
\caption{$\TDVM{3}{\alphabold}[\blank,\Mcal_3]$\\$\Mcal_3=(\Ibold,\Ibold,\Ibold)$\\PSNR: 25.70}
\end{subfigure}
\begin{subfigure}[t]{0.235\textwidth}\centering
\captionsetup{justification=centering}
\includegraphics[width=1\textwidth]{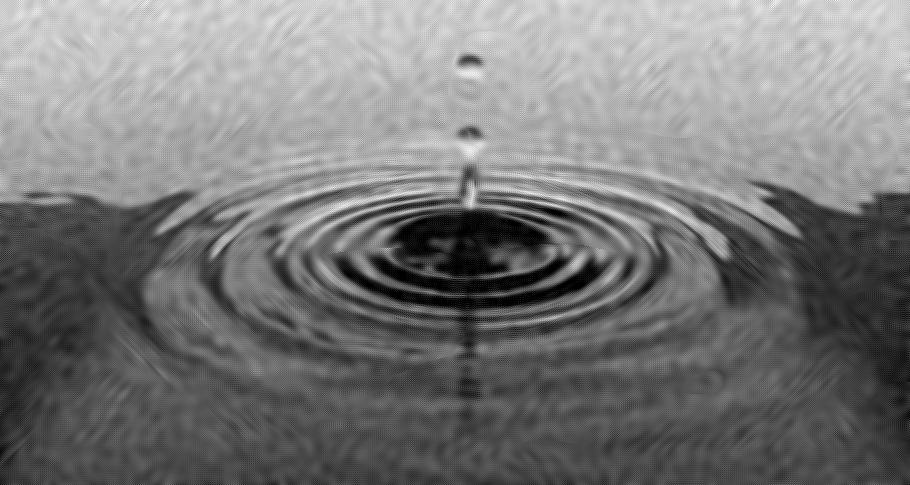}
\caption{$\TDVM{3}{\alphabold}[\blank,\Mcal_3]$\\$\Mcal_3=(\Mbold,\Ibold,\Ibold)$\\PSNR: 25.97}
\end{subfigure}
\begin{subfigure}[t]{0.235\textwidth}\centering
\captionsetup{justification=centering}
\includegraphics[width=1\textwidth]{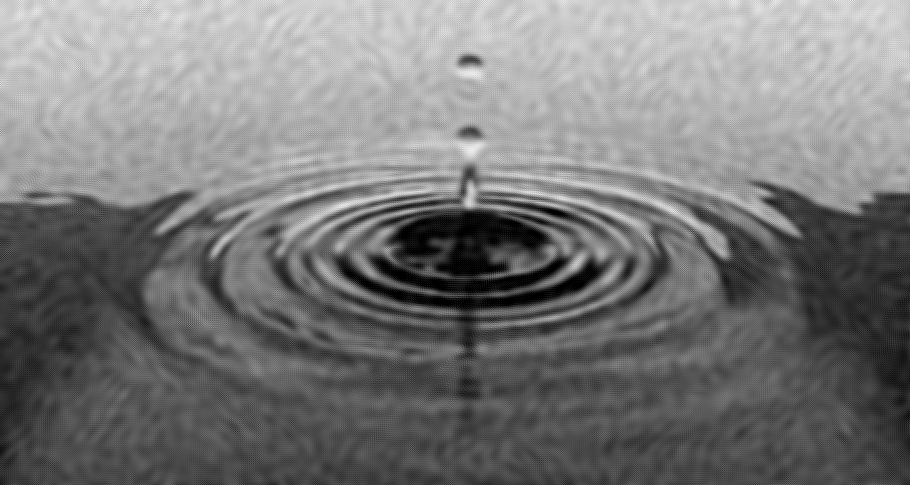}
\caption{$\TDVM{3}{\alphabold}[\blank,\Mcal_3]$\\$\Mcal_3=(\Ibold,\Mbold,\Ibold)$\\PSNR: 26.27}
\end{subfigure}
\begin{subfigure}[t]{0.235\textwidth}\centering
\captionsetup{justification=centering}
\includegraphics[width=1\textwidth]{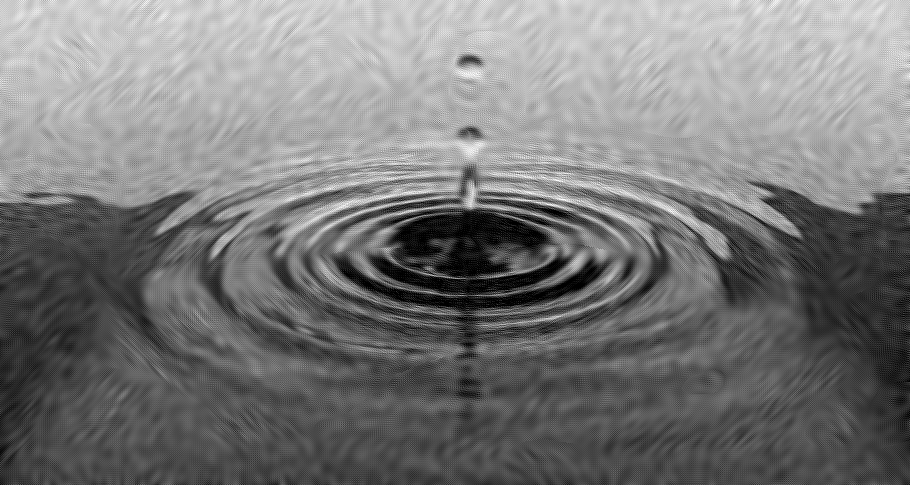}
\caption{$\TDVM{3}{\alphabold}[\blank,\Mcal_3]$\\$\Mcal_3=(\Ibold,\Ibold,\Mbold)$\\PSNR: 25.19}
\end{subfigure}
\\
\begin{subfigure}[t]{0.235\textwidth}\centering
\captionsetup{justification=centering}
\includegraphics[width=1\textwidth]{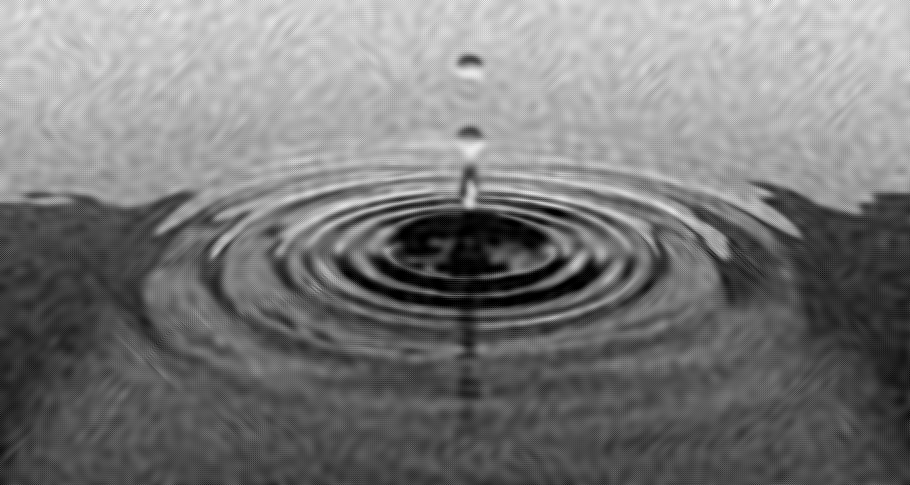}
\caption{$\TDVM{3}{\alphabold}[\blank,\Mcal_3]$\\$\Mcal_3=(\Mbold,\Mbold,\Ibold)$\\PSNR: 26.59}
\end{subfigure}
\begin{subfigure}[t]{0.235\textwidth}\centering
\captionsetup{justification=centering}
\includegraphics[width=1\textwidth]{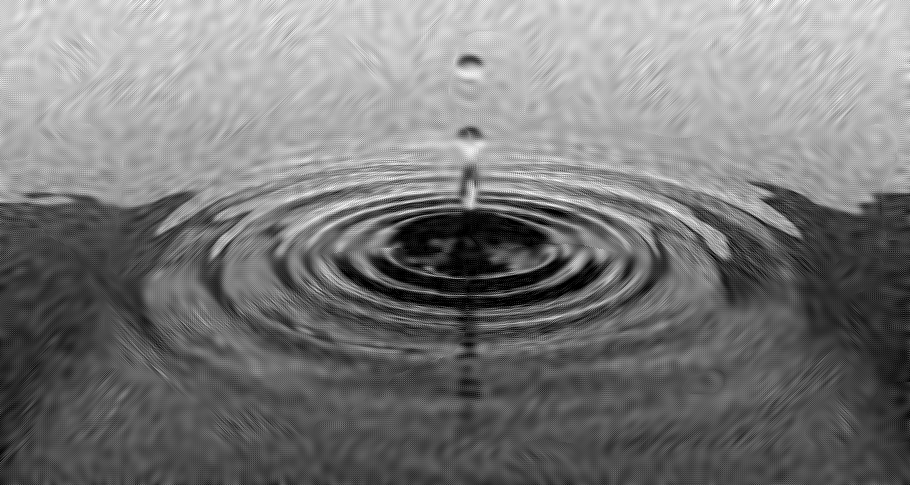}
\caption{$\TDVM{3}{\alphabold}[\blank,\Mcal_3]$\\$\Mcal_3=(\Mbold,\Ibold,\Mbold)$\\PSNR: 25.36}
\end{subfigure}
\begin{subfigure}[t]{0.235\textwidth}\centering
\captionsetup{justification=centering}
\includegraphics[width=1\textwidth]{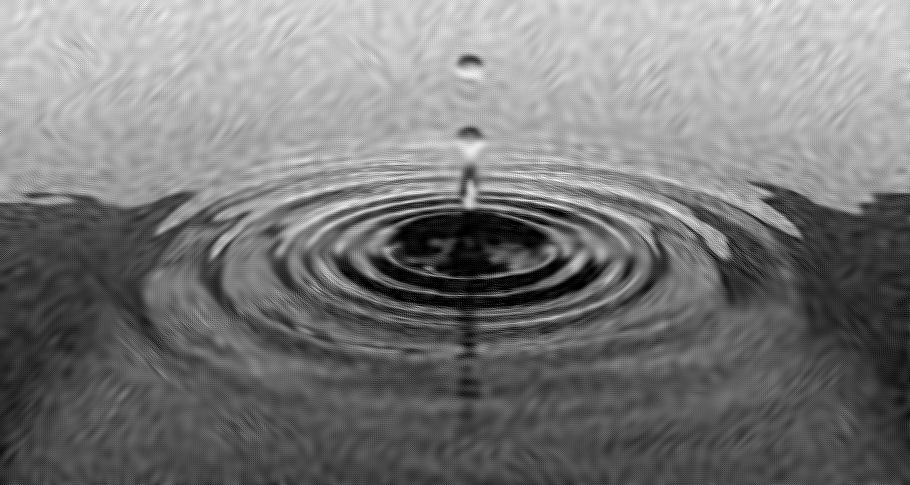}
\caption{$\TDVM{3}{\alphabold}[\blank,\Mcal_3]$\\$\Mcal_3=(\Ibold,\Mbold,\Mbold)$\\PSNR: 26.04}
\end{subfigure}
\begin{subfigure}[t]{0.235\textwidth}\centering
\captionsetup{justification=centering}
\includegraphics[width=1\textwidth]{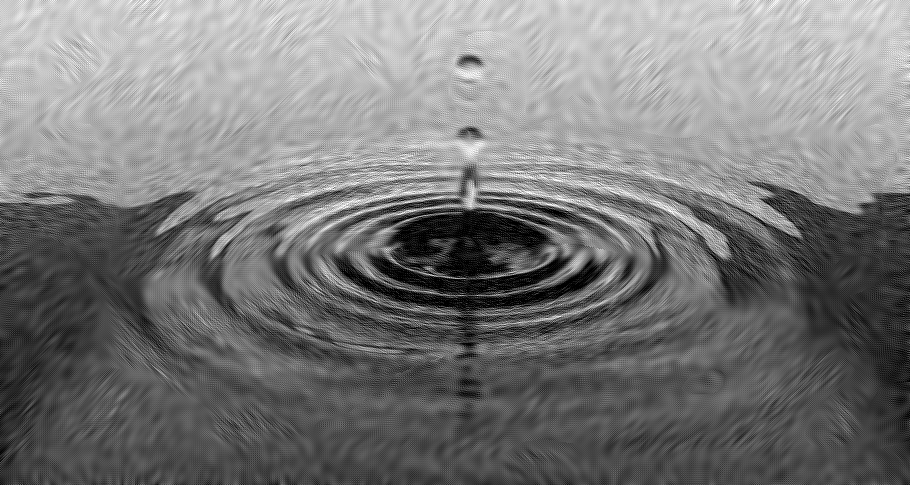}
\caption{$\TDVM{3}{\alphabold}[\blank,\Mcal_3]$\\$\Mcal_3=(\Mbold,\Mbold,\Mbold)$\\PSNR: 23.81}
\end{subfigure}
\caption{Here $\alphabold=(\alpha_j)_{j=0}^{q-1}$ with $\alpha_0 = 1$ and $\alpha_j=1.25$ for each $j>0$, $(\sigma,\rho)=(2,25)$ in the structure tensor computation, $\bbold=(1,\beta(\cdot))$ and fidelity $\eta = 1$.}
\label{fig: single tdv figure}
\end{figure}

In \cref{fig: single tdv figure} we report the results from the single $\TDVM{q}{\alphabold}-\LL{2}$ problem in \cref{eq: discretised problem fixed q}, with fixed $q=1,2,3$ and different choices for the anisotropy $\Mcal_q$. To compute $\Mbold$, we employ the spatially varying strategy based on the structure tensor eigen-decomposition described in \cref{sec: estimate v}, with $(\sigma,\rho)=(2,25)$ and $\bbold=(1,\beta(\xbold))$, where $\beta(\xbold)$ follows from \cref{eq: denoising weight 2,eq: compute b vary step 1,eq: compute b vary step 2}.
In these experiments, we keep fixed $\eta=1$, $\alpha_0=1$ and $\alpha_j=1.25$ for $j=1,\dots,q$ and we run \cref{alg: main primal-dual} for $\texttt{maxiter}=1000$.  We compare our results with $\TGV{\alphabold}{2}$ (with the same weights) from an online source\footnote{Code available at \url{www.gipsa-lab.grenoble-inp.fr/~laurent.condat/download/TGVdenoise.m}}: our strategy is able to encode spatially varying directional information at all the derivative orders.

\begin{figure}[htb]
\centering
\begin{subfigure}[t]{0.235\textwidth}\centering
\captionsetup{justification=centering}
\includegraphics[width=1\textwidth]{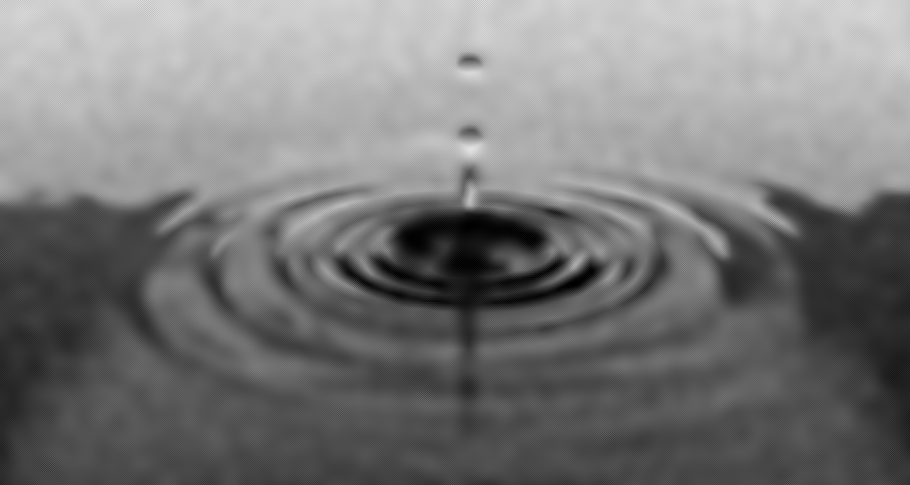}
\caption{$\TDVM{2}{\alphabold}[\blank,\Mcal_2]$\\$\Mcal_2=(\Ibold,\Ibold)$\\PSNR: 25.37}
\label{fig: 0307 a}
\end{subfigure}
\begin{subfigure}[t]{0.235\textwidth}\centering
\includegraphics[width=1\textwidth]{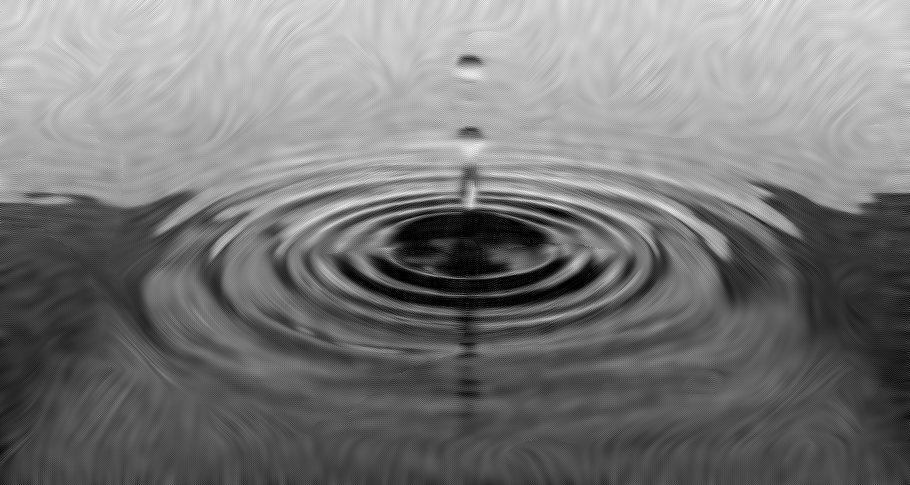}
\captionsetup{justification=centering}
\caption{$\TDVM{2}{\alphabold}[\blank,\Mcal_2]$\\$\Mcal_2=(\Mbold^{0.3},\Ibold)$\\PSNR: 27.16}
\label{fig: 0307 b}
\end{subfigure}
\begin{subfigure}[t]{0.235\textwidth}\centering
\captionsetup{justification=centering}
\includegraphics[width=1\textwidth]{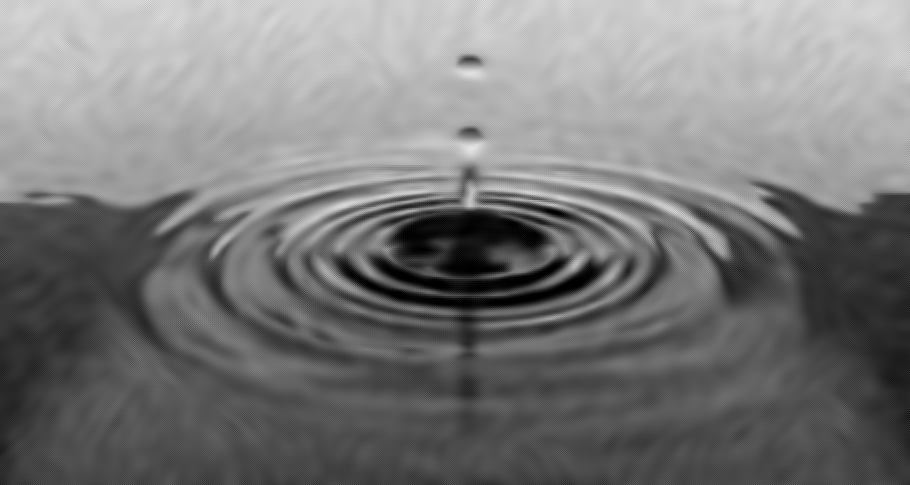}
\caption{$\TDVM{2}{\alphabold}[\blank,\Mcal_2]$\\$\Mcal_2=(\Ibold,\Mbold^{0.7})$\\PSNR: 27.57}
\label{fig: 0307 c}
\end{subfigure}
\begin{subfigure}[t]{0.235\textwidth}\centering
\captionsetup{justification=centering}
\includegraphics[width=1\textwidth]{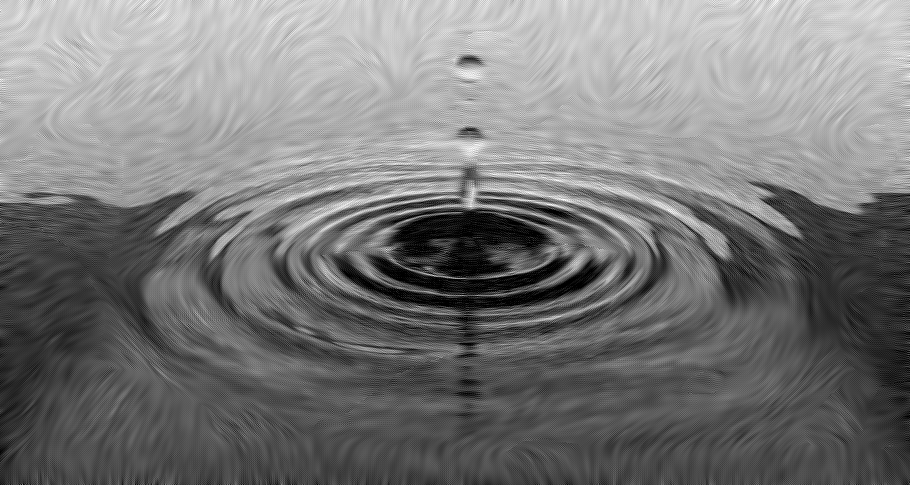}
\caption{$\TDVM{2}{\alphabold}[\blank,\Mcal_2]$\\$\Mcal_2=(\Mbold^{0.3},\Mbold^{0.7})$\\PSNR: 25.71}
\label{fig: 0307 d}
\end{subfigure}
\\
\begin{subfigure}[t]{0.235\textwidth}\centering
\captionsetup{justification=centering}
\includegraphics[width=1\textwidth]{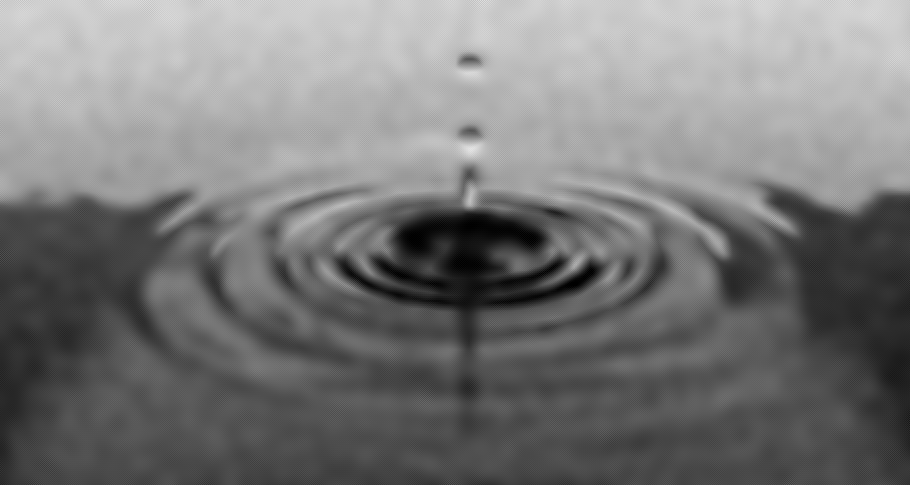}
\caption{$\TDVM{2}{\alphabold}[\blank,\Mcal_2]$\\$\Mcal_2=(\Ibold,\Ibold)$\\PSNR: 25.37}
\label{fig: 0703 a}
\end{subfigure}
\begin{subfigure}[t]{0.235\textwidth}\centering
\includegraphics[width=1\textwidth]{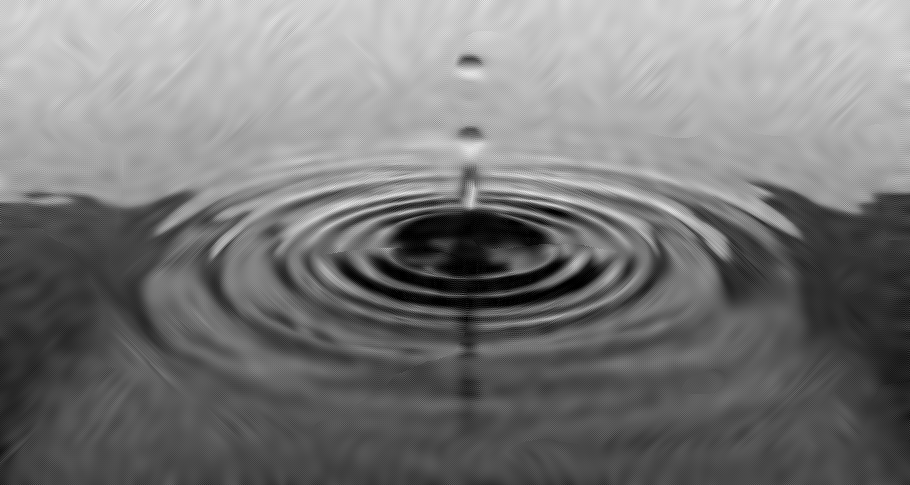}
\captionsetup{justification=centering}
\caption{$\TDVM{2}{\alphabold}[\blank,\Mcal_2]$\\$\Mcal_2=(\Mbold^{0.7},\Ibold)$\\PSNR: 27.45}
\label{fig: 0703 b}
\end{subfigure}
\begin{subfigure}[t]{0.235\textwidth}\centering
\captionsetup{justification=centering}
\includegraphics[width=1\textwidth]{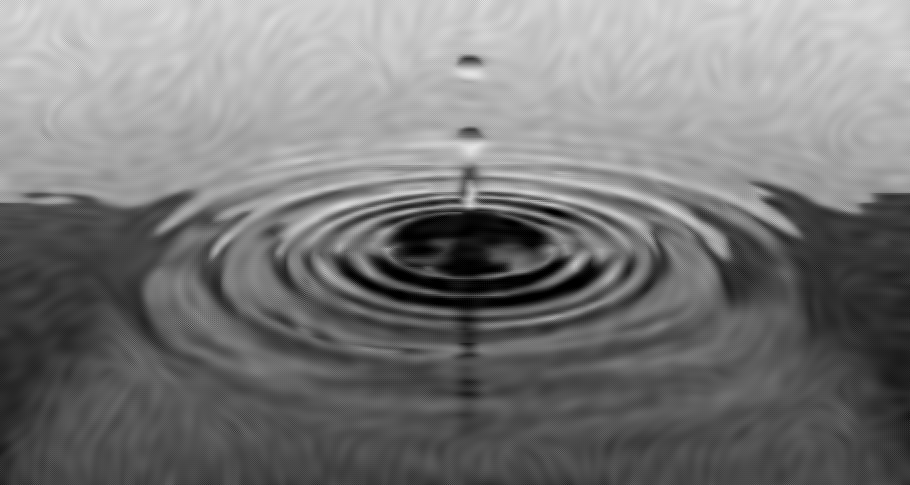}
\caption{$\TDVM{2}{\alphabold}[\blank,\Mcal_2]$\\$\Mcal_2=(\Ibold,\Mbold^{0.3})$\\PSNR: 27.77}
\label{fig: 0703 c}
\end{subfigure}
\begin{subfigure}[t]{0.235\textwidth}\centering
\captionsetup{justification=centering}
\includegraphics[width=1\textwidth]{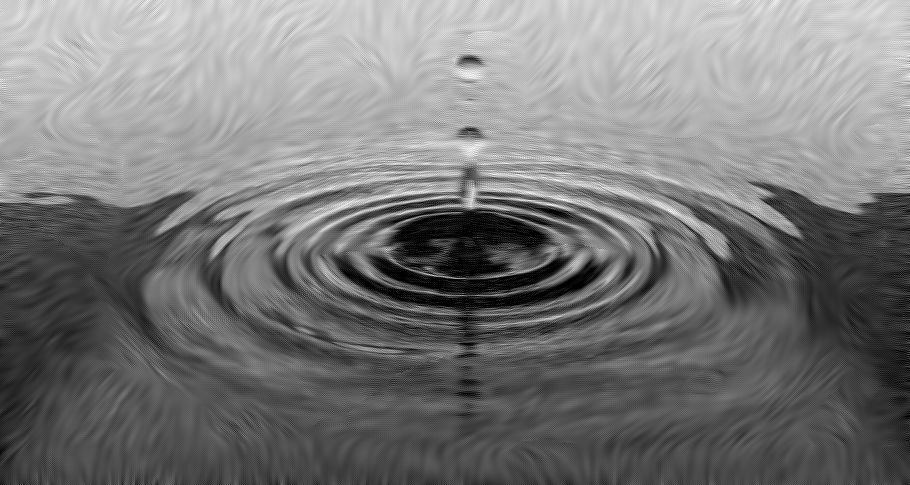}
\caption{$\TDVM{2}{\alphabold}[\blank,\Mcal_2]$\\$\Mcal_2=(\Mbold^{0.7},\Mbold^{0.3})$\\PSNR: 26.18}
\label{fig: 0703 d}
\end{subfigure}
\vspace{-1em}
\caption{Here $\alphabold=(\alpha_j)_{j=0}^{q-1}$ with $\alpha_0 = 1$ and $\alpha_j=1.25$ for each $j>0$, $(\sigma,\rho)=(2,25)$, $\eta = 1$. Starting from the anisotropy matrix $\Mbold$ of  \cref{fig: single tdv figure}, we define two variants: $\Mbold^{0.3}$ is obtained from $\Mbold$ by replacing $\bbold=(1,\beta(\xbold))$ with $\bbold_{0.3}=(1,0.3\beta(\xbold))$, and $\Mbold^{0.7}$ is obtained from $\Mbold$ by replacing $\bbold$ with $\bbold_{0.7}=(1,0.7\beta(\xbold))$. All other parameters are the same as in  \cref{fig: single tdv figure}. We report in \cref{fig: 0307 a,fig: 0307 b,fig: 0307 c,fig: 0307 d} the results obtained with four possible choices of $\Mcal_2$:  $(\Ibold,\Ibold)$, $(\Mbold^{0.3},\Ibold)$, $(\Ibold,\Mbold^{0.7})$ and $(\Mbold^{0.3},\Mbold^{0.7})$, respectively. The results in \cref{fig: 0703 a,fig: 0703 b,fig: 0703 c,fig: 0703 d} correspond to the choices $(\Ibold,\Ibold)$, $(\Mbold^{0.7},\Ibold)$, $(\Ibold,\Mbold^{0.3})$ and $(\Mbold^{0.7},\Mbold^{0.3})$, respectively.}
\label{fig: single tdv figure ext1}
\end{figure}

For the experiments of  \cref{fig: single tdv figure ext1}, we use the same anisotropy directions as in the experiments of  \cref{fig: single tdv figure}, but we modulate the anisotropy weight $\beta$ to test its impact on the result. More precisely, starting from the same anisotropy matrix $\Mbold$ of  \cref{fig: single tdv figure}, we define two variants: $\Mbold^{0.3}$ is obtained from $\Mbold$ by replacing $\bbold=(1,\beta(\xbold))$ with $\bbold_{0.3}=(1,0.3\beta(\xbold))$, and $\Mbold^{0.7}$ is obtained from $\Mbold$ by replacing $\bbold$ with $\bbold_{0.7}=(1,0.7\beta(\xbold))$. The results of \cref{fig: single tdv figure ext1} correspond to different anisotropic weighting of the first and second-order derivatives in the model $\TDVM{2}{\alphabold}-\LL{2}$, see the figure caption for details. 
Interestingly, a better PSNR is obtained with two orders of derivation and the following strategies: if the first derivative is anisotropically weighted, it is better to opt for an isotropic weight of the second-order derivative. And if the first derivative is isotropically weighted, it is better to opt for a strong anisotropic weight of the second-order derivative. This latter case seems the best in terms of PSNR but, of course, sharper results are obtained when anisotropic weights are chosen. Obviously, the PSNR metric does not capture all aspects of an image and a visual evaluation by the end user remains necessary for the choice of the most suitable combination of derivative orders and anisotropic weights.

\begin{remark}\label{rem: artefacts}
In our experiments we observed some checkerboard artifacts when employing strong anisotropies. This phenomenon has been attributed to spectral properties of finite differences \cite{WEICKERT2002103,Fehrenbach2014} and a non-negative stencil avoiding these issues has been introduced in \cite{Fehrenbach2014} for the case of directional Hessian. 
We leave a generalization of \cite{Fehrenbach2014} to our higher-order case for future work. The artifacts are not evident in the results of the next \cref{sec: jointmod}, where the joint model and a specific choice of the weights help in producing better quality results.
\end{remark}

Our approach has a reasonably large number of parameters depending on the noise level and the structural imaging information: the order of the derivatives, the spatially varying directional information from the structure tensor and the weights for the sparsity of the derivatives. This makes the model highly customisable, but also quite sensitive to parameters. 
In order to simplify the parameter choice, we detail in \cref{sec: jointmod} a simplified joint model with a suggested choice of parameters, particularly suited for images with large dominant structures, as well as for the reconstruction of images from scattered data, see \cref{sec: model-description-surface}.

\subsection{Joint minimisation model}\label{sec: jointmod}
Following the comments in \cref{rem: artefacts}, we now consider the joint model in
\cref{eq: minimization intro} and described in  \cref{sec: joint description}
for $\Qrm=3$, with $\Mcal_q=(\Ibold,\dots,\Ibold,\Mbold)$ and $\alphabold_q = (\alpha_{q,0},+\infty,\dots,+\infty)$ for each $q=1,2,3$. Since we consider the denoising problem, we will make use of the identity operator $\Scal=\Ibold$ in the fidelity term, aiming to solve:
\begin{equation}
u^\star
=
\argmin_{u\in\RR^{|\Omega^h|}} \sum_{q=1}^\Qrm \TDVM{q}{\alphabold_q}[u,\Mcal_q]+ \frac{\eta}{2} \norm{u -u^\diamond}_2^2.
\label{eq: energy denoising}
\end{equation}

Thus, the denoising problem in equation \cref{eq: energy denoising} can be simplified as
\begin{equation}
u^\star
=
\argmin_{u\in\RR^{|\Omega^h|}}  \left( \sum_{q=1}^3
\alpha_{q,0} \norm{\Mbold\Wcal^q\grad^q u}_{2,1} + \frac{\eta}{2} \norm{u -u^\diamond}_2^2
\right),
\label{eq: energy denoising simplified}
\end{equation}
and solved via the primal-dual \cref{alg: reduced main primal-dual} in which the computations for $u$ and $\vbold$ are performed alternatingly as described in \cref{alg: main denoising} so as to update the estimation of the directions $\vbold$. This updating strategy will allow a refinement for the computation of the directionality in images. 

\begin{algorithm}[tbhp]
\small
\caption{Discrete joint denoising problem \cref{eq: energy denoising simplified}}
\SetAlgoLined
\label{alg: main denoising}
\SetKwProg{Fn}{Function}{:}{}
\SetKwFunction{TDVdenoising}{TDV\_denoising}
\SetKwFunction{computevb}{compute\_v\_and\_b}
\SetKwFunction{buildKcal}{build\_K}
\SetKwFunction{primaldual}{PrimalDual\_reduced}
\SetKwFunction{computeopnorm}{compute\_operator\_norm}\SetKwFunction{updatesigmarho}{update\_sigmarho}
\SetKwData{maxitert}{$T$}
\SetKwData{maxiterk}{\texttt{maxiter}}
\SetKwData{F}{F}
\SetKwData{G}{G}
\SetKwData{K}{K}
\SetKwData{KS}{KS}
\SetKwData{ProxF}{ProxF}
\SetKwData{ProxFS}{ProxFS}
\SetKwData{ProxG}{ProxG}
\SetKwInOut{Input}{Input}
\SetKwInOut{ParametersModel}{Parameters for the model}
\SetKwInOut{ParametersV}{Parameters for $\vbold$}
\SetKwInOut{ParametersPD}{Parameters for Primal Dual}
\SetKwInOut{Operators}{Operators needed}
\SetKwInOut{Initialization}{Initialization}
\SetKwInOut{Update}{Update}
\BlankLine
\ParametersModel{$\abold=(\alpha_{1,0},\alpha_{2,0},\alpha_{3,0}),\,\bbold^h=(\Ibold,\bbold_2^h),\,\eta >0$, $\maxiterk$.}
\ParametersV{$\sigma_1,\,\rho_1 >0$.}
\BlankLine
\Fn{\TDVdenoising{$\ubold^{\diamond,h}$}}{
\BlankLine
\For{$k=1,\dots,\maxiterk$}{
\BlankLine
$[\vbold^{k+1,h},\bbold^{k+1,h}]$ = \computevb{$\ubold^{k},\sigma_k,\rho_k$}\tcp*[r]{from \cref{eq: structure tensor}-\cref{eq: compute b vary step 2}}
$\Kcal^{k+1,h}$ = \buildKcal($\abold$, $\vbold^{k+1,h}$, $\bbold^{k+1,h}$) \tcp*[r]{update $\Kcal=(\Kcal_q^h)_{q=1}^\Qrm$ with $\Kcal_q^h=\Mbold\Wcal^q\grad^{q,h}$}
$\ubold^{k+1,h}$ = \primaldual{$\ubold^{\diamond,h},\Kcal^{k+1,h},\abold,\eta$} \tcp*[r]{restart denoising with new $\Kcal^{k+1}$}
$[\sigma_{k+1},\rho_{k+1}]$ = \updatesigmarho($\sigma_{k},\rho_{k}$) \tcp*[r]{new structure tensor parameters}
\BlankLine
}
}
\Return{$\ubold^{k+1,h}$}
\BlankLine
\end{algorithm}

\subsubsection*{Numerical Results}
We discuss denoising results obtained with \cref{alg: main denoising} and different denoising approaches (the non-local method BM3D with normal-complexity profile and prior knowledge of the standard deviation of the Gaussian noise \cite{BM3D07} and the regularisers $\TV{}{}$ \cite{ROF}, $\TGV{}{2}$ \cite{BreKunPoc2010}, $\DTV{}{}$ and  $\text{DTGV}^{2}$ \cite{directionaltv}) for images with strong directional features.
We run the primal-dual \cref{alg: reduced main primal-dual} for $500$ iterations and we restarted the \cref{alg: main denoising} once the first denoised image is computed which serves as an oracle for a better estimation of $\vbold$. 
We show results for grey-scale images in \cref{fig: bamboo all orders} and for colour images \cref{fig: denoising results color,fig: denoising desert} and we discuss their PSNR. 

\begin{figure}[tbph]
\centering
\begin{subfigure}[t]{0.235\textwidth}\centering\captionsetup{justification=centering}
\includegraphics[width=0.85\textwidth]{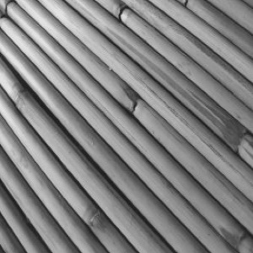}
\caption{Original $\ubold$\\ $253\times 253$ pixels}
\label{fig: bamboo original}
\end{subfigure}
\,
\begin{subfigure}[t]{0.235\textwidth}\centering\captionsetup{justification=centering}
\includegraphics[width=0.85\textwidth]{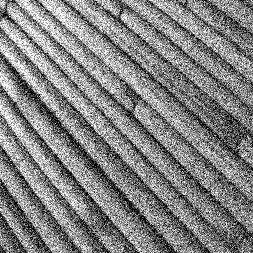}
\caption{Noisy $\ubold^\diamond$\\($20\%$ Gaussian)}
\label{fig: bamboo 20noise}
\end{subfigure}
\,
\begin{subfigure}[t]{0.235\textwidth}\centering\captionsetup{justification=centering}
\includegraphics[width=0.85\textwidth]{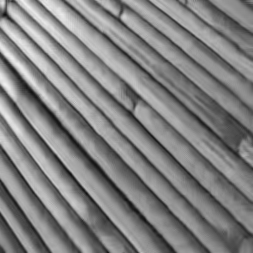}
\caption{BM3D\\
PSNR = 29.1 - (12 s.)}
\end{subfigure}
\\
\begin{subfigure}[t]{0.235\textwidth}\centering\captionsetup{justification=centering}
\includegraphics[width=0.85\textwidth]{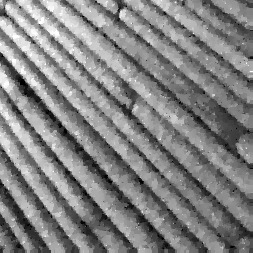}
\caption{$\TV{}$ (from \cite{directionaltv})\\ PSNR = 23.8 - (9 s.)}
\label{fig: bamboo TV}
\end{subfigure}
\,
\begin{subfigure}[t]{0.235\textwidth}\centering\captionsetup{justification=centering}
\includegraphics[width=0.85\textwidth]{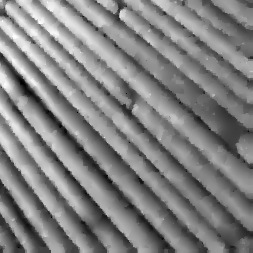}
\caption{$\TGV{}{2}$ (from \cite{directionaltv})\\ PSNR = 24.7 - (29 s.)}
\label{fig: bamboo TGV}
\end{subfigure}
\,
\begin{subfigure}[t]{0.235\textwidth}\centering\captionsetup{justification=centering}
\includegraphics[width=0.85\textwidth]{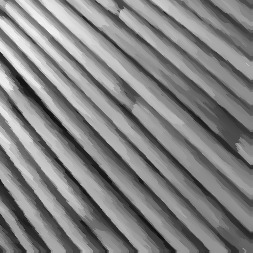}
\caption{$\DTV{}$ (from \cite{directionaltv})\\ PSNR = 26.8 - (12 s.)}
\label{fig: bamboo DTV}
\end{subfigure}
\,
\begin{subfigure}[t]{0.235\textwidth}\centering\captionsetup{justification=centering}
\includegraphics[width=0.85\textwidth]{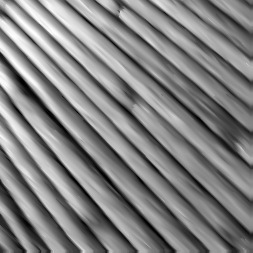}
\caption{$\DTGV{2}{}$ (from \cite{directionaltv})\\ PSNR = 28.2 - (37 s.)}
\label{fig: bamboo DTGV}
\end{subfigure}
\vspace{-1em}
\caption{Denoising of bamboo grey-scale image: current approaches.}
\label{fig: bamboo denoising current results}
\vspace{-1.5em}
\end{figure}

\paragraph{Bamboo image}
The grey-scale image in \cref{fig: bamboo original} shows a strong directional direction. 
In \cref{fig: bamboo 20noise} it has been corrupted by 20\% of Gaussian noise using the same random seed as in \cite{directionaltv}, see \cref{fig: bamboo 20noise}. 
In \cref{fig: bamboo denoising current results} we report the results from state of the art approaches, as reported in \cite{directionaltv}, where
$\DTV{}{}$ and $\DTGV{}{}$ are considered with a single fixed choice of anisotropy direction and minor semi-axis (our parameter $\beta$).
In \cref{fig: bamboo TV,fig: bamboo TGV} the staircasing effect is visible as expected while \cref{fig: bamboo DTV} and \cref{fig: bamboo DTGV} seem more promising, even if obtained with a single fixed direction only. 

\begin{figure}[tbp]
\centering
\begin{subfigure}[t]{0.235\textwidth}\centering\captionsetup{justification=centering}
\includegraphics[width=0.85\textwidth]{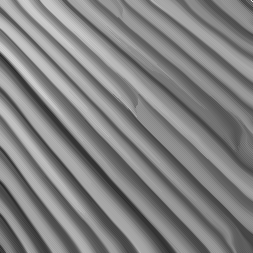}
\caption{$\TDVM{1}{}$, $\eta=0.25$\\PSNR: 24.15\\(13 s.\ per $k$ in Alg.\ \ref{alg: main denoising})}
\label{fig: bamboo eta0.25}
\end{subfigure}
\,
\begin{subfigure}[t]{0.235\textwidth}\centering\captionsetup{justification=centering}
\includegraphics[width=0.85\textwidth]{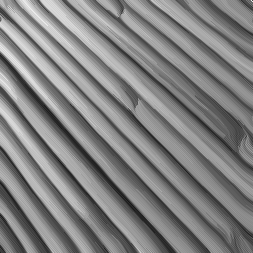}
\caption{$\TDVM{1}{}$, $\eta=0.75$\\PSNR: 27.23\\(14 s.\ per $k$ in Alg.\ \ref{alg: main denoising})}
\label{fig: bamboo eta0.75}
\end{subfigure}
\,
\begin{subfigure}[t]{0.235\textwidth}\centering\captionsetup{justification=centering}
\includegraphics[width=0.85\textwidth]{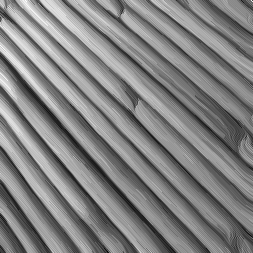}
\caption{$\TDVM{1}{}$, $\eta=1.25$\\PSNR: 26.73\\(11 s.\ per $k$ in Alg.\ \ref{alg: main denoising})}
\label{fig: bamboo eta1.25}
\end{subfigure}
\,
\begin{subfigure}[t]{0.235\textwidth}\centering\captionsetup{justification=centering}
\includegraphics[width=0.85\textwidth]{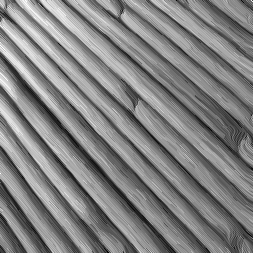}
\caption{$\TDVM{1}{}$, $\eta=1.75$\\PSNR: 25.34\\(11 s.\ per $k$ in Alg.\ \ref{alg: main denoising})}
\label{fig: bamboo eta1.75}
\end{subfigure}
\vspace{-1em}
\caption{Denoising of \cref{fig: bamboo 20noise} for different fidelity weights $\eta$ and first order $\TDVM{1}{\alphabold}$ regulariser with $\alphabold=\alpha=1$. Other parameters: \cref{alg: main denoising} with $\bbold=(1,0.02)$ for the construction of $\Mbold$, $(\sigma_1,\rho_1)=(1.8,2.8)$ and $(\sigma_2,\rho_2)=(1,1)$.
}
\label{fig: bamboo first sensitive eta}
\vspace{-1.5em}
\end{figure}

\begin{figure}[tbph]
\centering
\begin{subfigure}[t]{0.235\textwidth}\centering\captionsetup{justification=centering}
\includegraphics[width=0.85\textwidth]{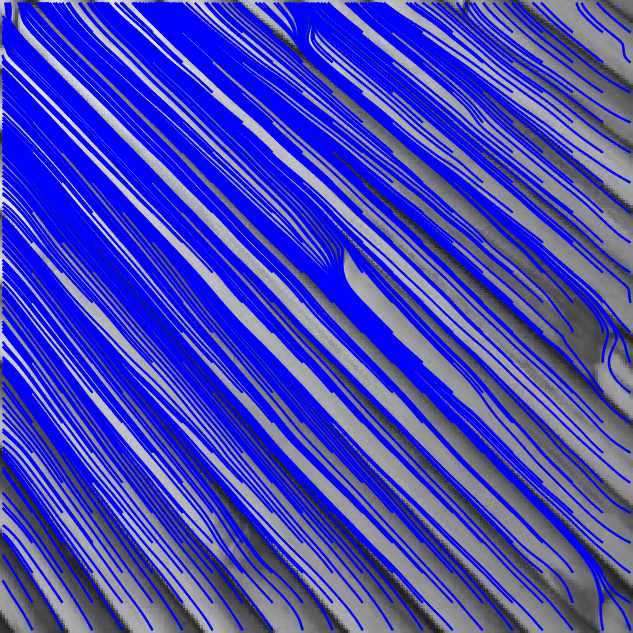}
\caption{streamlines of $\vbold$ \cref{eq: find v denoising} 
\\  
$(\sigma_1,\rho_1)=(1.8,2.8)$\\
$(\sigma_2,\rho_2)=(1,1)$\\
$\texttt{maxiter}=2$ in $k$
}
\label{fig: v bamboo}
\end{subfigure}
\,
\begin{subfigure}[t]{0.235\textwidth}\centering\captionsetup{justification=centering}
\includegraphics[width=0.85\textwidth]{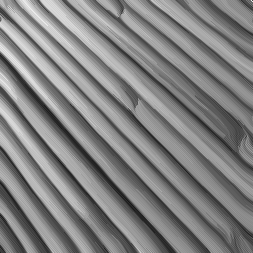}
\caption{$\TDVM{1}{\alphabold}$ $\eta = 0.75$\\PSNR = 27.23\\
(14 s.\ per $k$ in Alg.\ \ref{alg: main denoising})}
\end{subfigure}
\,
\begin{subfigure}[t]{0.235\textwidth}\centering\captionsetup{justification=centering}
\includegraphics[width=0.85\textwidth]{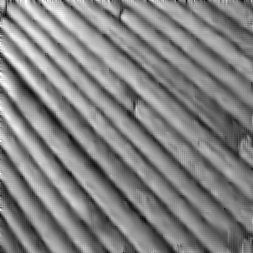}
\caption{$\TDVM{2}{\alphabold}$, $\eta = 3$\\PSNR = 25.43\\
(21 s.\ per $k$ in Alg.\ \ref{alg: main denoising})}
\end{subfigure}
\,
\begin{subfigure}[t]{0.235\textwidth}\centering\captionsetup{justification=centering}
\includegraphics[width=0.85\textwidth]{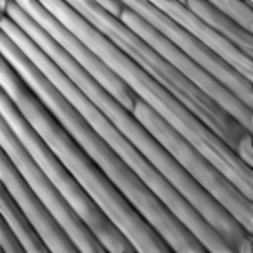}
\caption{$\TDVM{3}{\alphabold}$, $\eta = 1.75$\\PSNR = 27.99\\
(60 s.\ per $k$ in Alg.\ \ref{alg: main denoising})}
\end{subfigure}
\\
\begin{subfigure}[t]{0.235\textwidth}\centering\captionsetup{justification=centering}
\includegraphics[width=0.85\textwidth]{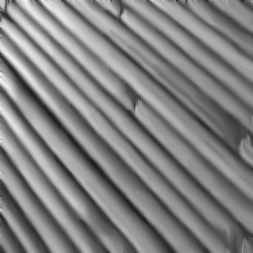}
\caption{$\TDVM{1}{}+\TDVM{2}{}$, $\eta = 4$\\PSNR = 26.47\\
(38 s.\ per $k$ in Alg.\ \ref{alg: main denoising})}
\end{subfigure}
\,
\begin{subfigure}[t]{0.235\textwidth}\centering\captionsetup{justification=centering}
\includegraphics[width=0.85\textwidth]{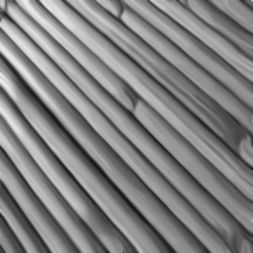}
\caption{$\TDVM{2}{}+\TDVM{3}{},\,\eta=3.75$\\PSNR = 25.95\\
(71 s.\ per $k$ in Alg.\ \ref{alg: main denoising})}
\end{subfigure}
\,
\begin{subfigure}[t]{0.235\textwidth}\centering\captionsetup{justification=centering}
\includegraphics[width=0.85\textwidth]{{./images/denoising/bamboo/bestorder/20noisy_eta3.5_b0.02_a1-0-1_sigma1.8_rho2.8_iter2_PSNR29.1984_time143.18}.png}
\caption{$\TDVM{1}{}+\TDVM{3}{}$, $\eta = 3.5$\\PSNR = 29.20\\
(75 s.\ per $k$ in Alg.\ \ref{alg: main denoising})}
\label{fig: bamboo 101}
\end{subfigure}
\,
\begin{subfigure}[t]{0.235\textwidth}\centering\captionsetup{justification=centering}
\includegraphics[width=0.85\textwidth]{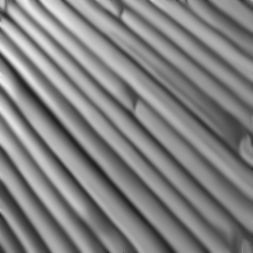}
\caption{$\sum_{q=1}^3\TDVM{q}{}$, $\eta = 4$\\PSNR = 26.19\\
(87 s.\ per $k$ in Alg.\ \ref{alg: main denoising})}
\end{subfigure}
\vspace{-1em}
\caption{
Denoising of \cref{fig: bamboo 20noise} with combinations of $\TDVM{q}{\alphabold_q}$ with $\alphabold_1=1$, $\alphabold_2=(1,\infty)$, and $\alphabold_3=(1,\infty,\infty)$ and with homogeneously fixed $\bbold=(1,0.02)$.}
\label{fig: bamboo all orders}
\vspace{-1.5em}
\end{figure}

In our approach we vary  the spatial directions estimating the vector field $\vbold$ as described in \cref{eq: structure tensor}--\cref{eq: compute b vary step 2} while we fix $\beta(\cdot)$, so as to fix the elliptic shape of the test functions. 
First, in \cref{fig: bamboo first sensitive eta} we report the sensitivity to the parameter $\eta\in[0.24,4]$ (with step-size increment of $0.25$) for the first order $\TDVM{}{}$ regulariser: from our experiments $\eta=0.75$ produces a better result than the first order regularisers, $\TV{}{}$ in \cref{fig: bamboo TV} and $\DTV{}{}$ in \cref{fig: bamboo DTV}, and the second order $\TGV{}{2}$ in \cref{fig: bamboo TGV}, showing less staircasing and directional artefacts. 
In this test, we fixed a-priori a number of parameters, i.e.\ $\bbold=(1,0.02)$, $(\sigma_1,\rho_1)=(1.8,2.8)$ and $(\sigma_2,\rho_2)=(1,1)$ for the anisotropic structure in building $\Mbold$ and $\texttt{maxiter}=2$ in \cref{alg: main denoising}. 

In \cref{fig: bamboo all orders} we report the best results obtained with the same fixed choice of parameters but now with all the possible combinations of first, second and third order regularisers, as well as the sketch of the streamlines of $\vbold$ in \cref{fig: v bamboo}.
We observed that the combination of the first and third order of regularisers outperforms the results from \cref{fig: bamboo denoising current results}.
In the next paragraph we comment about further choices for the parameters.

\paragraph{Selection of parameters for the bamboo image}
Let the choice of the regulariser orders $\abold$ and the $\texttt{maxiter}$ in \cref{alg: main denoising} be fixed. 
We would like to estimate good parameters that directly affect the image reconstruction, namely the fidelity $\eta$, the directional information provided by $(\sigma_1,\rho_1)$ of the structure tensor and $\beta$ from $\bbold=(1,\beta)$ for the anisotropic shape of the test functions (assumed for now fixed all over the imaging domain). 
Unfortunately, the structure tensor depends on the intrinsic image content and standard deviation of the noise.

\begin{figure}[tbph]
\centering
\begin{subfigure}[t]{0.23\textwidth}\centering\captionsetup{justification=centering}
\includegraphics[width=0.95\textwidth]{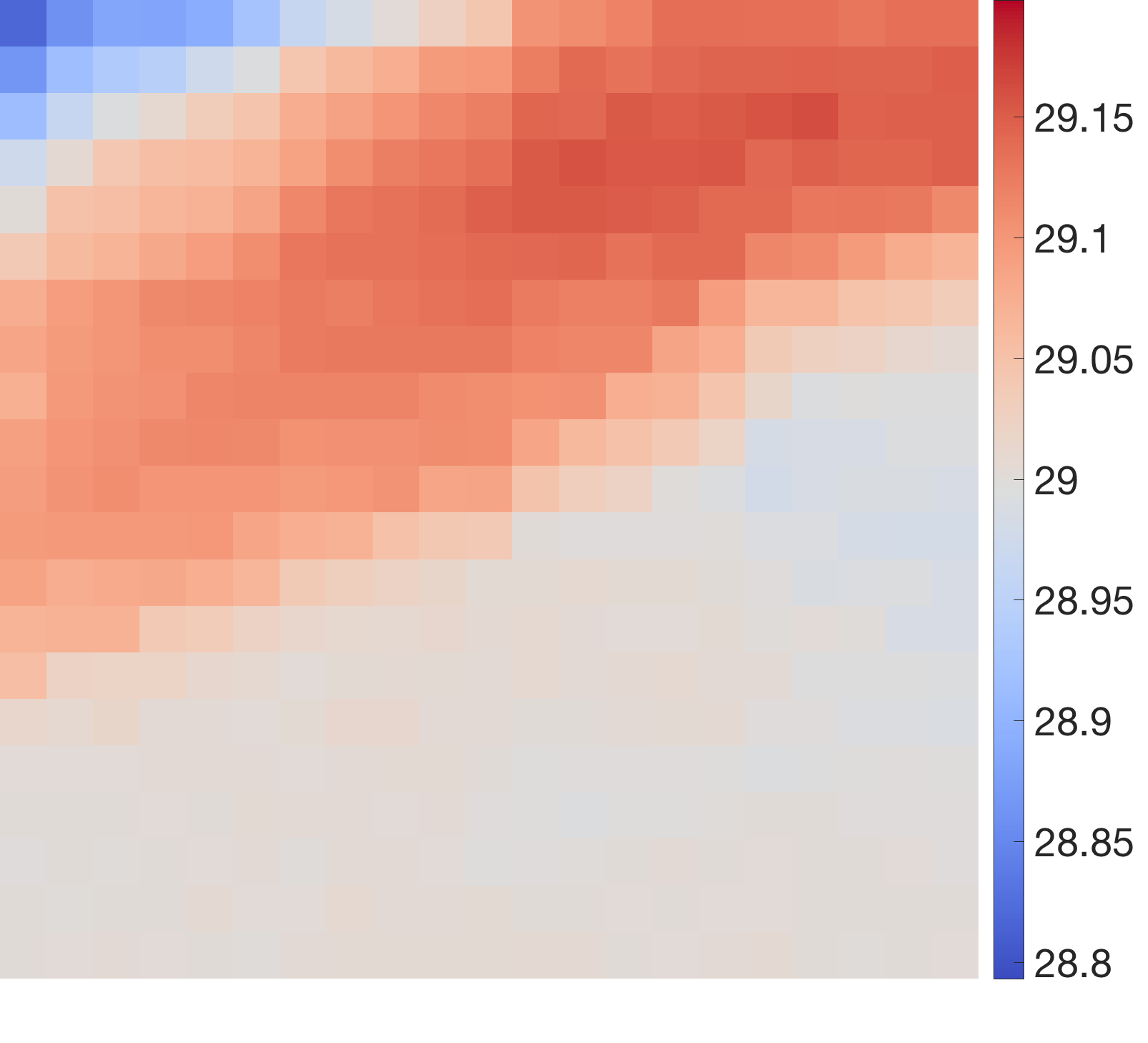}
\caption{$\bbold=(1,0)$}
\label{fig: imagesc denoising 000}
\end{subfigure}
\,
\begin{subfigure}[t]{0.235\textwidth}\centering\captionsetup{justification=centering}
\includegraphics[width=0.95\textwidth]{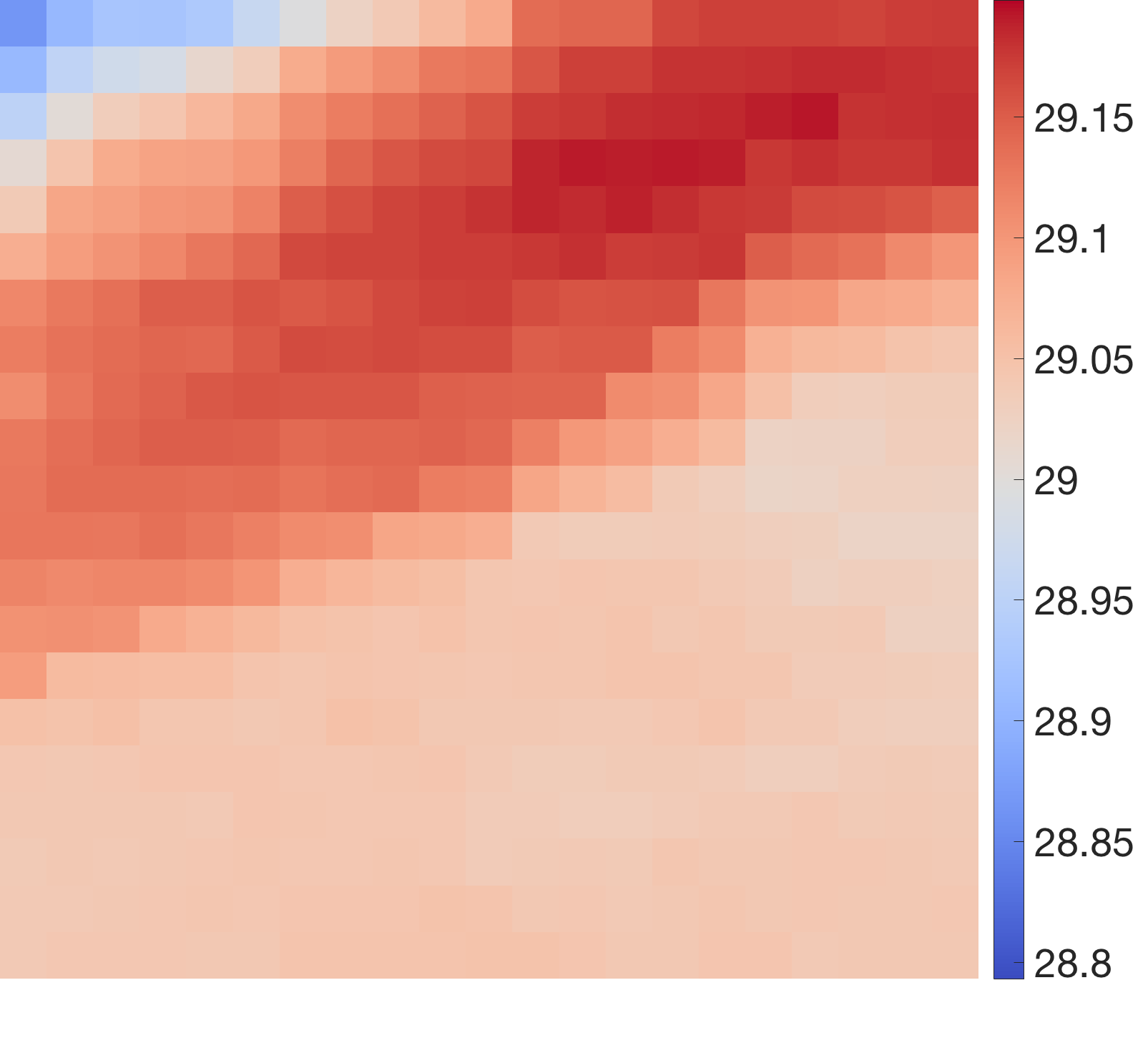}
\caption{$\bbold=(1,0.01)$}
\label{fig: imagesc denoising 001}
\end{subfigure}
\,
\begin{subfigure}[t]{0.235\textwidth}\centering\captionsetup{justification=centering}
\includegraphics[width=0.95\textwidth]{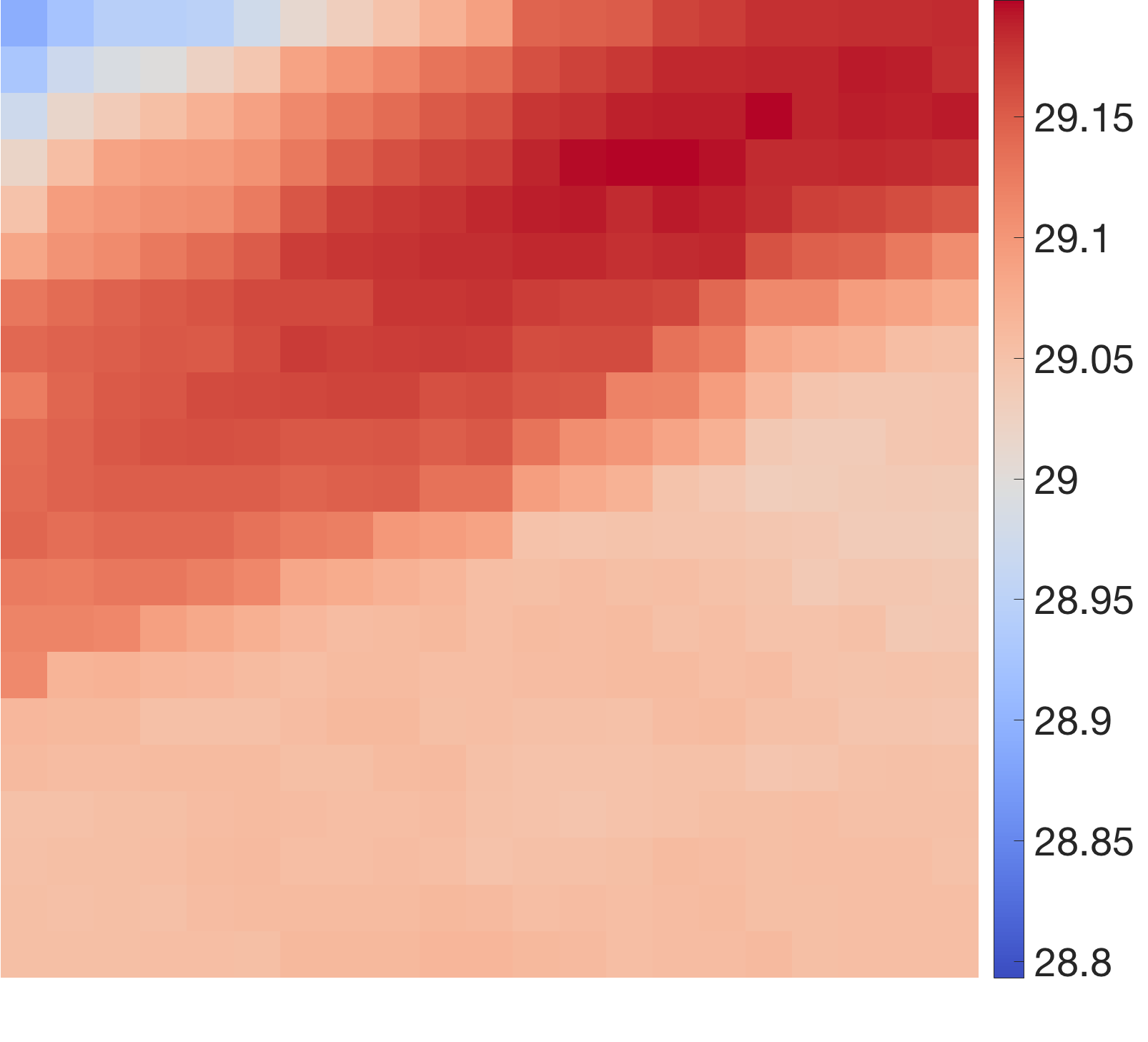}
\caption{$\bbold=(1,0.02)$}
\label{fig: imagesc denoising 002}
\end{subfigure}
\,
\begin{subfigure}[t]{0.235\textwidth}\centering\captionsetup{justification=centering}
\includegraphics[width=0.95\textwidth]{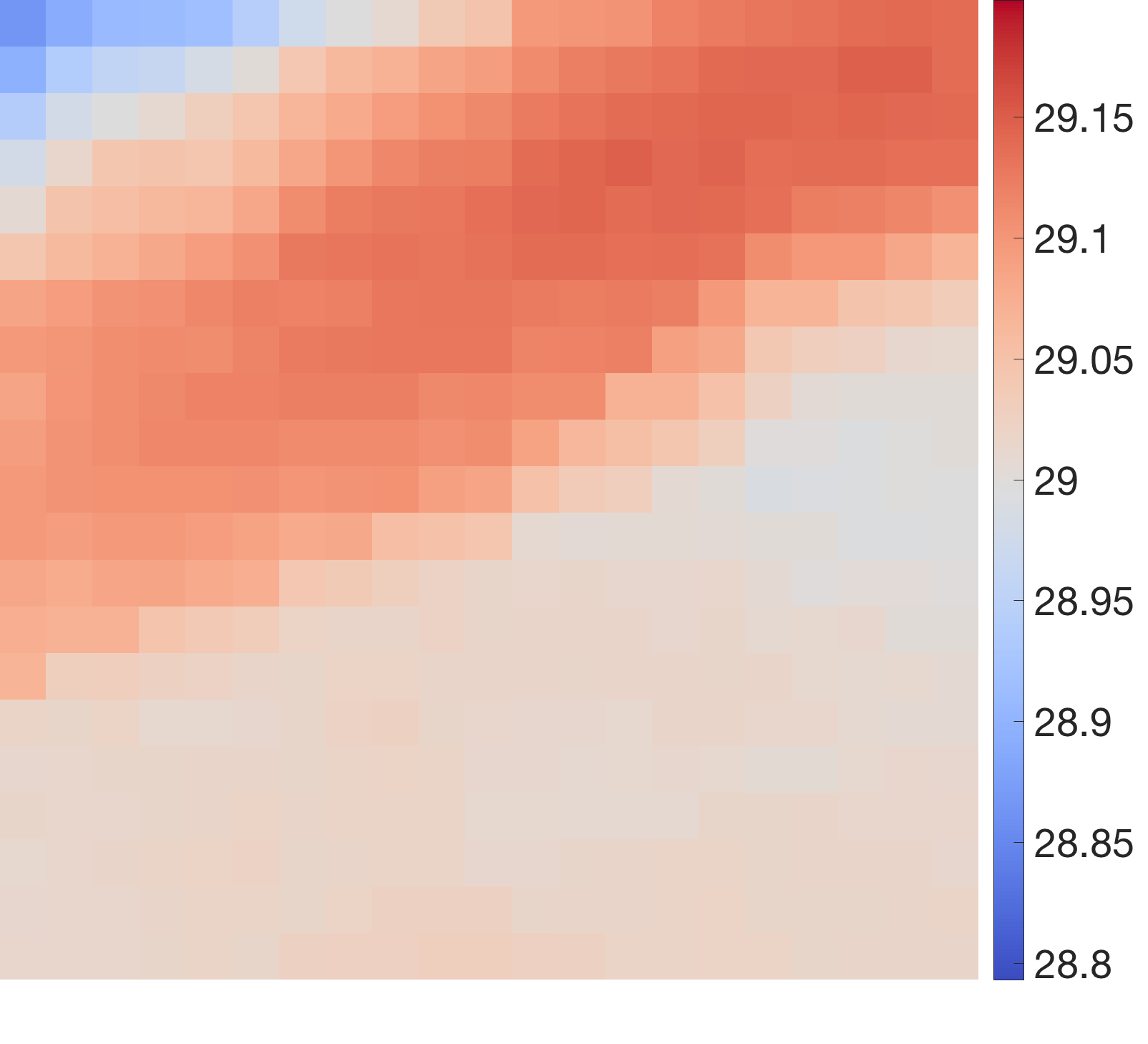}
\caption{$\bbold=(1,0.03)$}
\label{fig: imagesc denoising 003}
\end{subfigure}
\vspace{-1em}
\caption{PSNR for different $\sigma_1$ (y-axis), $\rho_1$ (x-axis) and $\bbold$ for \cref{fig: bamboo 20noise} with $\TDVM{1}{\alphabold_1}+\TDVM{3}{\alphabold_3}$, $\alphabold_1=(1)$, $\alphabold_3=(1,\infty,\infty)$ and $\eta=3.5$.  Here $k=2$ for  Alg.\ \ref{alg: main denoising}, with $(\sigma_2,\rho_2)=(1,1)$.}
\label{fig: denoising bamboo consistency test}
\vspace{-1.5em}
\end{figure}

\begin{figure}[tb]
\centering
\begin{subfigure}{0.325\textwidth}\centering
\captionsetup{justification=centering}
\includegraphics[width=0.9\textwidth]{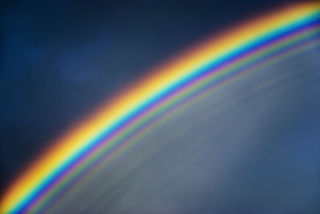}
\caption{Original $\ubold$, $320\times 214$ pixels}
\label{fig: rainbow original}
\end{subfigure}
\hfill
\begin{subfigure}{0.325\textwidth}\centering\captionsetup{justification=centering}
\includegraphics[width=0.9\textwidth]{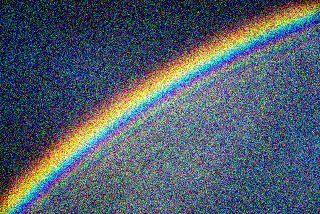}
\caption{Noisy $\ubold^\diamond$ ($20\%$ Gaussian)}
\label{fig: rainbow 20noise}
\end{subfigure}
\hfill
\begin{subfigure}{0.325\textwidth}\centering\captionsetup{justification=centering}
\includegraphics[width=0.9\textwidth]{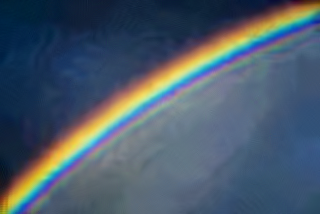}
\caption{BM3D, PSNR = 34.53}
\label{fig: BM3D rainbow}
\end{subfigure}
\\
\begin{subfigure}{0.325\textwidth}\centering\captionsetup{justification=centering}
\includegraphics[width=0.9\textwidth]{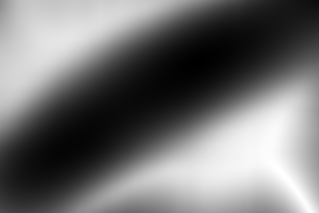}
\caption{$\bbold_2, (\sigma, \rho) = (2, 30)$,\\ from equation \cref{eq: compute b vary step 2}}
\label{fig: b rainbow}
\end{subfigure}
\hfill
\begin{subfigure}{0.325\textwidth}\centering\captionsetup{justification=centering}
\includegraphics[width=0.9\textwidth]{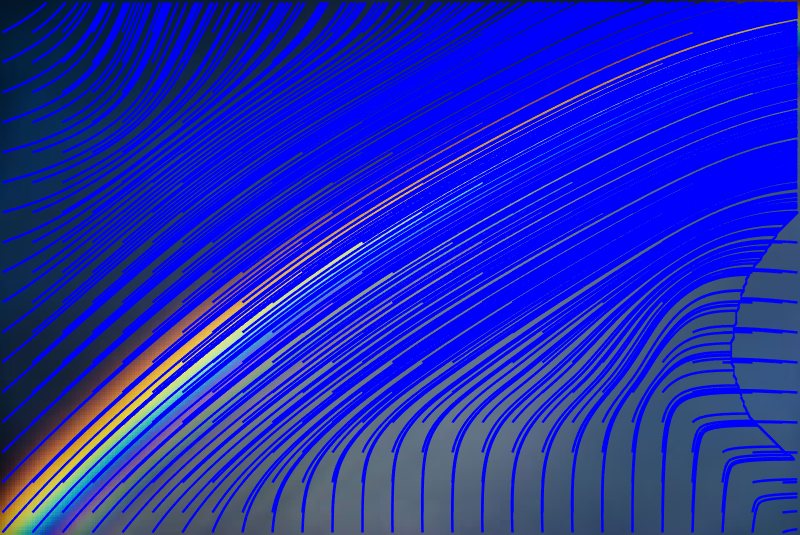}
\caption{streamlines of $\vbold$\\
$(\sigma, \rho) = (2, 30)$ from equation \cref{eq: find v denoising}}
\label{fig: v rainbow}
\end{subfigure}
\hfill
\begin{subfigure}{0.325\textwidth}\centering\captionsetup{justification=centering}
\includegraphics[width=0.9\textwidth]{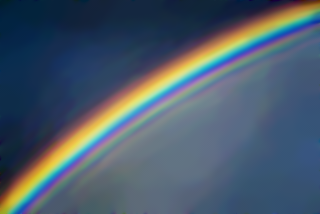}
\caption{
$\TDVM{1}{}+\TDVM{3}{}$, PSNR = 35.91\\
$\eta = 3.5$; $\bbold$, $\vbold$ from Figs.\ \ref{fig: b rainbow}-\ref{fig: v rainbow}.
}
\label{fig: u rainbow}
\end{subfigure}
\vspace{-1em}
\caption{Denoising of RGB rainbow (photo by M.P.\ Markkanen, CC-BY-SA-4.0 license). }
\label{fig: denoising results color}
\vspace{-1.5em}
\end{figure}

\begin{figure}[tbhp]
\begin{subfigure}[t]{0.16\textwidth}
\captionsetup{justification=centering}
\includegraphics[width=0.92\textwidth]{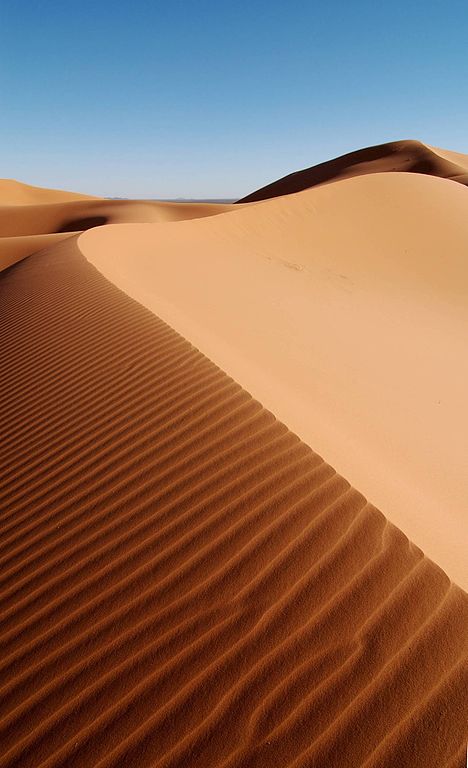}
\caption{Original $\ubold$ \\ $468\times 768$ pixels}
\label{fig: desert original}
\end{subfigure}
\hfill
\begin{subfigure}[t]{0.16\textwidth}
\captionsetup{justification=centering}
\includegraphics[width=0.92\textwidth]{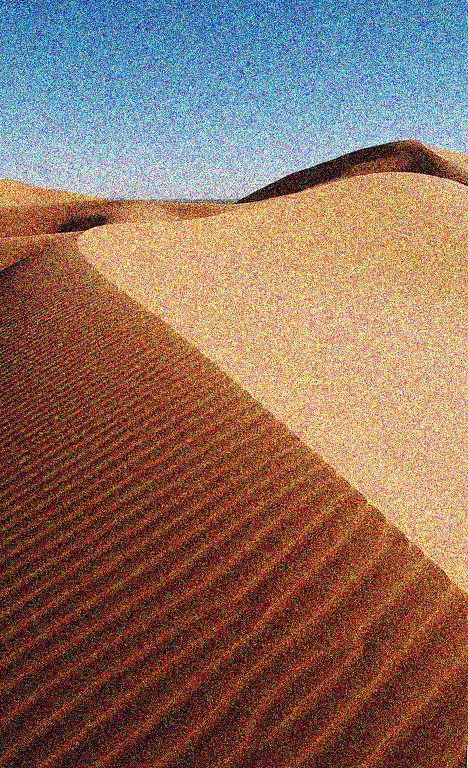}
\caption{Noisy $\ubold^\diamond$\\($20\%$ Gaussian)}
\label{fig: desert 20noise}
\end{subfigure}
\hfill
\begin{subfigure}[t]{0.16\textwidth}
\captionsetup{justification=centering}
\includegraphics[width=0.92\textwidth]{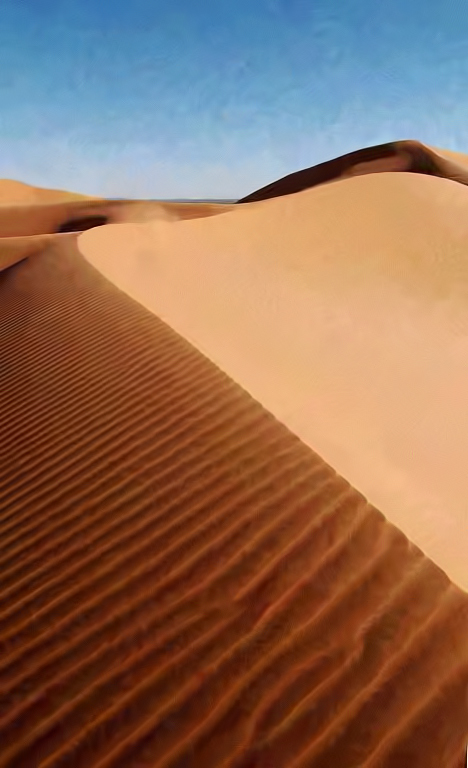}
\caption{BM3D\\PSNR = 32.45}
\label{fig: BM3D desert}
\end{subfigure}
\hfill
\begin{subfigure}[t]{0.16\textwidth}
\captionsetup{justification=centering}
\includegraphics[width=0.92\textwidth]{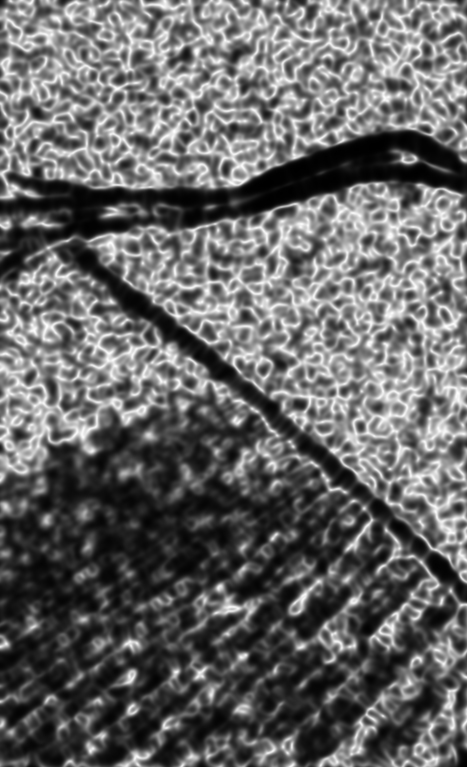}
\caption{$\bbold_2$,\\$(\sigma, \rho) = (3, 1.5)$, \\ from eq.\ \cref{eq: compute b vary step 2}}
\label{fig: b desert}
\end{subfigure}
\hfill
\begin{subfigure}[t]{0.16\textwidth}
\captionsetup{justification=centering}
\includegraphics[width=0.92\textwidth]{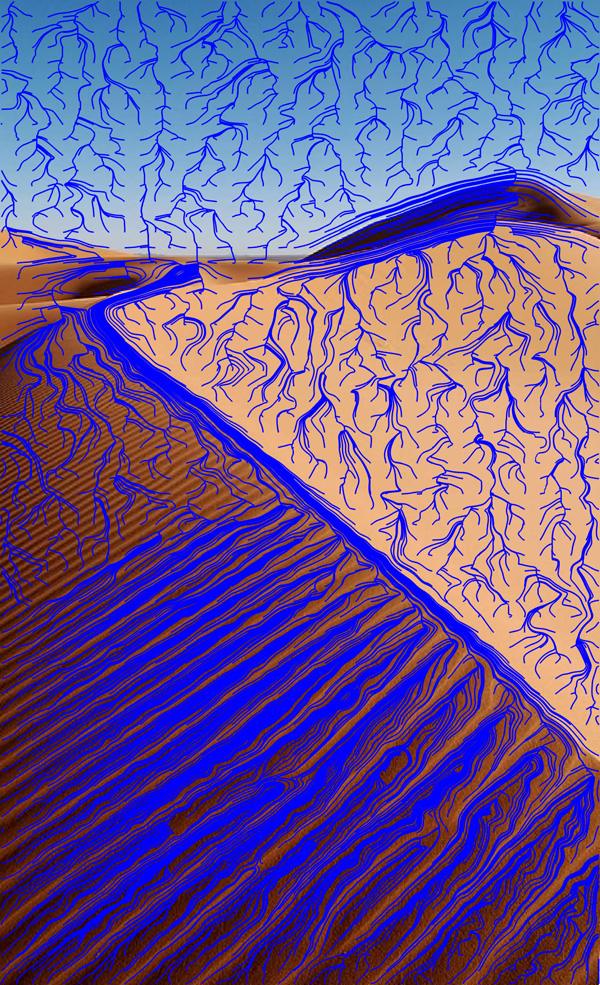}
\caption{$\vbold$ streamlines,\\$(\sigma, \rho) = (3, 1.5)$, \\ from eq. \cref{eq: find v denoising}}
\label{fig: v desert}
\end{subfigure}
\hfill
\begin{subfigure}[t]{0.16\textwidth}
\captionsetup{justification=centering}
\includegraphics[width=0.92\textwidth]{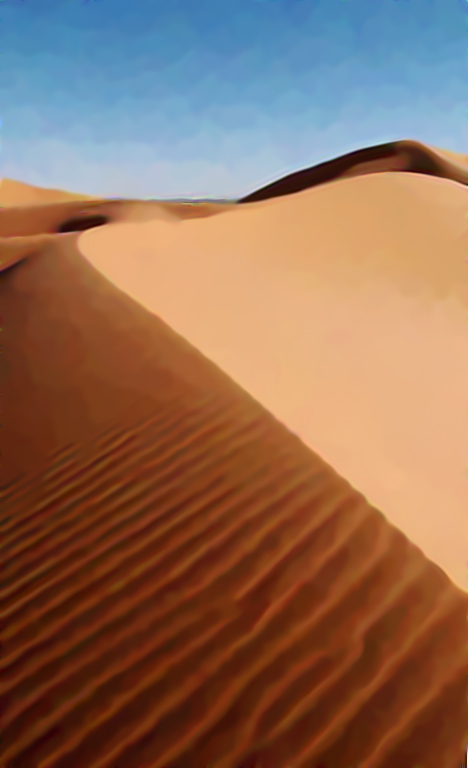}
\caption{
$\TDVM{1}{}+\TDVM{3}{}$\\ 
$\eta=3.5$\\ 
PSNR = 30.09 
}
\label{fig: u desert}
\end{subfigure}
\\
\begin{subfigure}[t]{0.495\textwidth}
\centering
\captionsetup{justification=centering}
\includegraphics[width=0.97\textwidth,trim=0 1.7cm 0 3cm, clip=true]{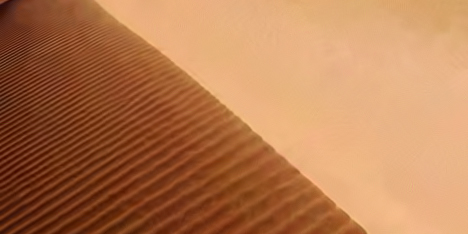}
\caption{
Zoom of BM3D in \cref{fig: BM3D desert}.
}
\label{fig: BM3D zoom desert}
\end{subfigure}
\hfill
\begin{subfigure}[t]{0.495\textwidth}
\centering
\captionsetup{justification=centering}
\includegraphics[width=0.97\textwidth,trim=0 1.7cm 0 3cm, clip=true]{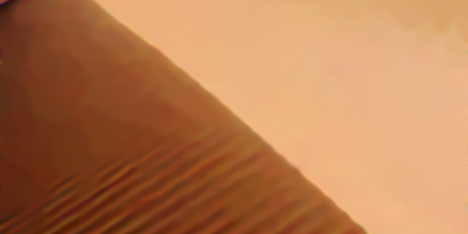}
\caption{
Zoom of our $\TDVM{1}{}+\TDVM{3}{}$ in \cref{fig: u desert}.
}
\label{fig: TDV zoom desert}
\end{subfigure}
\vspace{-1em}
\caption{Denoising of RGB desert (photo by Rosino, CC-BY-SA-2.0 license).
}
\label{fig: denoising desert}
\vspace{-1.5em}
\end{figure}

However, for $\eta=3.5$ (among the $\eta\in[0.25,4]$ tested) we obtained the best result for the combination of the first and third order regularisers, i.e.\ $\abold=(1,0,1)$. 
Therefore we inspect more the PSNR obtained by fixing $\eta,\abold$ and by changing both $\beta\in\{0,0.01,0.02,0.03\}$ and $(\sigma_1,\rho_1)$ in the range between $(1.5,3.5)$: results in \cref{fig: denoising bamboo consistency test} show that the optimal parameters for this case are the ones producing \cref{fig: bamboo 101} as output.
Other strategies for estimating the parameters, including $\eta$, would require to solve a bilevel problem, e.g.\ as in \cite{DelosReyes2017}, or to solve the problem by updating the parameters with a greedy line-search approach onto many directional images, so as to extract a rule of thumb for the choice,  e.g.\ as in \cite{ParSch2019}.

A natural question now is whether allowing also $\beta(\cdot)$ to change in $\Omega$ improves the performances. We are going to answer this question in the next two experiments.
However, since the scale of the directional texture in the next experiments are spatially different, we fix $\eta=3.5$ and focus on the promising case of combination of first and third order regularisers, i.e.\ $\abold=(1,0,1)$, with reasonable choice of $(\sigma_1,\rho_1)$. \paragraph{Rainbow image}
The rainbow in \cref{fig: rainbow original} has been corrupted by $20\%$ of Gaussian noise in each color channel, see \cref{fig: rainbow 20noise}. 
Due to the particular structure of the image, an isotropic approach  seems reasonable outside the rainbow while an anisotropic approach inside. This resulted in varying the $\beta$ parameter following equations \cref{eq: compute b vary step 1}-\cref{eq: compute b vary step 2}: in \cref{fig: b rainbow} the black pixels correspond to $\beta\approx 0$ and the white pixels to $\beta\approx 1$. Indeed for $\beta\approx 0$ we expect to denoise the image following the anisotropy induced by $\vbold$, while for $\beta\approx 1$ we expect to denoise the image isotropically in both $\vbold$ and $\vbold_\perp$ directions. 
In order to compute the vector field $\vbold$ as in \cref{fig: v rainbow}, we did not apply the regularisation step \cref{eq: find v denoising} and we did not iterate \cref{alg: main denoising} since both the resulting $\vbold$ and $\bbold$ seem good enough for our purposes, performing better than BM3D in \cref{fig: BM3D rainbow}, with less wavy artefacts and a smoother global structure.

 \paragraph{Desert image}
The desert image in \cref{fig: desert original} is a mix of anisotropic and isotropic information. 
We denoised \cref{fig: desert 20noise}, corrupted again with $20\%$ of Gaussian noise in each color channel, with $(\sigma,\rho)=(3,1.5)$ and $\gamma=0.1$ to estimate $\vbold$ in \cref{fig: v desert} as described before. We also allowed $\beta(\cdot)$ to vary across the domain and we did not refine $\vbold$ with further iterations. Here, BM3D in \cref{fig: BM3D desert} performed slightly better than our approach due to the wrong estimation of $\vbold$ along some dune waves: this is clearly visible in Figures \ref{fig: b desert} and \ref{fig: v desert} where both the wrong estimation of $\vbold$ and the isotropy requirement of $\beta$ result in a smoothing performance, as shown in the zooms of Figures \ref{fig: BM3D zoom desert} and \ref{fig: TDV zoom desert}. However, we recall that BM3D is a non-local method and better performances than local methods are expected. 
Nevertheless, since this kind of images have patterns and structures at different scales such that they cannot be captured by a single global structure tensor, we expect that the reconstruction quality can be improved, e.g.\ by a structure tensor workflow with locally adapting parameters $(\sigma,\rho)$.

\section{Other imaging applications with the joint model}
In what follows we focus on the joint minimisation model \cref{eq: minimization intro} for the applications of image zooming and surface interpolation.
\label{sec: imaging applications}
\subsection{Wavelet Zooming}\label{sec: model-description-wavelet}
\begin{figure}[!tb]
\begin{subfigure}[t]{0.24\textwidth}\centering\captionsetup{justification=centering}
\includegraphics[width=0.88\textwidth]{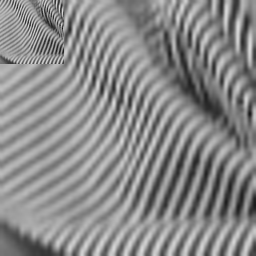}
\caption{4X Lanczos 2 filter}
\label{fig: wavelet box}
\end{subfigure}
\hfill
\begin{subfigure}[t]{0.24\textwidth}\centering\captionsetup{justification=centering}
\includegraphics[width=0.88\textwidth]{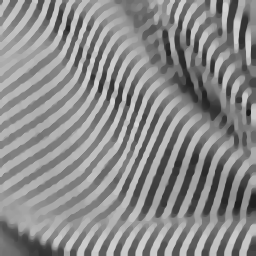}
\caption{$\TGV{}{2}$, \cite{Bredies2013}}
\label{fig: wavelet TGV}
\end{subfigure}
\hfill
\begin{subfigure}[t]{0.24\textwidth}\centering\captionsetup{justification=centering}
\includegraphics[width=0.88\textwidth]{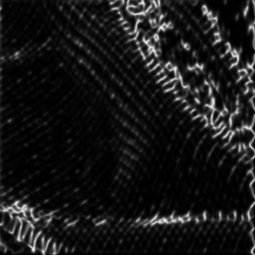}
\caption{$\bbold_2$, $(\sigma, \rho) = (1, 1)$}
\end{subfigure}
\hfill
\begin{subfigure}[t]{0.24\textwidth}\centering\captionsetup{justification=centering}
\includegraphics[width=0.88\textwidth]{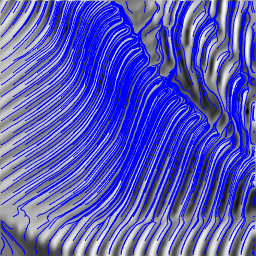}
\caption{$\vbold$, $(\sigma, \rho) = (1, 1)$}
\end{subfigure}
\\
\begin{subfigure}[t]{0.24\textwidth}\centering\captionsetup{justification=centering}
\includegraphics[width=0.88\textwidth]{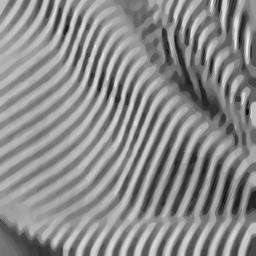}
\caption{$\TDVM{1}{}$, $\bbold_2$, $\vbold$}
\end{subfigure}
\hfill
\begin{subfigure}[t]{0.24\textwidth}\centering\captionsetup{justification=centering}
\includegraphics[width=0.88\textwidth]{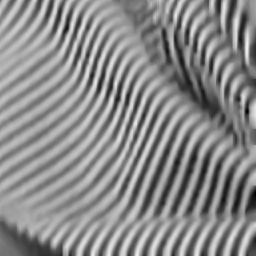}
\caption{$\TDVM{2}{}$, $\bbold_2$, $\vbold$}
\end{subfigure}
\hfill
\begin{subfigure}[t]{0.24\textwidth}\centering\captionsetup{justification=centering}
\includegraphics[width=0.88\textwidth]{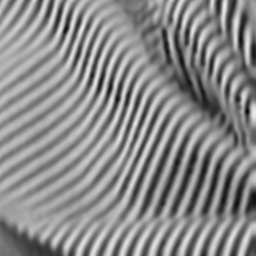}
\caption{$\TDVM{3}{}$, $\bbold_2$, $\vbold$}
\end{subfigure}
\hfill
\begin{subfigure}[t]{0.24\textwidth}\centering\captionsetup{justification=centering}
\includegraphics[width=0.88\textwidth]{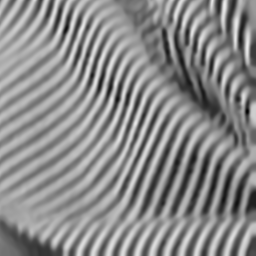}
\caption{$\TDVM{1}{}+\TDVM{3}{}$, $\bbold_2$, $\vbold$}
\end{subfigure}
\vspace{-1em}
\caption{Wavelet-based zooming with CDF 9/7.}
\label{fig: wavelet2 zooming}
\end{figure}

\begin{figure}[!tb]\centering
\begin{subfigure}[t]{0.146\textwidth}\centering\captionsetup{justification=centering}
\includegraphics[width=1\textwidth]{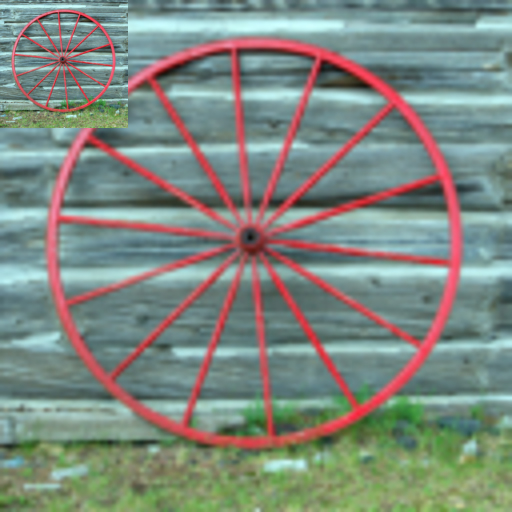}
\caption{4X Bicubic interpolation}
\label{fig: wavelet3 bicubic}
\end{subfigure}
\hfill
\begin{subfigure}[t]{0.146\textwidth}\centering\captionsetup{justification=centering}
\includegraphics[width=1\textwidth]{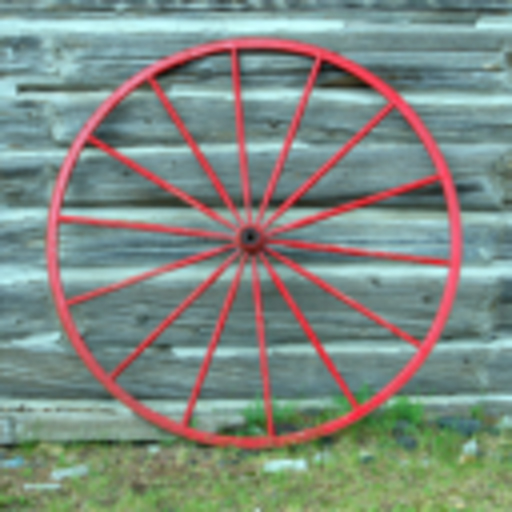}
\caption{4X Lanczos 2}
\label{fig: wavelet3 lanczos}
\end{subfigure}
\hfill
\begin{subfigure}[t]{0.146\textwidth}\centering\captionsetup{justification=centering}
\includegraphics[width=1\textwidth]{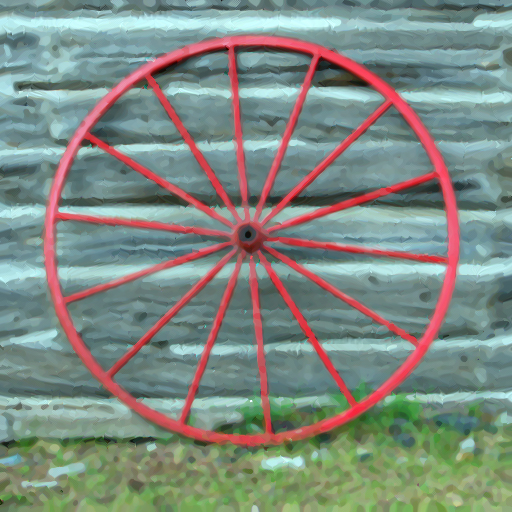}
\caption{$\TDVM{1}{}$}
\label{fig: waveletrgb100}
\end{subfigure}
\hfill
\begin{subfigure}[t]{0.146\textwidth}\centering\captionsetup{justification=centering}
\includegraphics[width=1\textwidth]{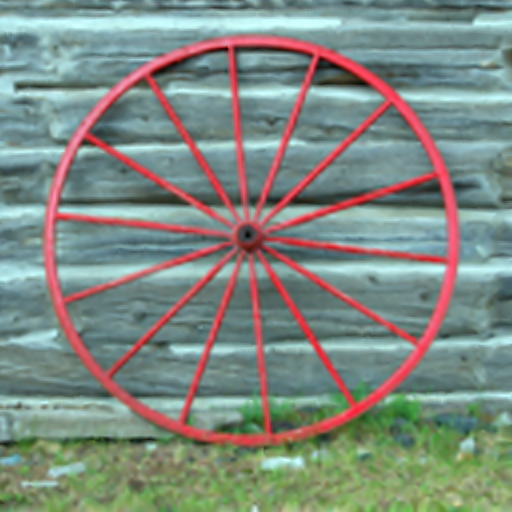}
\caption{$\TDVM{2}{}$}
\label{fig: waveletrgb010}
\end{subfigure}
\hfill
\begin{subfigure}[t]{0.146\textwidth}\centering\captionsetup{justification=centering}
\includegraphics[width=1\textwidth]{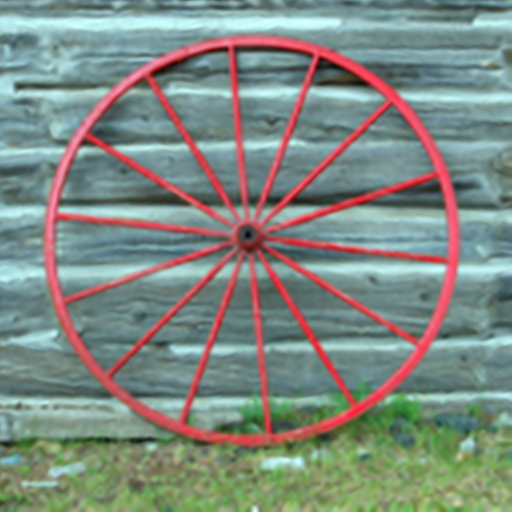}
\caption{$\TDVM{3}{}$}
\label{fig: waveletrgb001}
\end{subfigure}
\hfill
\begin{subfigure}[t]{0.146\textwidth}\centering\captionsetup{justification=centering}
\includegraphics[width=1\textwidth]{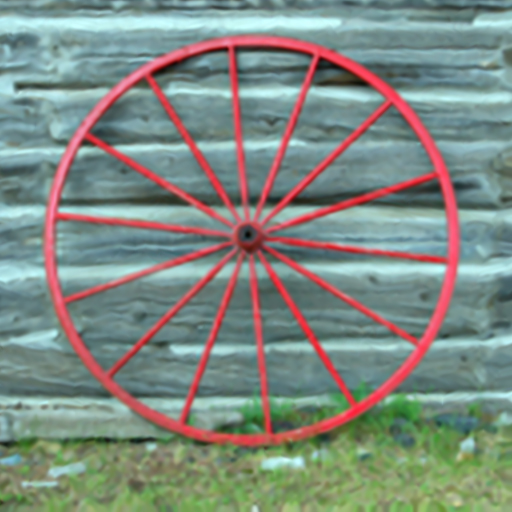}
\caption{$\TDVM{1}{}+\TDVM{3}{}$}
\label{fig: waveletrgb101}
\end{subfigure}
\\
\begin{subfigure}[t]{0.146\textwidth}\centering\captionsetup{justification=centering}
\includegraphics[width=1\textwidth,trim=5.85cm 6.15cm 6.15cm 6.15cm,clip=true]{{./images/wavelet/tyre4x/u_bicubic}.jpg}
\caption{Zoom of \cref{fig: wavelet3 bicubic}.}
\end{subfigure}
\hfill
\begin{subfigure}[t]{0.146\textwidth}\centering\captionsetup{justification=centering}
\includegraphics[width=1\textwidth,trim=5.85cm 6.15cm 6.15cm 6.15cm,clip=true]{{./images/wavelet/tyre4x/u_lanczos2}.jpg}
\caption{Zoom of \cref{fig: wavelet3 lanczos}.}
\end{subfigure}
\hfill
\begin{subfigure}[t]{0.146\textwidth}\centering\captionsetup{justification=centering}
\includegraphics[width=1\textwidth,trim=5.85cm 6.15cm 6.15cm 6.15cm,clip=true]{{./images/wavelet/tyre4x/4magification_bvary1_a1-0-0_sigma1.2_rho2.5_iter1_PSNR0}.jpg}
\caption{Zoom of \cref{fig: waveletrgb100}.}
\end{subfigure}
\hfill
\begin{subfigure}[t]{0.146\textwidth}\centering\captionsetup{justification=centering}
\includegraphics[width=1\textwidth,trim=5.85cm 6.15cm 6.15cm 6.15cm,clip=true]{{./images/wavelet/tyre4x/4magification_bvary1_a0-1-0_sigma1.2_rho2.5_iter1_PSNR0}.jpg}
\caption{Zoom of \cref{fig: waveletrgb010}.}
\end{subfigure}
\hfill
\begin{subfigure}[t]{0.146\textwidth}\centering\captionsetup{justification=centering}
\includegraphics[width=1\textwidth,trim=5.85cm 6.15cm 6.15cm 6.15cm,clip=true]{{./images/wavelet/tyre4x/4magification_bvary1_a0-0-1_sigma1.2_rho2.5_iter1_PSNR0}.jpg}
\caption{Zoom of \cref{fig: waveletrgb001}.}
\end{subfigure}
\hfill
\begin{subfigure}[t]{0.146\textwidth}\centering\captionsetup{justification=centering}
\includegraphics[width=1\textwidth,trim=5.85cm 6.15cm 6.15cm 6.15cm,clip=true]{{./images/wavelet/tyre4x/4magification_bvary1_a1-0-1_sigma1.2_rho2.5_iter1_PSNR0}.jpg}
\caption{Zoom of \cref{fig: waveletrgb101}.}
\end{subfigure}
\vspace{-1em}
\caption{
Wavelet-based zooming with CDF 9/7. Image from the MorgueFile archive.
}
\label{fig: wavelet3 zooming}
\vspace{-1.5em}
\end{figure}

In this section we apply our regularisation to wavelet-based image zooming as in \cite{Bredies2013}. 
Here, the data fidelity term is modelled by a wavelet transformation operator. 
Let $\phi\in\LL{2}[\RR]$, $\psi\in\LL{2}[\Omega]$ be the scaling and mother wavelet function, respectively. Then, a Riesz basis of $\LL{2}[\Omega]$ is obtained from translations and rotations of $\phi$ and $\psi$. Here, we will consider only functions $\phi$ with compact support, yielding a compactly supported basis elements. Let $R\in\ZZ$ be a resolution level and $M_R, (L_j)_{j\leq R}$ be finite index sets in $\ZZ^2$, then:
\begin{itemize}
\item a Riesz basis of $\LL{2}[(0,1)\times(0,1)]$ is $
(\phi_{R,\kbold})_{\kbold\in M_R}$, $(\psi_{j,\kbold})_{j\leq R,\kbold\in L_j}$;
\item  the dual Riesz basis of the above is defined as $(\widetilde{\phi}_{R,\kbold})_{\kbold\in M_R}$, $(\widetilde{\psi}_{j,\kbold})_{j\leq R,\kbold\in L_j}$.
\end{itemize}
Thus, the following decomposition holds:
\[
u = \sum_{\kbold\in M_R} (u,\phi_{R,\kbold})_2\widetilde{\phi}_{R,\kbold} + \sum_{j\leq R,\kbold\in L_j} (u,\psi_{j,\kbold})_2 \widetilde{\psi}_{j,\kbold}.
\]
Let $u_0\in\mathrm{span}\left\{\widetilde{\phi}_{R,\kbold}\,\lvert\, \kbold\in M_R\right\}$ be a low resolution version of $u$ given by $((u_0,\phi_{R,\kbold})_2)_{\kbold\in M_R}$, where the unknown $u$ is such that $
(u,\phi_{R,\kbold})_2 = (u_0,\phi_{R,\kbold})_2,\text{ for all }\kbold \in M_R,
$ and
\[
u\in\LL{2}[\Omega] 
= 
\mathrm{span}
\left( 
\left\{ \widetilde{\phi}_{R,\kbold}\,\lvert\, \kbold\in M_R \}\cup \{ \widetilde{\psi}_{j,\kbold}\,\lvert\, j\leq R,\kbold\in L_j
\right\}
\right).
\]
The wavelet-based zooming problem with higher order total directional regularisers reads as
\[
u^\star \in \argmin_{u\in\RR^{|\Omega^h|}} \sum_{q=1}^\Qrm \TDVM{q}{\alphabold_q}[u,\Mcal] + \Ical_{U_D}(u),
\]
where
$
U_D = \left\{u\in\LL{2}[\Omega]\,\lvert\, (u,\phi_{R,\kbold})_2 = (u_0,\phi_{R,\kbold})_2,\text{ for all }\kbold \in M_R\right\}
$ and $\Ical_{U_D}$ is the convex indicator function w.r.t.\ $U_D$, see \cite{Bredies2013} for more details. Since we did not downsample the original image, we avoided artefacts introduced by algorithms for reducing the image but at the same time no ground truth is available. Results based on CDF 9/7 wavelet are shown in \cref{fig: wavelet2 zooming} (grey-scale) and in \cref{fig: wavelet3 zooming} (colour). 

 \subsection{Surface Interpolation}\label{sec: model-description-surface}

In this experiment, we aim to reconstruct a surface from scattered height data available in $\Omega$. 
The available data lies on partially occluded isolines or on random points in $\Omega$ and the challenge is to interpolate them 
by
preserving the anisotropic features via the reconstruction of a suitable vector field $\vbold$.
Before presenting our approach for this problem using $\TDVM{q}{\alphabold}$, we briefly review the state-of-the-art for surface interpolation.

\subsubsection*{Related works}
The reconstruction of surfaces from scattered height values has been approached in two different ways in the literature: based on explicit and implicit models. 
Surface interpolation is sometimes also addressed as \emph{digital elevation map} (DEM) problem.

In this paper we focus on implicit surface interpolation which has the advantage of being independent with respect to parametrization. 
Here the surface is an implicit function of height values over the domain.
Two prominent methods in this range are the Thin Plate Spline (TPS) \cite{Mei1984} and the \emph{Absolute Minimizing Lipschitz Extension} (AMLE) \cite{almansa} approach. 
TPS is a flexible approach since it can embed both grey values and gradient information. However, it has the drawback to be a fourth order isotropic method and the resulting interpolated surface is isotropically smooth. AMLE, on the other hand, is able to interpolate data given in isolated points and on curves but it fails to interpolate slopes of a surface, resulting in $\CCCspace{1}$, see \cite{Sav2005}.

For interpolating surfaces with sharp features, e.g.\ strong creases, and possibly non-smooth features, e.g.\ corners in a pyramid, it seems promising therefore to consider (higher-order) total variation ($\TV{}$) regularisers for surface interpolation.

Our main model approach here is  \cite{LelMorSch2013}, where a third-order directional total variation regulariser has been proposed that reads for a given vector field $\vbold$ as
\begin{equation}
E(u) = \int_\Omega \norm{\grad_\vbold(\grad^2u)}_2,
\label{eq: third order directional jan}
\end{equation}
where $\grad_\vbold^3u = \grad(\grad^2 u)\cdot \vbold$
is the directional derivative of the Hessian of $u$, along $\vbold$. 
Note that this is a special case of $\TDVM{\Qrm}{\alphabold}$ with $\Qrm=3$, $a=1$, $\bbold=(1,0)$, $\vbold = (\cos\thetabold,\sin\thetabold)$, i.e.\ $\Mcal=(\Ibold,\Ibold,\Mbold)$ and $\Mbold =\Lambdabold_\bbold (\Rbold_\thetabold)^\T$, leading to $\grad_\vbold^3u\equiv\Mbold\grad\otimes(\grad^2 u)$.

The estimation of $\vbold$ is crucial to obtain a good quality result. In \cite{LelMorSch2013}, $\vbold$ has been computed as a two step minimization-regularisation problem by solving firstly
\[
\widetilde{\vbold} = \argmin_{ \substack{\ybold\,:\,\norm{\ybold}_2 = 1 }} \norm{ K_\sigma \ast \grad \left( \frac{\grad u}{\abs{\grad u}}\right)(\xbold) \,\ybold}_2
\]
and then applying to $\widetilde{\vbold}$ the same regularisation step in \cref{eq: find v denoising},
where $w(\xbold)$ is a weight chosen as the largest singular value of $K_\sigma \ast \grad \left( \frac{\grad u}{\abs{\grad u}}\right)$ and $\gamma$ is a regularising parameter smoothing the vector field where $u$ is almost planar and preserving the local values of $\widetilde{\vbold}$ for level lines of large curvature. 
As last step, $\vbold$ is normalised to be unitary.

Another directional interpolation model for $u$ and $\vbold$ appears in \cite{reportBinWu}: differently from our approach in this paper, it requires knowledge of the vector field $\vbold$ prior to the interpolation.

In this section, we generalize the approach of \cite{LelMorSch2013} for the reconstruction of a surface, given scattered height values lying (possibly) on partial contour lines.  Differently from \cref{sec: model-description-denoising}, the unitary vector field $\vbold$ is reconstructed in the missing domain as follows.

Let $\Omega^h$ be a 2D domain ($d=2$)
and $u^{\diamond,h}$ be sparse sampled height values. 
In the following, the projection onto the data available $u^{\diamond,h}$ is identified by the operator $\Scal$. 
We aim to find the interpolated surface $u\in\RR^{|\Omega^h|}$, driven by the unitary directions $\vbold\in\RR^{|\Gamma^h|}$. Let $\Mcal=(\Mcal_1,\dots,\Mcal_\Qrm)$ be a collection of weighting fields, where for a fixed $q$ each collection $\Mcal_q$ is defined as $\Mcal_q=(\Ibold,\dots,\Ibold,\Mbold)$  with explicit dependence on $\vbold$. 
We solve by \cref{alg: outer_min}, alternatingly:
\begin{align}
u^\star 
&\in 
\argmin_{u\in\RR^{|\Omega^h|}} \sum_{q=1}^\Qrm \TDVM{q}{\alphabold_q}[u,\Mcal_q^a(\vbold)] + \frac{\eta}{2}\norm{\Scal u - u^\diamond}_2^2,
\label{eq: surface minimization energy u}
\\
\vbold^\star 
&\in 
\argmin_{\vbold\in\RR^{|\Gamma^h|}} \mu \TV{}[\vbold] + \zeta \int_\Omega \left(1-\vbold \cdot \frac{\grad u}{\abs{\grad u}}\right)^2 \diff\xbold,
\label{eq: surface minimization energy v}
\end{align}
with the primal-dual in \cref{alg: reduced main primal-dual} for \cref{eq: surface minimization energy u} and a classic primal-dual for \cref{eq: surface minimization energy v}. 
In particular, in \cref{eq: surface minimization energy v} we identify $F(\vbold)=\TV{}[\vbold]$ for regularising the vector field $\vbold$ and $G(\vbold)= \norm{1-\vbold \cdot \frac{\grad u}{\abs{\grad u}}}_{\LL{2}}^2$ for normalising $\vbold$ in the direction of the normalised gradient \cite{BalBerCasSapVer2001}.

\begin{algorithm}[H]
\small
\caption{Alternating scheme for the surface interpolation problem in \cref{eq: surface minimization energy u,eq: surface minimization energy v}}
\label{alg: outer_min}
\SetKwData{maxiter}{maxiter}
\SetKwInOut{Input}{Input}
\SetKwInOut{Parameters}{Parameters}
\SetKwInOut{Identification}{Identification}
\SetKwInOut{Initialization}{Initialization}
\SetKwInOut{Update}{Update}
\Input{the sparse data $u^\diamond$ in $\Omega$.}
\Parameters{$\alphabold_q=(\alpha_j)_{j=0}^{q-1}$ for $q=1,\dots,\Qrm$ with $\alpha_0=1$ and $+\infty$ otherwise, $\eta>0$, $\mu,\zeta>0$.}
\Initialization{random $u^0$ and $\vbold^0$, $\Scal(u^0) = u^\diamond$, $t=0$.}
\Update{$u^t$ and $\vbold^{t}$ as follows.}
\While{\texttt{stopping criterion is not satisfied}}
{
$
\displaystyle
u^{t+1} \in\argmin_{u\in\RR^{|\Omega^h|}}\, \sum_{q=1}^\Qrm   \TDVM{q}{\alphabold_q}[u,\Mcal_q(\vbold^{t})]  + \frac{\eta}{2}\norm{\Scal u -u^\diamond}^2$\tcp*{Minimization w.r.t.\ $u$}
$\displaystyle
\vbold^{t+1} \in\argmin_{\vbold\in\RR^{|\Gamma^h|}} \mu \TV{}[\vbold] + \zeta \int_\Omega \left(1-\vbold \cdot \frac{\grad u^{t+1}}{\abs{\grad u^{t+1}}}\right)^2 \diff\xbold$\tcp*{Minimization w.r.t.\ $\vbold$}
Update $\Mcal_q(\vbold^{t+1})$ for each $q=1,\dots,\Qrm$\tcp*{update the weights}
}
\Return $(u^{\star},\vbold^\star) = (u^{t+1},\,\vbold^{t+1})$.
\end{algorithm}

\subsubsection*{Minimization with respect to \texorpdfstring{$u$}{u}}
Fixing an unitary vector field $\vbold^t$, the minimization problem \cref{eq: surface minimization energy u} is convex with respect to $u$ and the minimization problem 
can be solved via primal-dual \cref{alg: reduced main primal-dual} without acceleration due to the lack of strong convexity of the projection map $\Scal$, which results in a non-smooth dual problem.

\subsubsection*{Minimization with respect to \texorpdfstring{$\vbold$}{v}}
Fixing $u^{t+1}$, the minimization problem \cref{eq: surface minimization energy v} can be solved by the primal-dual algorithm with 
\[
F(\vbold) = \mu \TV{}[\vbold]\quad\text{and}\quad
G(\vbold) = \zeta \int_\Omega \left(1-\vbold \cdot \frac{\grad u^{t+1}}{\abs{\grad u^{t+1}}}\right)^2 \diff\xbold.
\]

Let $\sbold = K\vbold$, $K=\grad$ and $K^\ast=-\div$. Then, the proximal of $F^\ast$, with $F(\vbold)=\mu \TV{}[\vbold]$, is the projection onto the polar ball:
\[
\prox_{\sigma F^\ast}(\sbold) 
=
\frac{\sbold}{\max\left(1, \mu^{-1}\norm{\sbold}_2\right)}.
\]

The proximal map of $G$ at $\widehat{\vbold}=(\widehat{v}_1,\widehat{v}_2)$, for $\pbold=\grad u^{t+1} / \abs{\grad u^{t+1}}=(p_1,p_2)$, reads as
\[
\prox_{\tau G}(\widehat{\vbold}) = \argmin_{\vbold\in\RR^2} \zeta \norm{1-\vbold\cdot\pbold}_2^2 + \frac{1}{2\tau} \norm{\vbold-\widehat{\vbold}}_2^2,
\]
thus
\[
\prox_{\tau G}(\widehat{\vbold}) 
=
\begin{pmatrix}
2\zeta p_1^2 + \tau^{-1}   & 2\zeta p_1p_2 \\
2\zeta p_1p_2 & 2\zeta p_2^2 + \tau^{-1} 
\end{pmatrix}^{-1}
\begin{pmatrix}
2\zeta p_1 + \tau^{-1}\widehat{v}_1\\
2\zeta p_2 + \tau^{-1}\widehat{v}_2
\end{pmatrix}.
\]
Since $G$ is strongly-convex, we use the accelerated scheme \cref{eq: primal-dual accelerated}, with $K$ instead of $\Kcal$.

\subsubsection*{Numerical Results}\label{sec: results}
We tested \cref{alg: outer_min} in MATLAB on synthetic and real surfaces. 
Differently from \cite{LelMorSch2013}, we did not use CVX or MOSEK, making our approach suitable for larger surfaces, beyond the variable size limit imposed by CVX. 
In what follows, we will use \cref{eq: energy denoising simplified} for solving \cref{eq: surface minimization energy u} and we will test both single and joint directional regularisers, namely $\abold = (0,1,0)$, $\abold=(0,0,1)$  and $\abold = (0,\alpha_{2,0},\alpha_{3,0})$, with $\alpha_{2,0}$ and $\alpha_{3,0}$ to be chosen appropriately for the situation. 
For a better visualization of the results, a divergence RGB colormap in the range $([0.230, 0.299, 0.754],[0.706, 0.016, 0.150])$ has been applied.

\begin{figure}[!ht]
\centering
\begin{subfigure}[t]{0.194\textwidth}\centering
\includegraphics[width=\textwidth]{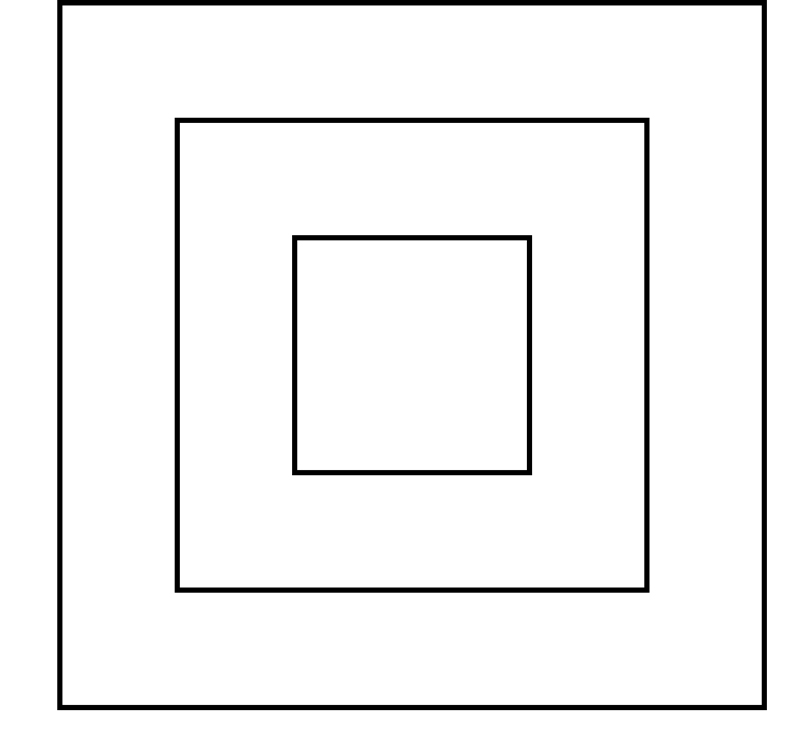}
\end{subfigure}
\hfill
\begin{subfigure}[t]{0.194\textwidth}\centering
\includegraphics[width=\textwidth]{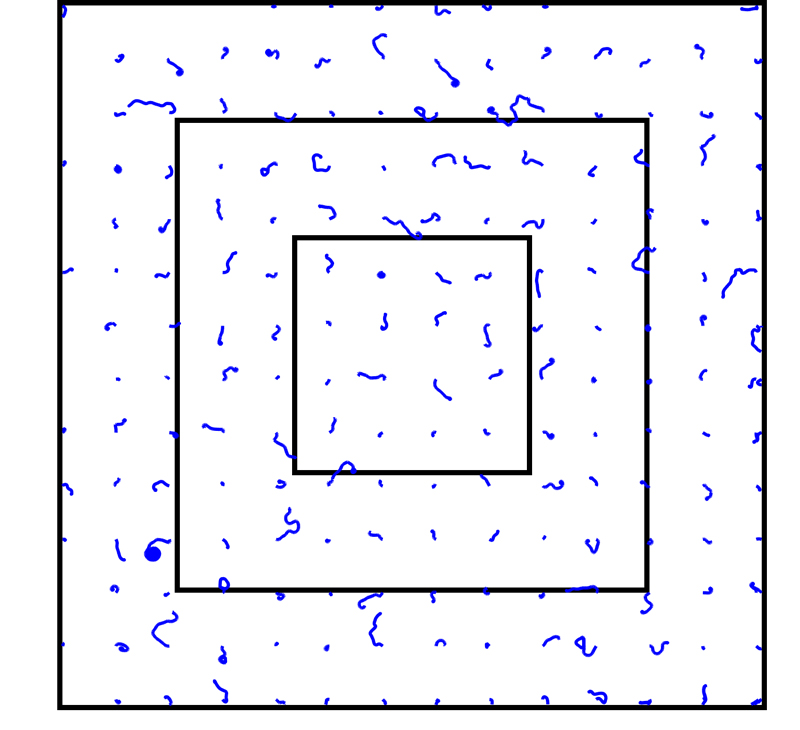}
\end{subfigure}
\hfill
\begin{subfigure}[t]{0.194\textwidth}\centering
\includegraphics[width=\textwidth]{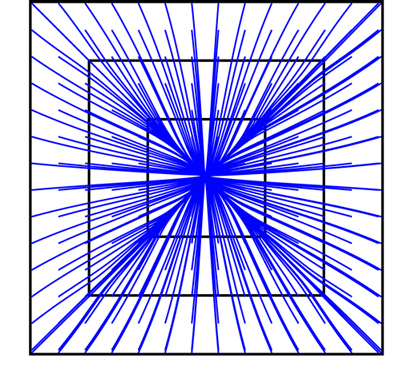}
\end{subfigure}
\hfill
\begin{subfigure}[t]{0.194\textwidth}\centering
\includegraphics[width=\textwidth]{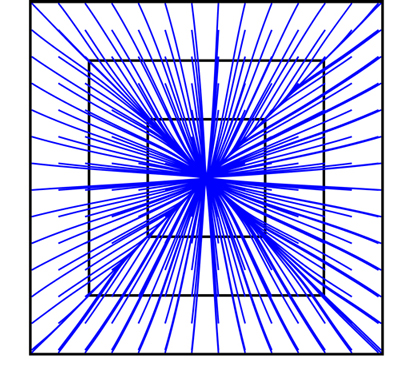}
\end{subfigure}
\hfill
\begin{subfigure}[t]{0.194\textwidth}\centering
\includegraphics[width=\textwidth]{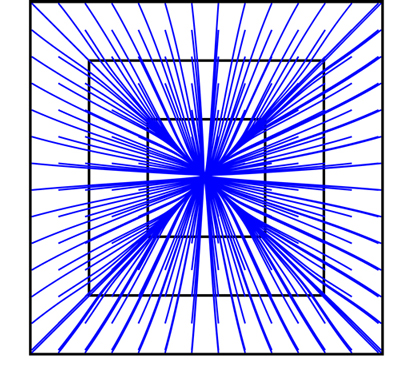}
\end{subfigure}
\\
\begin{subfigure}[t]{0.194\textwidth}\centering
\includegraphics[width=\textwidth]{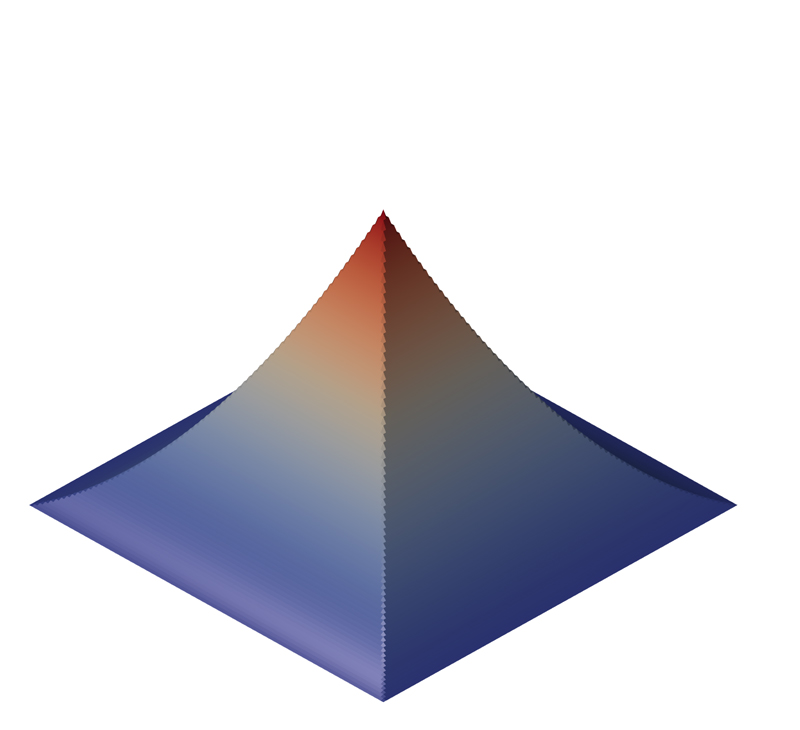}
\caption*{Ground Truth}
\end{subfigure}
\hfill
\begin{subfigure}[t]{0.194\textwidth}\centering
\includegraphics[width=\textwidth]{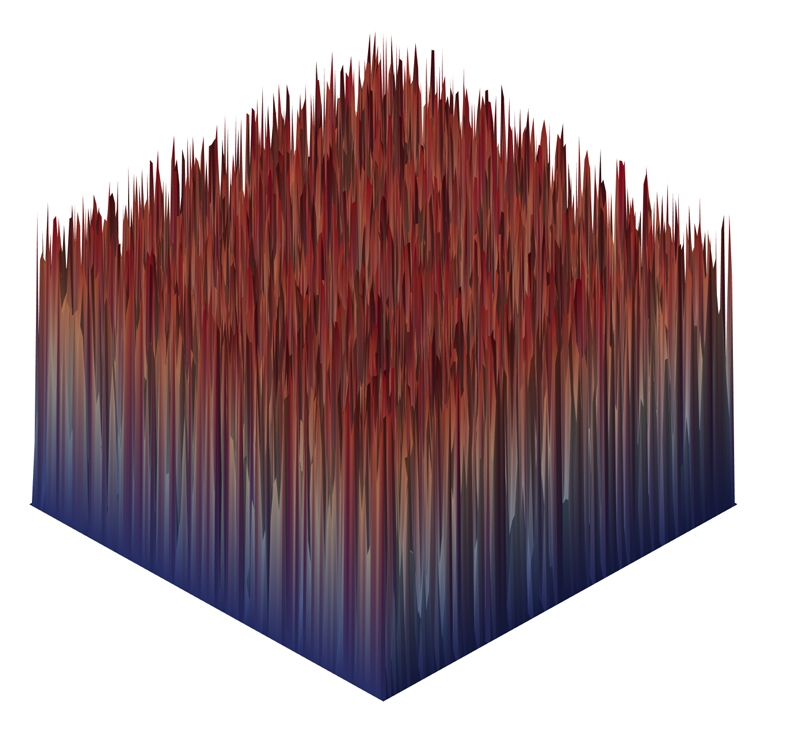}
\caption*{Input}
\end{subfigure}
\hfill
\begin{subfigure}[t]{0.194\textwidth}\centering
\includegraphics[width=\textwidth]{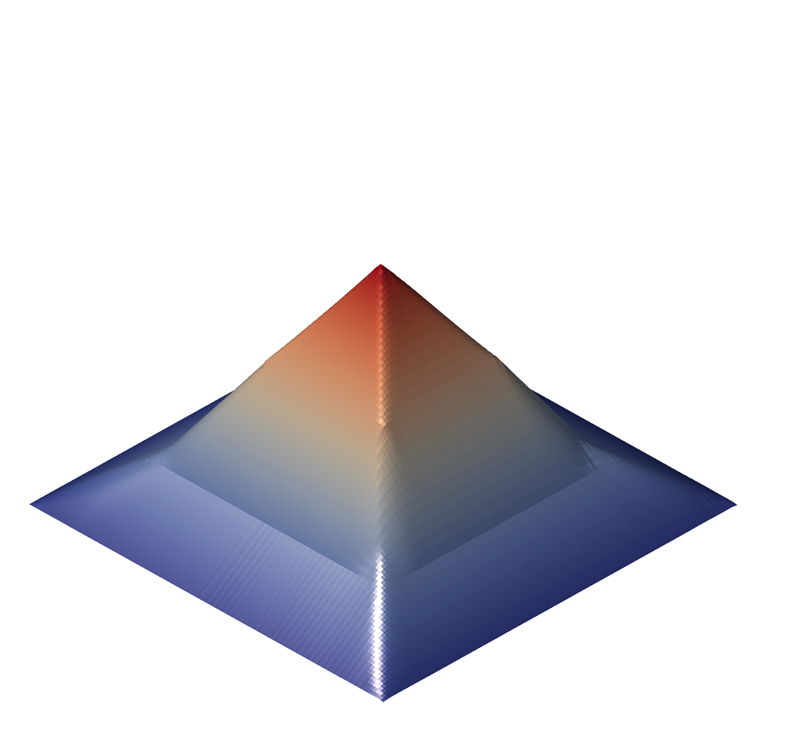}
\caption*{$\TDVM{2}{}$}
\end{subfigure}
\hfill
\begin{subfigure}[t]{0.194\textwidth}\centering
\includegraphics[width=\textwidth]{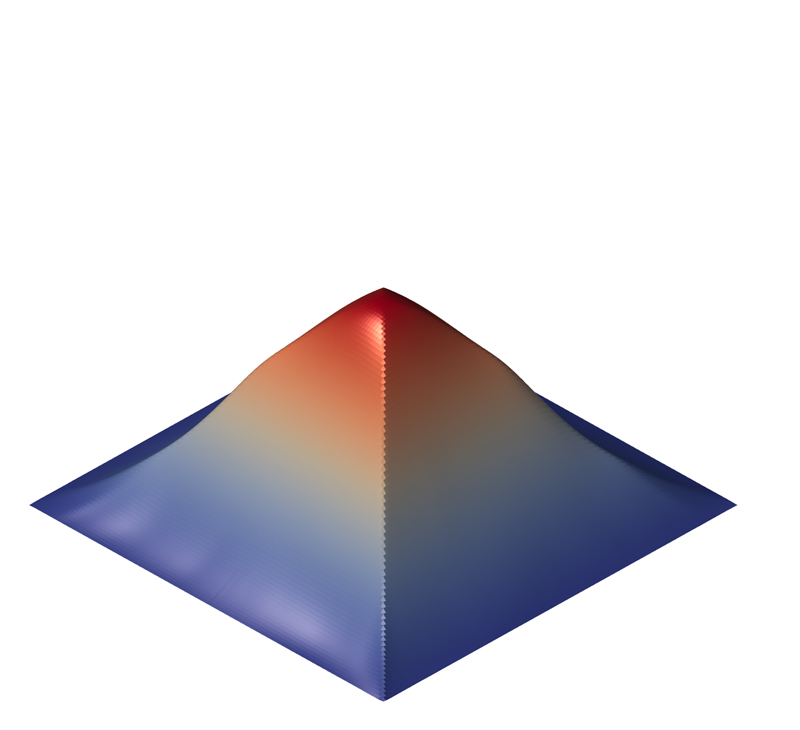}
\caption*{$\TDVM{3}{}$}
\end{subfigure}
\hfill
\begin{subfigure}[t]{0.194\textwidth}\centering
\includegraphics[width=\textwidth]{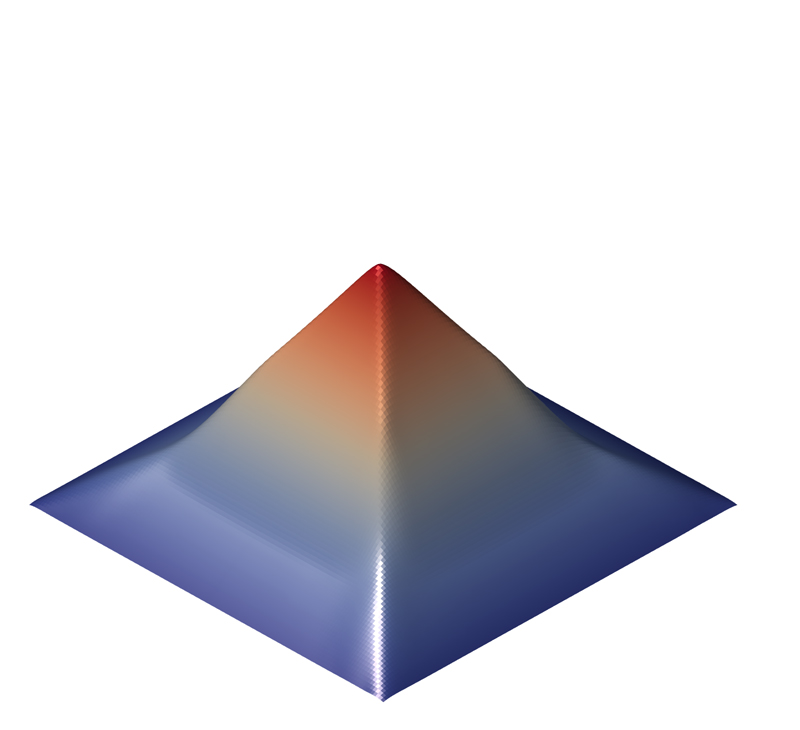}
\caption*{$\TDVM{2}{}+\TDVM{3}{}$}
\end{subfigure}
\caption{Pyramid results: 10 iterations of \cref{alg: outer_min}. Parameters: $\eta=100$, $\mu=\zeta=1$.}
\label{fig: dataset pyramid}
\end{figure}

\paragraph{Pyramid dataset from \cite{LelMorSch2013}}\label{sec: syntetic dataset}
A pyramid with height data available on three contour lines and no extra information on the tip is given, so as to test whether our model can reconstruct it.
We initialize $u^0$ and $\vbold^0$ randomly. 
In \cref{fig: dataset pyramid}, we report in the first column the location of the available data (top) and the ground truth (bottom); in the second column the random initialization of $\vbold^0$ (top) and $u^0$ (bottom); in the third, fourth and fifth columns we report the results from different orders of directional regularisers, namely $\abold=(0, 1, 0)$, $\abold=(0, 1, 0)$ and $\abold=(0, 1, 1)$, with one level of anisotropy $a=1$.
The similarity of the resulting vector fields in \cref{fig: dataset pyramid}, despite the different derivative orders involved in the minimisation with respect to $u$, shows the robustness of the computation of $\vbold$ for such problem.
Visual results suggest that a combination between second and third order directional regularisers, e.g.\ $\abold=(0,\alpha_{2,0},\alpha_{3,0})$, is desirable since it smooths the second-order result without loosing its features. 

\paragraph{SRTM dataset from \cite{SRTM}}
This dataset is part of the \emph{Shuttle Radar Topography Mission} (SRTM) \cite{SRTM} NASA mission so as to obtain elevation data for most areas of the world. 
We download \texttt{.hgt} ``height'' binary data files from \cite{readhgt}, where by selection of latitude and longitude coordinates we get 1x1 degree tiles of 1-arc seconds resolution (around $\SI{30}{\metre}$ per pixel). 
We selected some famous mountains within Italy in \cref{fig: dataset SRTM}: Etna volcano (Sicily), Baldo Mountain (Verona), Vesuvio volcano (Naples),  Brenner border (Sterzing) and Gran Sasso (L'Aquila), with image size domain of $250\times 250$ pixels. As input,  we randomly sampled about $7\%$ of sparse data on level lines and isolated points, to be interpolated with $0.1\TDVM{2}{}+\TDVM{3}{}$.
\begin{figure}[tb]
\centering
\begin{subfigure}[t]{0.19\textwidth}\centering
\includegraphics[width=\textwidth,trim={0.45cm 0cm 0.45cm 0cm},clip=true]{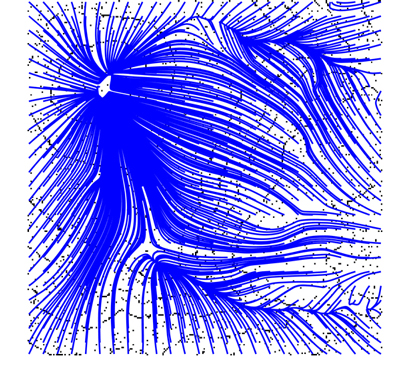}
\end{subfigure}
\hfill
\begin{subfigure}[t]{0.19\textwidth}\centering
\includegraphics[width=\textwidth,trim={0.45cm 0cm 0.45cm 0cm},clip=true]{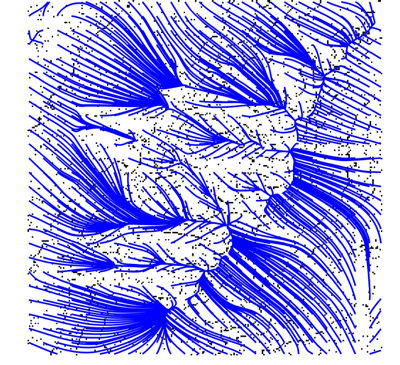}
\end{subfigure}
\hfill
\begin{subfigure}[t]{0.19\textwidth}\centering
\includegraphics[width=\textwidth,trim={0.45cm 0cm 0.45cm 0cm},clip=true]{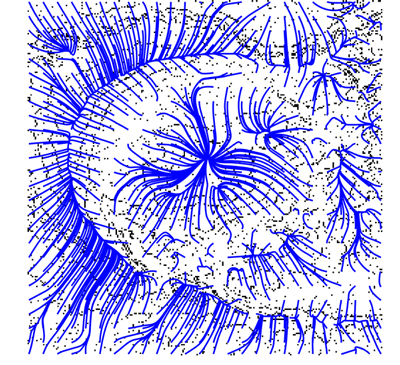}
\end{subfigure}
\hfill
\begin{subfigure}[t]{0.19\textwidth}\centering
\includegraphics[width=\textwidth,trim={0.45cm 0cm 0.45cm 0cm},clip=true]{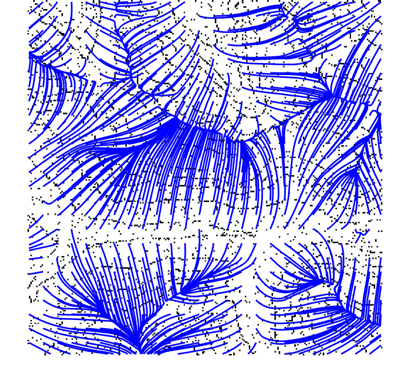}
\end{subfigure}
\hfill
\begin{subfigure}[t]{0.19\textwidth}\centering
\includegraphics[width=\textwidth,trim={0.45cm 0cm 0.45cm 0cm},clip=true]{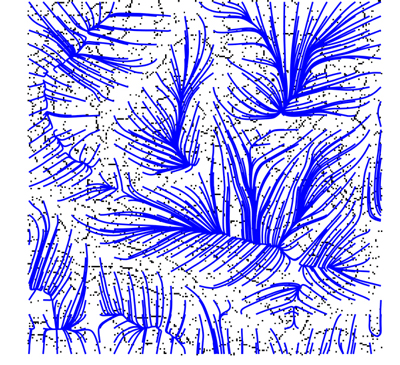}
\end{subfigure}
\\
\vspace{-0.8em}
\begin{subfigure}[t]{0.19\textwidth}\centering
\includegraphics[width=\textwidth,trim={0.45cm 0cm 0.45cm 0cm},clip=true]{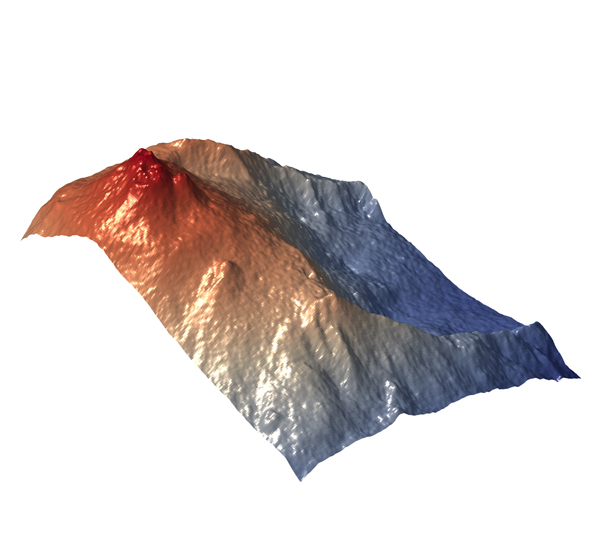}
\end{subfigure}
\hfill
\begin{subfigure}[t]{0.19\textwidth}\centering
\includegraphics[width=\textwidth,trim={0.45cm 0cm 0.45cm 0cm},clip=true]{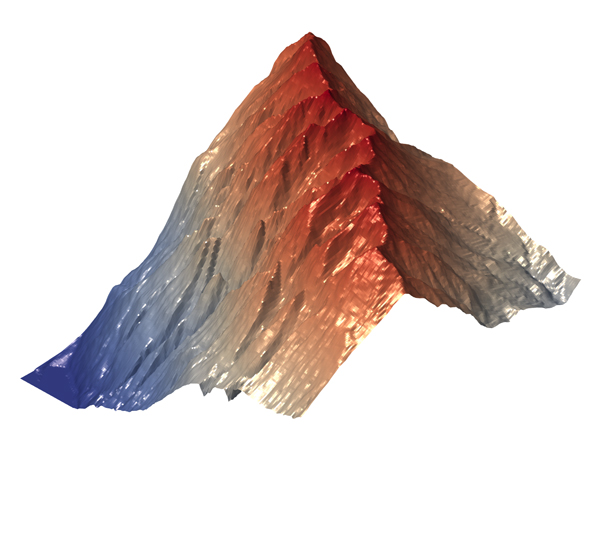}
\end{subfigure}
\hfill
\begin{subfigure}[t]{0.19\textwidth}\centering
\includegraphics[width=\textwidth,trim={0.45cm 0cm 0.45cm 0cm},clip=true]{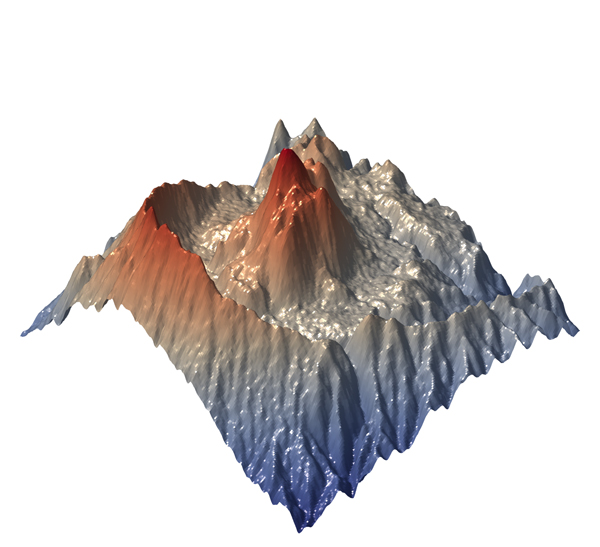}
\end{subfigure}
\hfill
\begin{subfigure}[t]{0.19\textwidth}\centering
\includegraphics[width=\textwidth,trim={0.45cm 0cm 0.45cm 0cm},clip=true]{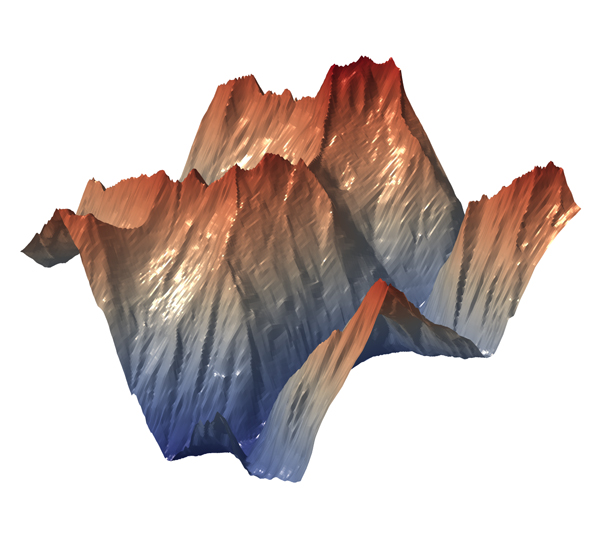}
\end{subfigure}
\hfill
\begin{subfigure}[t]{0.19\textwidth}\centering
\includegraphics[width=\textwidth,trim={0.45cm 0cm 0.45cm 0cm},clip=true]{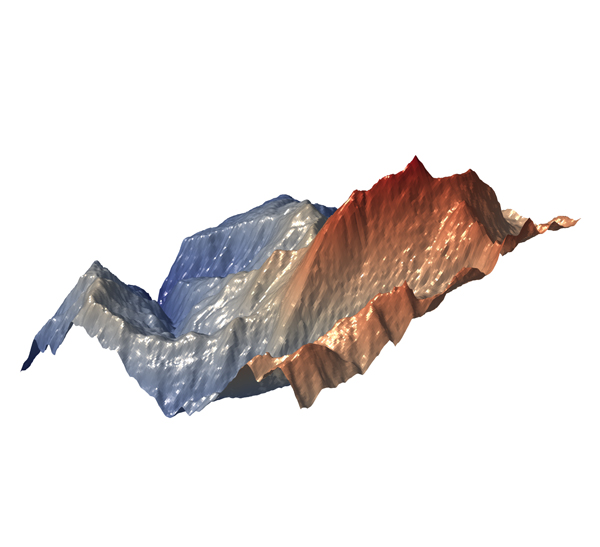}
\end{subfigure}
\\
\vspace{-0.8em}
\begin{subfigure}[t]{0.19\textwidth}\centering
\includegraphics[width=\textwidth,trim={0.45cm 0cm 0.45cm 0cm},clip=true]{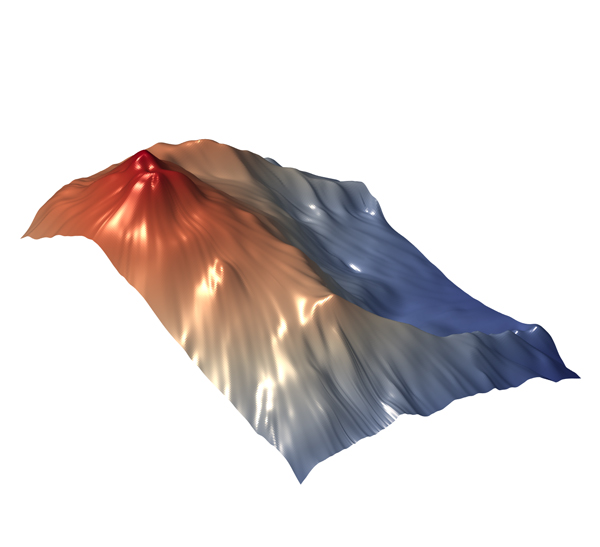}
\caption{Etna}
\label{fig: result dataset1 surf 6}
\end{subfigure}
\hfill
\begin{subfigure}[t]{0.19\textwidth}\centering
\includegraphics[width=\textwidth,trim={0.45cm 0cm 0.45cm 0cm},clip=true]{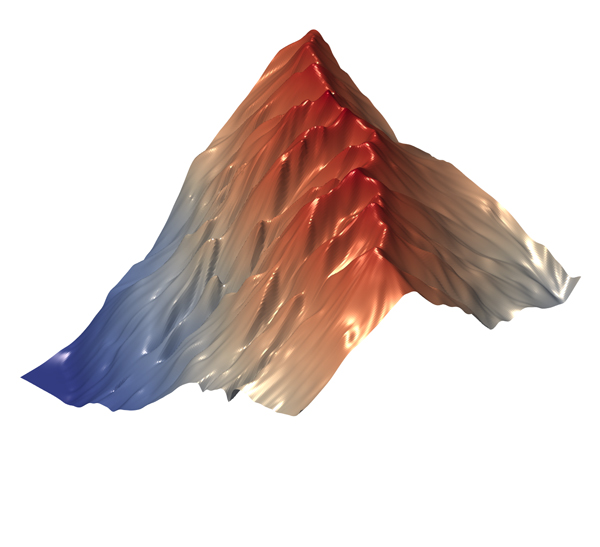}
\caption{Baldo Mountain}
\label{fig: result dataset1 surf 7}
\end{subfigure}
\hfill
\begin{subfigure}[t]{0.19\textwidth}\centering
\includegraphics[width=\textwidth,trim={0.45cm 0cm 0.45cm 0cm},clip=true]{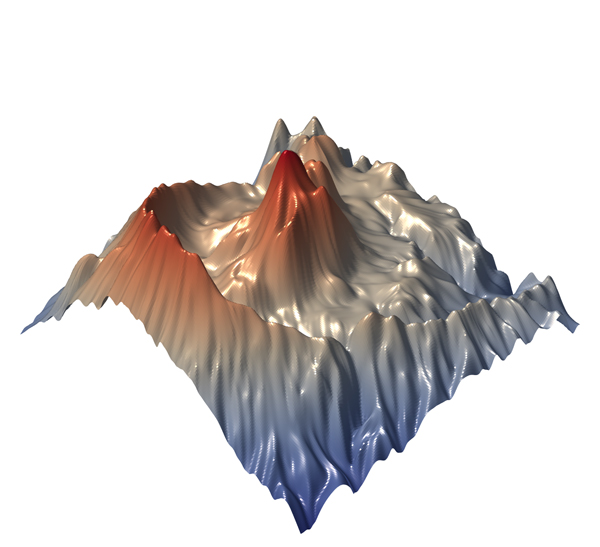}
\caption{Vesuvio}
\label{fig: result dataset1 surf 8}
\end{subfigure}
\hfill
\begin{subfigure}[t]{0.19\textwidth}\centering
\includegraphics[width=\textwidth,trim={0.45cm 0cm 0.45cm 0cm},clip=true]{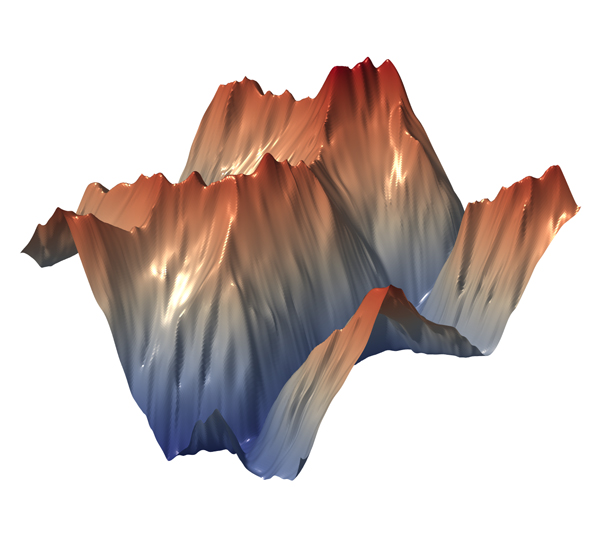}
\caption{Brenner}
\label{fig: result dataset1 surf 9}
\end{subfigure}
\hfill
\begin{subfigure}[t]{0.19\textwidth}\centering
\includegraphics[width=\textwidth,trim={0.45cm 0cm 0.45cm 0cm},clip=true]{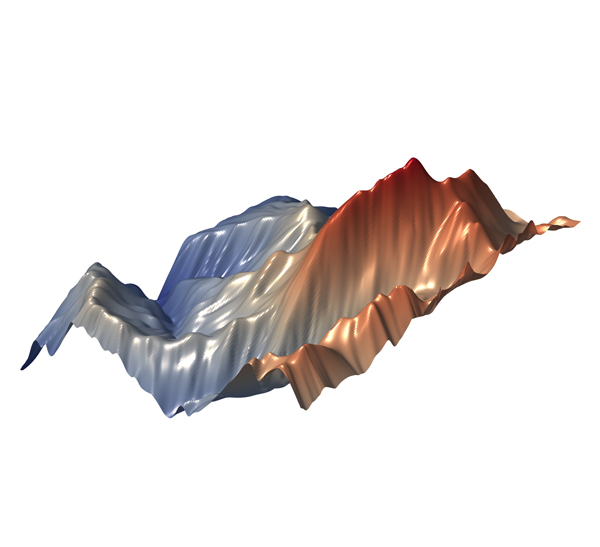}
\caption{Gran Sasso}
\label{fig: result dataset1 surf 10}
\end{subfigure}
\caption{The dataset from \cite{SRTM}.  First row: height  location (black), streamlines of $\vbold$ (blue). Second row: ground truth. Third Row: $0.1\TDVM{2}{}+\TDVM{3}{}$ results, $\mu=\zeta=1$, $\eta=1000$. }
\label{fig:  dataset SRTM}
\end{figure} 
\paragraph{Atomic Force Microscopy dataset from \cite{christian_rankl_2015_17573}}
Atomic force microscopy (AFM), or scanning probe microscopy (SPM), is a topography imaging technique commonly used in the detection of cancer cells in cellular biology: it scans objects at high resolution while recording their topographical information.
In \cite{OXVIG20171}, the study of a compressed sensing approach on AFM images was motivated by the reduction of the image acquisition time for multiple reasons, e.g.\ to minimize the operator time spent at the equipment \cite{Hansma601}, to allow time-dependent dynamic processes \cite{SCHITTER200840} and to minimize the interaction of instruments with specimens so as to reduce potential risks of damages \cite{MULLER2011461}.
Therefore, the authors proposed to speed up the sampling procedure by scanning height data on spirals rather than exploring pixel by pixel, so as to reconstruct the missing data via compressed sensing. The authors define the under-sampling ratio as
$\rho = L/L_{ref}$,
where $L$ is the length of the spiral path followed by the probe for acquiring the data and $L_{ref}$ is the distance travelled by the probe in pixels during the raster scan. 
Note that $L$ also counts the path outside the imaging domain due to smooth movement requirement of the probe, while $L_{ref}$ is approximated by the value $2\cdot\#\text{pixels}$ and the factor of two is due by the usual approach to acquire two topography buffers, even if only one is used for the visualization.
In order to test our reconstruction method based on the directional regularisers, we downloaded the open source AFM \texttt{.mi} dataset of $512\times 512$ height values from \cite{christian_rankl_2015_17573}, exported in ASCII text via the open-source software Gwyddion and imported in MATLAB. 
Our input are AFM surfaces of size $256\times 256$ obtained by slicing the orginal surface, for comparison purposes following \cite{OXVIG20171}. 
We show the results in \cref{fig: AFM} for the ground truth image in \cref{fig: afm gt}, with different under-sampling ratio $\rho$, see \cref{fig: afm 0.10,fig: afm 0.25}. 
In \cref{fig: afm comparison} we compare the structural similarity index (SSIM) \cite{SSIM} for our results (producing the black line of scores on the top of the graph) with \cite[Figure 7]{OXVIG20171}, where iterative hard thresholding (IHT), iterative soft tresholding (IST), their weighted version (w-IHT and w-IST) and Basis Pursuilt Denoising (BPDN) were tested: we conclude that our approach is robust throughout different under-sampled data, with good quality surfaces in terms of SSIM.
\begin{figure}[!htbp]
\centering
\begin{subfigure}[t]{0.245\textwidth}\centering
\captionsetup{justification=centering}
\includegraphics[width=\textwidth]{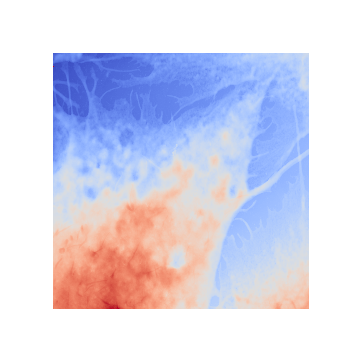}
\caption{Ground truth, \cite{christian_rankl_2015_17573}}
\label{fig: afm gt}
\end{subfigure}
\hfill
\begin{subfigure}[t]{0.735\textwidth}\centering
\centering
\includegraphics[width=0.7\textwidth]{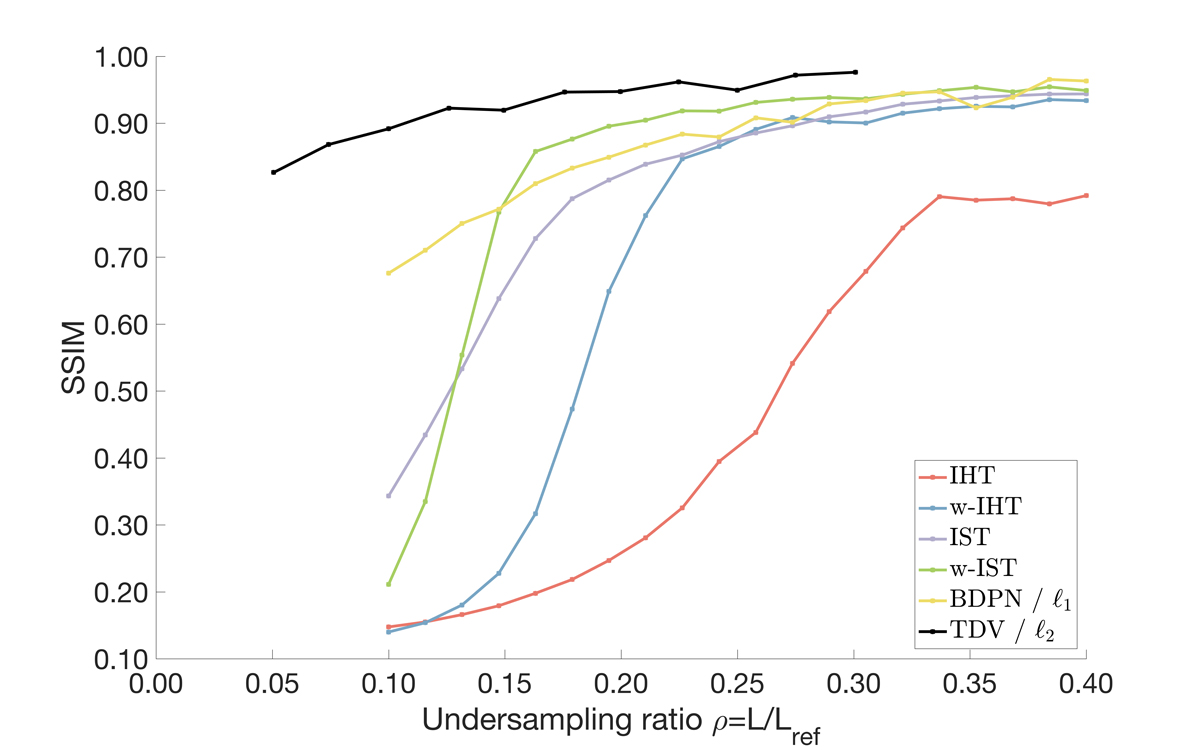}
\caption{Comparison of SSIM with \cite[Figure 7]{OXVIG20171}.}
\label{fig: afm comparison}
\end{subfigure}
\\
\begin{subfigure}[t]{0.495\textwidth}\centering
\captionsetup{justification=centering}
\settowidth{\imagewidth}{\includegraphics{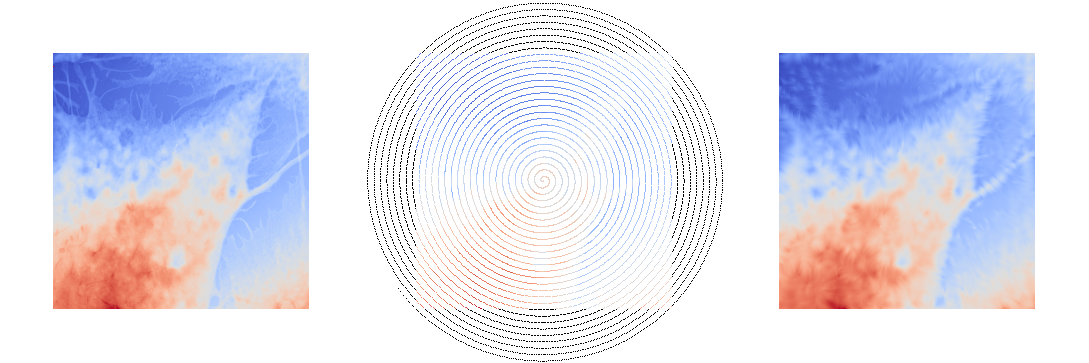}}
\includegraphics[width=\textwidth,trim = {0.333\imagewidth{}} 0 0 0,clip=true]{{./images/surface/datasetA/AFM_40_256_10}.png}
\caption{$0.1\TDVM{2}{}+\TDVM{3}{}$.\\Undersampling $\rho=0.10$.\\Input on spirals (left), result (right).\\
Parameters: $\eta=10000$, $\mu,\zeta=1$. \\
SSIM = 0.892}
\label{fig: afm 0.10}
\end{subfigure}
\hfill
\begin{subfigure}[t]{0.495\textwidth}\centering
\captionsetup{justification=centering}
\settowidth{\imagewidth}{\includegraphics{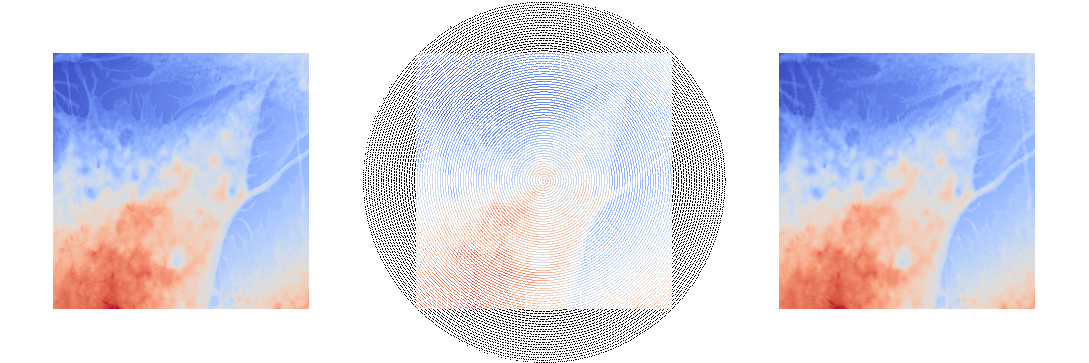}}
\includegraphics[width=\textwidth,trim = {0.333\imagewidth{}} 0 0 0,clip=true]{{./images/surface/datasetA/AFM_46_256_25}.png}
\caption{$0.1\TDVM{2}{}+\TDVM{3}{}$.\\
Undersampling $\rho=0.25$.\\
Input on spirals (left), result (right).\\
Parameters: $\eta=10000$, $\mu,\zeta=1$.\\
SSIM = 0.949}
\label{fig: afm 0.25}
\end{subfigure}
\caption{AFM reconstruction with \cref{alg: outer_min} from CC BY 4.0 data in \cite{christian_rankl_2015_17573}.}
\label{fig: AFM}
\end{figure}

\section{Conclusions}
In this work, we have shown that embedding anisotropic directional information into higher order derivatives improves the performance of total variation regularisation in many imaging applications where anisotropy plays a crucial role. In particular, we presented results for image denoising, image zooming and interpolation of scattered measurements, with details on the numerical discretisation and the solution via a primal-dual hybrid gradient algorithm. 
Among the range of experiments provided, we emphasise that our approach is particularly suitable for the reconstruction task from scattered data, motivating the interest in studying the proposed energy.
With this we provided a precise discrete framework which extends the works 
\cite{LelMorSch2013,directionaltv,BreKunPoc2010,Bredies2013}, bringing higher-order total variation together with spatially-varying anisotropy. The continuous model is analysed in the companion paper~\cite{ParMasSch18analysis}, 
\section*{Acknowledgements}
The authors are grateful to Dr.\ Martin Holler, University of Graz (Austria) for the useful discussions 
and 
to 
Prof.\ Thomas Arildsen, Aalborg University (Denmark) for the AFM data.

\medskip
\medskip
\bibliographystyle{siamplain}
\bibliography{surface}

\end{document}